\documentstyle[amssymb,amsfonts]{amsart}

\def\be#1{\begin{equation} \label{#1}}
\def\bi{\begin{itemize}}
\def\bs{\begin{split}}
\def\es{\end{split}}
\def\ba{\begin{align}}
\def\bas{\begin{align*}}
\def\ea{\end{align}}
\def\eas{\end{align*}}

\def\sgn{{\hbox{\rm sgn}}}

\def\supp{{\hbox{\rm supp}}}

\def\C{{\hbox{\bf C}}}
\def\R{{\hbox{\bf R}}}
\def\Z{{\hbox{\bf Z}}}
\def\T{{\hbox{\bf T}}}

\def\eps{\varepsilon}
\newenvironment{proof}{\noindent {\bf Proof} }{\endprf\par}
\def \endprf{\hfill  {\vrule height6pt width6pt depth0pt}\medskip}
\def\emph#1{{\it #1}}
\def\textbf#1{{\bf #1}}

\theoremstyle{plain}

  \newtheorem{proposition}[subsection]{Proposition}
  \newtheorem{lemma}[subsection]{Lemma}
  \newtheorem{corollary}[subsection]{Corollary}

\theoremstyle{remark}

\theoremstyle{definition}
  \newtheorem{definition}[subsection]{Definition}

\include{psfig}

\begin{document}

\title[Multilinear $L^2$ convolution]{Multilinear weighted convolution of $L^2$ functions, and applications to non-linear dispersive equations}
\author{Terence Tao}
\address{Department of Mathematics, UCLA, Los Angeles CA 90095-1555}
\email{tao@@math.ucla.edu}

\vspace{-0.3in}
\begin{abstract}
The $X^{s,b}$ spaces, as used by Beals, Bourgain, Kenig-Ponce-Vega, Klainerman-Machedon and others, are fundamental tools to study the low-regularity behaviour of non-linear dispersive equations.  It is of particular interest to obtain bilinear or multilinear estimates involving these spaces.  By Plancherel's theorem and duality, these estimates reduce to estimating a weighted convolution integral in terms of the $L^2$ norms of the component functions.  In this paper we systematically study weighted convolution estimates on $L^2$.  As a consequence we obtain sharp bilinear estimates for the KdV, wave, and Schr\"odinger $X^{s,b}$ spaces. 
\end{abstract}

\maketitle

\section{Introduction}

Let $Z$ be any abelian additive group with an invariant measure $d\xi$.  For instance, $Z$ could be Euclidean space $\R^{d+1}$ with Lebesgue measure, or the space $\Z^d \times \R$ with the product of counting and Lebesgue measure.  

For any integer $k \geq 2$, we let $\Gamma_k(Z)$ denote the ``hyperplane''
$$ \Gamma_k(Z) := \{ (\xi_1, \ldots, \xi_k) \in Z^k: \xi_1 + \ldots + \xi_k = 0 \}$$
with we endow with the obvious measure
$$ \int_{\Gamma_k(Z)} f := \int_{Z^{k-1}} f(\xi_1, \ldots, \xi_{k-1}, -\xi_1 - \ldots - \xi_k)\ d\xi_1\ldots d\xi_{k-1}.$$
Note that this measure is symmetric with respect to permutation of the co-ordinates.

We define a \emph{$[k;Z]$-multiplier} to be any function $m: \Gamma_k(Z) \to \C$.  If $m$ is a $[k;Z]$-multiplier, we define $\|m\|_{[k;Z]}$ to be the best constant such that the inequality
\be{czk-def}
| \int_{\Gamma_k(Z)} m(\xi) \prod_{j=1}^k f_i(\xi_i) | \leq \|m\|_{[k;Z]}
\prod_{j=1}^k \|f_i\|_{L^2(Z)}.
\end{equation}
holds for all test functions $f_i$ on $Z$.  It is clear that $\|m\|_{[k;Z]}$ determines a norm on $m$, for test functions at least; we are interested obtaining good bounds on this norm.  We will also define $\| m \|_{[k;Z]}$ in situations when $m$ is defined on all of $Z^k$ by restricting to $\Gamma_k(Z)$.

This general problem occurs frequently in the study of non-linear dispersive equations in both the periodic and non-periodic setting.  For instance, let $G$ be either $\R^d$ or $\T^d$ for some $d \geq 1$, and let $H$ be given by a real Fourier multiplier $h(\xi)$ on the dual group $Z$ (either $\R^d$ or $\Z^d$), i.e.
$$ \widehat {Hf}(\xi) := h(\xi) \hat f(\xi)$$
where the Fourier transform $\hat f$ is defined\footnote{We recommend that the reader ignore all factors of $2\pi$ which appear in the sequel.} by
$$ \hat f(\xi) := \int e^{-2\pi i x \cdot \xi} f(x)\ dx.$$
We consider non-linear Cauchy problems of the form
\be{cauchy}
\phi_t = 2\pi iH\phi + F(\phi); \quad \phi(0) = \phi_0
\end{equation}
where $\phi = \phi(x,t)$ is a field on $G \times \R$ which can either be scalar or vector-valued, the initial data $\phi_0$ lives in some Sobolev space $H^s$, and $F$ is a nonlinearity containing second-order and higher order terms.  We call the equation $\tau = h(\xi)$ the \emph{dispersion relation} of the Cauchy problem.

Examples of such problems include the modified KdV family of Cauchy problems
$$ u_t + \frac{1}{4\pi^2} u_{xxx} + u^{k-1} u_x = 0; \quad u(0) = u_0,$$
in which $G = \R$ or $\T$, $h(\xi) := \xi^3$, and non-linear Schr\"odinger Cauchy problems
$$ \phi_t + \frac{i}{2\pi}\Delta \phi + F(\phi, \overline{\phi}, \nabla \phi, \overline{\nabla \phi}) = 0; \quad u(0) = u_0$$
in which $h(\xi) := |\xi|^2$, and $F$ is some polynomial in the indicated variables.  Non-linear wave (and Klein-Gordon equations) can also be placed in this framework, by writing a second-order wave equation as a first order system and setting $h(\xi) = \pm |\xi|$.

Experience has shown that if the regularity $H^s$ of the initial data is sub-critical (i.e. if $s > s_c$, where $\dot H^{s_c}$ is the scale-invariant regularity), then the Cauchy problem \eqref{cauchy} can often be satisfactorily studied using the $X^{s,b}_{\tau = h(\xi)}(G \times \R)$ spaces\footnote{These spaces appear for the wave equation in \cite{beals:xsb}, and were applied to local well-posedness problems by \cite{borg:hsd} and (implicitly) in \cite{klainerman:nulllocal}.}, which are spaces of functions on $G \times \R$ defined via the Fourier transform as
\be{xsb-def}
\| \phi \|_{X^{s,b}_{\tau = h(\xi)}(G \times \R)} :=
\| \langle \xi \rangle^s \langle \tau - h(\xi) \rangle^b \hat \phi(\xi,\tau) \|_{L^2(G^* \times \R)}
\end{equation}
where $\langle \xi \rangle := (1 + |\xi|^2)^{1/2}$ and $G^*$ is the dual group of $G$.
For brevity we shall often abbreviate $X^{s,b}_{\tau = h(\xi)}(G \times \R)$ as
$X^{s,b}_{\tau = h(\xi)}$ or even $X^{s,b}$.

Indeed, one can usually obtain local well-posedness in \eqref{cauchy} by the method of Picard iteration provided that one can prove a multilinear estimate such as
\be{multi}
\| F(\phi) \|_{X^{s,b-1}_{\tau = h(\xi)}}
\lesssim \| \phi \|_{X^{s,b}_{\tau = h(\xi)}}^k
\end{equation}
for some $b > 1/2$, assuming that $F(\phi) = F(\phi, \ldots, \phi)$ is a $k$-linear function of $\phi$.  See e.g.
\cite{borg:hsd}, \cite{kpv:kdv}, \cite{kman.barrett} for examples of this technique, and \cite{ginibre:survey} for a general discussion.  (Of course, $F$ may be anti-linear in some of the variables, e.g. $F(\phi) = \phi \overline \phi \phi$, in which case the following discussion must be modified slightly).

It is thus of interest to obtain estimates of the form \eqref{multi}.  By the duality of the spaces $X^{s,b-1}_{\tau = h(\xi)}$ and $X^{-s,1-b}_{\tau = -h(-\xi)}$, it suffices to show the $k+1$-linear form estimate
\be{multi-form}
| \int_{G \times \R} F(\phi, \ldots, \phi) \psi | \lesssim 
\| \phi \|_{X^{s,b}_{\tau = h(\xi)}}^k
\| \psi \|_{X^{-s,1-b}_{\tau = -h(-\xi)}}.
\end{equation}
In many applications the multilinear operator $F$ is translation-invariant, and can be given by a multilinear multiplier $m$.  This means that the left-hand side of \eqref{multi-form} can be rewritten using the Fourier transform on $G \times \R$ as
$$
| \int_{\Gamma_{k+1}(G^* \times \R)} m(\xi_1, \ldots, \xi_{k+1}) 
(\prod_{j=1}^k \hat \phi(\xi_j,\tau_j)) \hat \psi(\xi_{k+1}, \tau_{k+1})|$$
where we have parameterized the co-ordinates of $\Gamma_{k+1}(G^* \times \R)$ as
$(\xi_j, \tau_j)_{j=1}^{k+1}$.  From \eqref{xsb-def}, \eqref{czk-def}, and some change of variables, we thus reduce to showing that
\be{mess}
\| 
\frac{ m(\xi) \langle \xi_{k+1} \rangle^s 
\langle \tau_{k+1} + h(-\xi_{k+1}) \rangle^{b-1} }
{ \prod_{j=1}^k \langle \xi_j \rangle^s 
\langle \tau_j - h(\xi_j) \rangle^b } \|_{[k+1; G^* \times \R]}
\end{equation}
is finite.  Thus we are led back to the problem of bounding expressions of the form $\| m \|_{[k;Z]}$.

There are two approaches to computing these quantities in the literature.  One approach proceeds using the Cauchy-Schwarz inequality, this reducing matters to integrating certain weights on intersections of hypersurfaces $\tau = h(\xi)$; the other utilizes dyadic decomposition and orthogonality before resorting to Cauchy-Schwarz.  We shall rely exclusively on the latter approach.  The advantages of dyadic decomposition are that one can re-use the estimates on dyadic blocks to prove other estimates, and the nature of interactions between different scales of frequency is more apparent.  At first glance it may appear that a dyadic decomposition may cause a logarithmic loss of exponents, but in practice this loss is either essential, or can be removed by orthogonality techniques.  A comparison between the two techniques can be found in \cite{kl-mac:null}; see also \cite{borg:book}, \cite{delort:fang}, \cite{damiano:null}, \cite{tataru:wave1}.

In this paper we systematically study the quantities $\| m \|_{[k;Z]}$, especially in the case $k=3$, which corresponds to bilinear $X^{s,b}$ estimates.  In the first half of the paper we shall discuss the $[k;Z]$ norm in a very general setting, with few assumptions on $Z$ or on the structure of $m$.  We develop some elementary but ubiquitous tools to compute these norms efficiently.  Many of these tools have appeared implicitly elsewhere in the literature, although the induction on scales techniques in Section \ref{split-sec} appear to be new.  Then, in the second half of the paper, we apply these tools to obtain sharp estimates for the dyadic components of expressions such as \eqref{mess} in the specific contexts of the KdV, Schr\"odinger and wave dispersion relations.  This reduces the verification of estimates on \eqref{mess} to that of showing that a certain explicit dyadic summation converges (although at the endpoint cases one needs some additional orthogonality arguments to eliminate logarithmic divergences).  In principle, this gives a complete characterization of all estimates on \eqref{mess} for the contexts listed above.  Unfortunately, these dyadic summations are somewhat tedious to compute and split into many cases depending on the relative sizes and signs of the frequencies $\xi_j$ and $\tau_j$.  We have tried to develop some tools (based on averaging arguments, conjugation, symmetry, etc.) to reduce the number of cases, but each such tool has a limited range of application and so can only be applied on an \emph{ad hoc} basis.

We have not attempted to produce comprehensive tables of all possible $X^{s,b}$ estimates for the KdV, wave, and Schr\"odinger equations, given that in practice one often introduces modifications to \eqref{mess} tailored to the specific application.  Instead, we focus on the estimates on dyadic blocks of \eqref{mess}, which are usable for many applications, and then present some selected applications of these estimates to prove multilinear estimates and local well-posedness results.  Some of these results have appeared before, but others seem to be new. 

In the KdV context, we derive in Section \ref{kdv-sec} sharp estimates for the dyadic blocks of \eqref{mess}, and use this to prove some bilinear and trilinear estimates of Kenig, Ponce and Vega \cite{kpv:kdv} in the periodic and non-periodic setting, as well as the $L^4$ periodic Strichartz estimate of Bourgain \cite{borg:hsd}.  More recent quadrilinear and higher estimates have been developed and applied to global well-posedness for periodic and non-periodic equations of KdV type: see \cite{ckstt}.

In the wave equation context, we derive in Section \ref{wave-sec} sharp estimates for the dyadic blocks of \eqref{mess}.  Our work here is somewhat in the spirit of \cite{damiano:null} (see also some $L^p$ variants in \cite{tv:cone2}, \cite{tao:cone}).  We then apply these estimates to prove some $X^{s,b}$ estimates relating to the Maxwell-Klein-Gordon and Yang-Mills gauge theories, generalizing some of the three-dimensional estimates of Cuccagna \cite{cuccagna}.  We remark that these estimates are 1/4 of a derivative away from reaching the critical regularity, and this seems to be due to an inherent limitation of the $X^{s,b}$ method for these equations, at least at the level of bilinear estimates.  We also indicate some methods to eliminate logarithmic divergences at the endpoint cases, although it is known that some of these divergences are essential, especially at critical regularities (see e.g. \cite{kman.barrett} for a discussion on this).

In the Schr\"odinger context, a new phenomenon arises in dimensions $d \geq 2$, namely that there is a complicated set of frequencies $(\xi_1, \xi_2, \xi_3)$ for which the denominators $\tau_j \mp |\xi_j|^2$ simultaneously vanish.  Specifically, this can occur when two of the $\xi_j$ are orthogonal.  This means that the problem of estimating \eqref{mess} accurately is akin to that of obtaining good $L^2$ bounds on a bilinear spherical Radon transform.  We have been able to obtain nearly sharp estimates in Section \ref{schro-sec} on these types of expressions by an induction-on-scales argument, which shares some intriguing similarities to some techniques in restriction theory (see e.g. \cite{borg:kakeya}, \cite{tvv:bilinear}, \cite{wolff:cone}).  The methods used here should have application to other situations (Zakharov, KP-I, KP-II, etc.) where the denominators vanishes for a complicated set of frequencies.  As a sample application we present a low-regularity local well-posedness result for a quadratic non-linear Schr\"odinger equation; this is a three-dimensional version of some results in \cite{staff:quadratic}, \cite{cdks}.

The author thanks Jim Colliander, Jean-Marc Delort, Mark Keel, Sergiu Klainerman, and Gigliola Staffilani for helpful conversations.  The author is also indebted to the referees for their excellent suggestions.  The author is supported by grants from the Packard and Sloan foundations.

\section{Notation}\label{notation-sec}

We use $A \lesssim B$ to denote the statement that $A \leq C B$ for some large constant $C$ which may vary from line to line and depend on various parameters such as the dimension $d$, and similarly use $A \ll B$ to denote the statement $A \leq C^{-1} B$.  We use $A \sim B$ to denote the statement that $A \lesssim B \lesssim A$.

Any summations over capitalized variables such as $N_j$, $L_j$, $H$ are presumed to be dyadic, i.e. these variables range over numbers of the form $2^k$ for $k \in \Z$.  We shall frequently be computing dyadic summations of positive algebraic expressions in the sequel.  To evaluate these expressions, we recommend using the heuristic that a dyadic summation is usually comparable to the largest term in the summation, which usually occurs at one of the two ends of the summation (or occasionally in an intermediate point if there is a $\min$ in the numerator, or a $\max$ or sum in the denominator).  If several terms are of comparable magnitude then one usually loses an additional logarithmic factor.

In addition to the usual notation $\chi_E$ for characteristic functions, we define $\chi_P$ for statements $P$ to be 1 if $P$ is true and 0 otherwise, e.g. $\chi_{1 \leq |\xi| \leq 2}$.  We adopt the usual convention of ignoring sets of measure zero, thus the disclaimer ``almost everywhere'' is implicit in many of our statements.

If $E \subset Z$ is a set, we use $|E|$ to denote the measure of $E$ with respect to the measure on $Z$, which may be Lebesgue measure (if $Z = \R^d$), counting measure (if $Z = \Z^d$), or some combination of the two (e.g. if $Z = \Z^d \times \R$).

Let $N_1, N_2, N_3 > 0$.  It will be convenient to define the quantities $N_{max} \geq N_{med} \geq N_{min}$ to be the maximum, median, and minimum of $N_1$, $N_2$, $N_3$ respectively.  Similarly define $L_{max} \geq L_{med} \geq L_{min}$ whenever $L_1, L_2, L_3 > 0$.  The quantities $N_j$ will measure the magnitude of frequencies of our waves, while $L_j$ measures how closely our waves approximate a free solution.  We will concentrate on the $k=3$ case, which explains why there are three $N_j$ and $L_j$.  We shall sometimes refer to the $N_j$ and $L_j$ as the $j^{th}$ \emph{frequency} and \emph{modulation} respectively.

We adopt the following summation conventions.  Any summation of the form $L_{max} \sim \ldots$ is a sum over the three dyadic variables $L_1, L_2, L_3 \gtrsim 1$, thus for instance
$$ \sum_{L_{max} \sim H} := \sum_{L_1, L_2, L_3 \gtrsim 1: L_{max} \sim H}.$$
Similarly, any summation of the form $N_{max} \sim \ldots$ sum over the three dyadic variables $N_1, N_2, N_3 > 0$, thus for instance
$$ \sum_{N_{max} \sim N_{med} \sim N} := \sum_{N_1, N_2, N_3 > 0: N_{max} \sim N_{med} \sim N}.$$

If $\tau$, $\xi$, and $h()$ are given, we adopt the convention that $\lambda$ is short-hand for
$$ \lambda := \tau - h(\xi).$$
Similarly we have
\be{lj-def}
\lambda_j := \tau_j - h_j(\xi_j).
\end{equation}
The quantity $\lambda_j$ thus measures how close in frequency the $j^{th}$ factor is to a free solution.  Generally, we shall use $N_j$ to denote the magnitude of $\xi_j$ and $L_j$ to denote the magnitude of $\lambda_j$.

\section{Basic properties}\label{elem-sec}

In this section we collect some simple properties about the operator norm $\|m\|_{[k;Z]}$.  The results in this section are implicit at various places in the literature, but we have gathered them here for explicitness.

When $k=2$, it is easy to see that the operator norm is just the $L^\infty$ norm:
$$ \|m\|_{[2;Z]} = \| m \|_{L^\infty(\Gamma_2(Z))}.$$
Thus the first non-trivial norm occurs when $k=3$.  This is the norm used to prove bilinear estimates, as indicated previously.  Henceforth we shall always assume $k > 2$.

The following comparison principle is extremely useful:

\begin{lemma}[Comparison principle]\label{comparison}
If $m$ and $M$ are $[k;Z]$ multipliers, and $|m(\xi)| \leq M(\xi)$ for all $\xi \in \Gamma_k(Z)$, then $\|m\|_{[k;Z]} \leq \| M\|_{[k;Z]}$.
Also, if $m$ is a $[k;Z]$ multiplier, and $a_1, \ldots, a_k$ are functions from $Z$ to $\R$, then
$$ \| m(\xi) \prod_{j=1}^k a_j(\xi_j) \|_{[k;Z]} \leq \| m \|_{[k;Z]} \prod_{j=1}^k \| a_j \|_\infty.$$
\end{lemma}

\begin{proof}
The first claim follows by replacing everything by absolute values in \eqref{czk-def}.  The second claim follows by applying \eqref{czk-def} with $f_j(\xi_j)$ replaced by $f_j(\xi_j) a_j(\xi_j)$.
\end{proof}

This comparison principle, combined with the triangle inequality, allows one to easily decompose the support of a multiplier into various regions upon which the analysis is easier.  For instance, one could perform dyadic decompositions of the frequencies $\xi_i$, or partition depending on the relative sizes of $\xi_i$ and $\xi_j$, etc.  This principle is also useful for controlling null forms and similar expressions.

The first version of the comparison principle ignores the possible effects of cancellation if $m$ fluctuates in sign.  As an example of cancellation, consider the expression
$$ \| p.v. \frac{1}{\xi_1 + \xi_2} \|_{[4;\R]}.$$
This expression is bounded; indeed, if one applies \eqref{czk-def} and uses the Fourier transform, the estimate reduces to
$$ C|\int\int \hat f_1(x) \hat f_2(x) \sgn(x-y) \hat f_3(y) \hat f_4(y)\ dx dy|
\lesssim \prod_{j=1}^4 \|f_j\|_2,$$
and the claim follows from replacing everything by absolute values and then applying Cauchy-Schwarz and Plancherel.  Without the cancellation, the quantity
$$ \| \frac{1}{|\xi_1 + \xi_2|} \|_{[4;\R]}$$
is infinite, as the expression in \eqref{czk-def} is not integrable even for bump functions $f_j$.  (The Bilinear Hilbert Transform estimates in \cite{lacey:thiele} can be considered as an instance of this cancellation effect in the $L^p$ setting).

Fortunately  one does not need to exploit cancellation in $m$ in many sub-critical applications.  Indeed, we will not exploit any such cancellation in this paper. 

The second version of the comparison principle can be used to imply that $\| m(\xi) a(\xi) \|_{[k;Z]} \leq C_a \|m\|_{[k;Z]}$ whenever $a$ is a sufficiently smooth bump function, by decomposing $a$ as a Fourier series in the $\xi_j$; cf. \cite{tao:boch-rest}. 

From the Comparison principle and multilinear complex interpolation we obtain a convexity theorem:

\begin{corollary}[Convexity]\label{convex}
If $m_1$, $m_2$ are non-negative $[k;Z]$ multipliers and $0 \leq \theta \leq 1$, then
$$ \| m_1^\theta m_2^{1-\theta} \|_{[k;Z]} \leq \| m_1\|_{[k;Z]}^\theta \| m_2\|_{[k;Z]}^{1-\theta}.$$
\end{corollary}

We now give three elementary propositions, whose proof we omit.

\begin{lemma}[Symmetry]\label{symmetry} The norm $\|m\|_{[k;Z]}$ is invariant under permutations of the indices $\xi_1, \ldots, \xi_k$.
\end{lemma}

\begin{lemma}[Translation invariance / Averaging]\label{translation}
For any $\xi^0 \in \Gamma_k(Z)$ and any $[k;Z]$ multiplier $m$, we have
$$ \| m(\xi) \|_{[k;Z]} = \| m(\xi + \xi_0) \|_{[k;Z]}.$$
From this and Minkowski's inequality, we thus have the averaging estimate
\be{convolve}
\| m * \mu \|_{[k;Z]} \leq \| m\|_{[k;Z]} \| \mu \|_{L^1(\Gamma_k(Z))}
\end{equation}
for any finite measure $\mu$ on $\Gamma_k(Z)$.
\end{lemma}

\begin{lemma}[Scaling]\label{scaling}
Let $L: Z \to Z$ be an automorphism on $Z$.  Suppose there exists a number $\det(L)$ such that the change of variables formula
$$ \int_Z f(L\xi)\ d\xi = |\det(L)|^{-1} \int_Z f(\xi)\ d\xi$$
holds for all test functions $f$.  Let $m$ be a $[k;Z]$ multiplier, and let $m \circ L$ denote the multiplier
$$ m \circ L(\xi_1, \ldots, \xi_k) := m(L \xi_1, \ldots, L \xi_k).$$
Then 
$$ \| m \circ L \|_{[k;Z]} = |\det(L)|^{\frac{k}{2}-1} \|m\|_{[k;Z]}.$$
\end{lemma}

In applications of scaling, $Z$ shall usually be a Euclidean space $\R^d$ and $L$ an invertible linear transformation on $Z$, in which case $\det(L)$ is just the familiar determinant.

If a multiplier splits as the tensor product arising from smaller groups $Z_1$, $Z_2$, then one can split the norm similarly:

\begin{lemma}[Direct and semi-direct tensor products]\label{tensor}
Let $Z_1, Z_2$ be abelian groups, with $Z_1 \times Z_2$ parameterized by $(\xi^1, \xi^2)$, and let $m_1$, $m_2$ be $[k;Z_1]$ and $[k;Z_2]$ multipliers respectively.  Define the tensor product $m_1 \otimes m_2$ to be the $[k;Z_1 \times Z_2]$ multiplier
$$ m_1 \otimes m_2( (\xi_1^1,\xi_1^2), \ldots, (\xi_k^1, \xi_k^2))
:= m_1(\xi_1^1, \ldots, \xi_k^1) m_2(\xi_1^2, \ldots, \xi_k^2).$$
Then we have
\be{direct}
\| m_1 \otimes m_2 \|_{[k; Z_1 \times Z_2]} = \| m_1\|_{[k;Z_1]} \| m_2\|_{[k;Z_2]}.
\end{equation}
More generally, if $m$ is a $[k;Z_1 \times Z_2]$ multiplier, define the $[k;Z_2]$ multiplier $m(\xi^1)$ for all $\xi^1 \in \Gamma_k(Z_1)$ by
$$ m(\xi^1)(\xi^2) := m( (\xi_1^1,\xi_1^2), \ldots, (\xi_k^1, \xi_k^2)).$$
Then we have
\be{semi-direct}
\| m \|_{[k;Z_1 \times Z_2]} \leq \| \| m(\xi^1) \|_{[k;Z_2]} \|_{[k;Z_1]}.
\end{equation}
\end{lemma}

\begin{proof}
The estimate \eqref{semi-direct} follows from \eqref{czk-def}, Fubini's theorem, and the identity
$$ \| f_j \|_{L^2(Z_1 \times Z_2)} = \| \| f_j(\xi^1_j)\|_{L^2(Z_2)} \|_{L^2(Z_1)}.$$  
From \eqref{semi-direct} we see that the left-hand side of \eqref{direct} is less than or equal to the right-hand side.  To prove the reverse inequality, apply \eqref{czk-def} for $m_1 \otimes m_2$ and for functions $f_j$ which split as tensor products of functions on $Z_1$ and functions on $Z_2$.
\end{proof}

We can also compose two multilinear estimates to obtain a multilinear estimate of higher order:  

\begin{lemma}[Composition and $TT^*$]\label{compose}
If $k_1, k_2 \geq 1$ and $m_1$, $m_2$ are functions on $Z^{k_1}$ and $Z^{k_2}$ respectively, then
\be{comp-est}
\begin{split}
\| m_1(\xi_1, \ldots, \xi_{k_1}) &m_2(\xi_{k_1+1}, \ldots, \xi_{k_1+k_2}) \|_{[k_1+k_2;Z]} \\
&\leq \| m_1(\xi_1, \ldots, \xi_{k_1}) \|_{[k_1+1;Z]}
\| m_2(\xi_1, \ldots, \xi_{k_2}) \|_{[k_2+1;Z]}.
\end{split}
\end{equation}
As a special case we have the $TT^*$ identity
\be{tt-est}
\| m(\xi_1, \ldots, \xi_k) \overline{m(-\xi_{k+1}, \ldots, -\xi_{2k})} \|_{[2k;Z]} = \| m(\xi_1, \ldots, \xi_k) \|_{[k+1;Z]}^2
\end{equation}
for all functions $m: Z^k \to \R$.
\end{lemma}

\begin{proof}
For $s=1,2$, let $L_s$ be the $k_s$-linear operator defined by
$$ \langle L_s(f_1, \ldots, f_{k_s}), f_{k_s+1} \rangle = 
\int_{\Gamma_{k_s+1}(Z)} m_1(\xi) \prod_{j=1}^{k_s+1} f_j(\xi_j).$$
From duality we have
$$ \| L_s(f_1, \ldots, f_{k_s}) \|_2
\leq \| m_s \|_{[k_s+1;Z]} \prod_{j=1}^{k_s} \|f_j\|_2.$$
By Cauchy-Schwarz we thus have
$$
\int_Z L_1(f_1, \ldots, f_{k_1}) L_2(f_{k_1+1}, \ldots, f_{k_1+k_2})
\leq \| m_1 \|_{[k_1+1;Z]} \| m_2 \|_{[k_2+1;Z]}
 \prod_{j=1}^{k_1+k_2} \|f_j\|_2.$$
This gives \eqref{comp-est}.  This (together with Lemma \ref{scaling}) implies that the left-hand side of \eqref{tt-est} is less than or equal to the right-hand side.  To prove the reverse inequality, apply \eqref{czk-def} with $f_{k+j}(\xi) := \overline{f_j(-\xi)}$ for $1 \leq j \leq k$ and use duality.
\end{proof}

As an immediate consequence of \eqref{tt-est} and symmetry (Lemma \ref{symmetry}) we have

\begin{corollary}[Conjugation]\label{conjugate}
If $m_1: Z^{k_1} \to \R$ and $m_2: Z^{k_2} \to \R$ are functions for some $k_1, k_2 \geq 1$, then
\bas \| m_1(\xi_1, \ldots, \xi_{k_1}) &m_2(\xi_{k_1+1}, \ldots, \xi_{k_1+k_2}) \|_{[k_1+k_2+1;Z]} \\ &= 
\| \overline{m_1(-\xi_1, \ldots, -\xi_{k_1})} m_2(\xi_{k_1+1}, \ldots, \xi_{k_1+k_2}) \|_{[k_1+k_2+1;Z]}.
\end{align*}
\end{corollary}

This Corollary can be viewed as a restatement of the trivial identity $\| u v\|_2 = \| \overline{u} v \|_2$, and can be quite effective when combined with the Cauchy-Schwartz inequality (see below).

For any $1 \leq j \leq k$ and $\eta_j \in Z$, let $\Gamma_k(Z; \xi_j = \eta_j)$ denote the set
$$ \Gamma_k(Z; \xi_j = \eta_j) := \{ (\xi_1, \ldots, \xi_k) \in \Gamma_k(Z): \xi_j = \eta_j) \}.$$
In other words, $\Gamma_k(Z; \xi_j = \eta_j)$ is the section of $\Gamma_k(Z)$ when the $\xi_j$ variable is frozen at $\eta_j$.
We endow this space with the induced measure from $\Gamma_k(Z)$, thus
$$ \int_{\Gamma_k(Z;\xi_k = \eta_k)} f := \int_{Z^{k-2}} f(\xi_1, \ldots, \xi_{k-2}, -\xi_1 - \ldots - \xi_{k-2} - \eta_k, \eta_k)\ d\xi_1\ldots d\xi_{k-2}$$
and similarly for other values of $j$.

The following application of the Cauchy-Schwarz inequality has been used implicitly in many places.

\begin{lemma}[Cauchy-Schwarz estimate]\label{cauchy-schwarz}
If $m$ is a $[k;Z]$ multiplier and $1 \leq j \leq k$, then
$$ \| m \|_{k;Z} \leq \sup_{\eta_j \in Z} (\int_{\Gamma_k(Z; \xi_j = \eta_j)}
|m(\xi)|^2)^{1/2}.$$
\end{lemma}

\begin{proof}
By Lemma \ref{symmetry} we may take $j=k$.  We can rewrite the left-hand side of \eqref{czk-def} as
$$ \int_Z f_k(\eta_k) (\int_{\Gamma_k(Z; \xi_k = \eta_k)}
m(\xi) \prod_{j=1}^{k-1} f_j(\xi_j))\ d\eta_k.$$
From Fubini's theorem we have
$$ \int_Z (\int_{\Gamma_k(Z; \xi_k = \eta_k)}
\prod_{j=1}^{k-1} |f_j(\xi_j)|^2)\ d\eta_k = \prod_{j=1}^{k-1} \|f_j\|_2^2.$$
By Cauchy-Schwarz, the left-hand side of \eqref{czk-def} is thus less than or equal to
$$  \prod_{j=1}^{k-1} \|f_j\|_2
(\int_Z |f_k(\eta_k)|^2
(\int_{\Gamma_k(Z; \xi_j = \eta_j)}
|m(\xi)|^2)\ d\eta_k)^{1/2}.$$
The claim follows.
\end{proof}

In some cases this Lemma is essentially sharp.  For instance, we have

\begin{corollary}\label{cz-cor}
For any complex functions $m_1(\xi)$, $m_2(\xi)$ on $\Z$ we have 
\be{m-bold}
\frac{ \| |m_1|^2 * |m_2|^2 \|_2 }{ \| |m_1|^2 * |m_2|^2 \|_1^{1/2}} 
\leq \| m_1(\xi_1) m_2(\xi_2) \|_{[3;Z]} \leq
\| |m_1|^2 * |m_2|^2 \|_\infty^{1/2}.
\end{equation}
In particular, for any subsets $A,B$ of $Z$ we have
the characteristic function estimate
\be{char-zar}
\| \chi_A(\xi_1) \chi_B(\xi_2) \|_{[3;Z]} \leq
| \{ \xi_1 \in A: \xi - \xi_1 \in B \} |^{1/2}
\end{equation}
for some $\xi \in Z$.
\end{corollary}

\begin{proof}
The right-hand side of \eqref{m-bold} (and thus \eqref{char-zar}) follows immediately from Lemma \ref{cauchy-schwarz}.  The left-hand side follows from \eqref{czk-def} and setting $f_1 = \overline{m_1}$, $f_2 = \overline{m_2}$, $f_3 = |m_1|^2 * |m_2|^2$.
\end{proof}

\begin{figure}[htbp] \centering
\ \psfig{figure=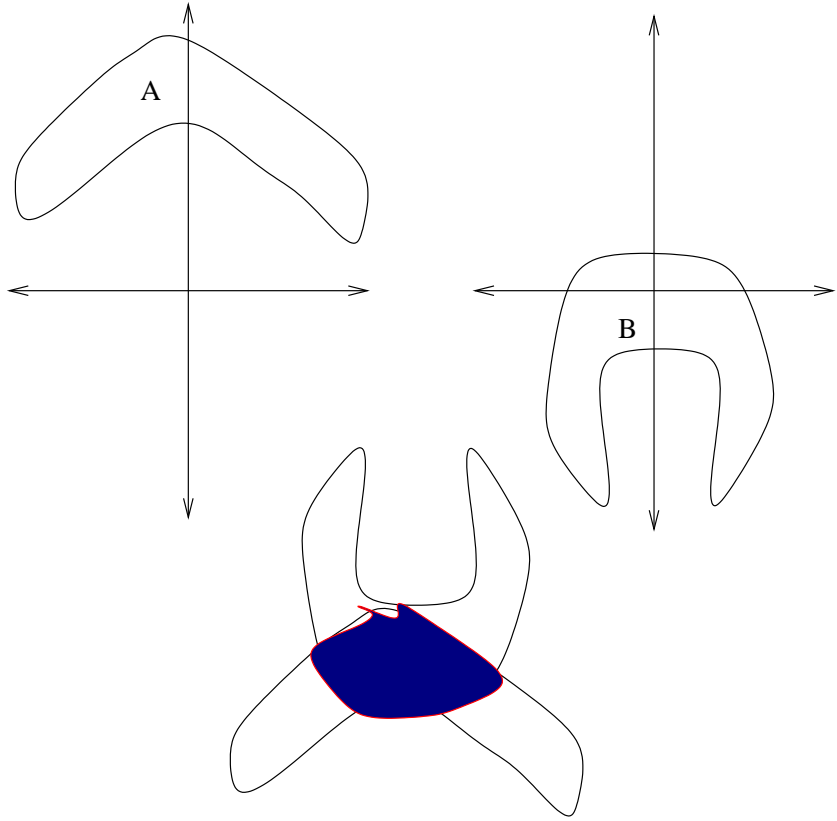}
\caption{A depiction of \eqref{char-zar}.  When $\xi_1$ is restricted to $A$ and $\xi_2$ is restricted to $B$, the multiplier norm in \eqref{char-zar} is controlled by the square root of the largest intersection of $A$ and some inverted translate of $B$.
 }
\end{figure}

This lemma gives good control on $\| \chi_A(\xi_1) \chi_B(\xi_2)\|_{[3;Z]}$ when $\chi_A * \chi_B$ resembles a constant multiple of a characteristic function (i.e. no sharp ``spikes'' which would disrupt the $L^\infty$ norm).  Occasionally these spikes can be eliminated by applying Corollary \ref{conjugate} to flip $A$ or $B$; cf. the ``doubling arguments'' in e.g. \cite{damiano:null}.  We shall be able to show that \eqref{char-zar} is sharp when one of $A$, $B$ is a box; see Corollary \ref{boxes}.

In applications the sets $A$, $B$ in \eqref{char-zar} will often be neighbourhoods of hypersurfaces.  To compute the right-hand side one therefore has to consider how one of these hypersurfaces intersects (an inverted translate of) the other hypersurface.  Thus one expects \eqref{char-zar} to be efficient when these two hypersurfaces are transverse; we shall make this idea rigorous in Section \ref{transverse-sec}.

Although the Cauchy-Schwarz estimate (Lemma \ref{cauchy}) is efficient for many simple situations (especially when the multiplier $m$ has a tensor product structure, e.g. $m = m_1(\xi_1) m_2(\xi_2)$), it does not give sharp results in all cases, and one must often perform some additional decompositions as well as orthogonality arguments.  Our main orthogonality tool which will be the following application of Schur's test.

If $m$ is a $[k;Z]$ multiplier and $1 \leq j \leq k$, we define the \emph{$j$-support} $\supp_j(m) \subset Z$ of $m$ to be the set
$$ \supp_j(m) := \{ \eta_j \in Z: \Gamma_k(Z; \xi_j = \eta_j) \cap \supp(m) \neq \emptyset \}.$$ 
More generally, if $J$ is a non-empty subset of $\{1,\ldots, k\}$, we define the set $\supp_J(m) \subset Z^J$ by
$$ \supp_J(m) := \prod_{j \in J} \supp_j(m).$$
Note that $\supp_{\{1,\ldots,k\}}(m)$ can be much larger than $\supp(m)$.

\begin{lemma}[Schur's test]\label{schur}
Let $J_1, J_2$ be disjoint non-empty subsets of $\{1, \ldots, k\}$ and $A_1,A_2 > 0$.  Suppose that $(m_\alpha)_{\alpha \in I}$ is a collection of $[k;Z]$ multipliers such that 
$$ \# \{ \alpha \in I: \xi \in \supp_{J_s}(m_\alpha) \} \leq A_s$$
for all $\xi \in Z^{J_s}$ and $s=1,2$.  Then
$$ \| \sum_{\alpha \in I} m_\alpha \|_{[k;Z]} \leq (A_1 A_2)^{1/2} 
\sup_{\alpha \in I} \| m_\alpha \|_{[k;Z]}.$$
In particular, if the $m_\alpha$ are non-negative and $A_1, A_2 \sim 1$, then we have the orthogonality estimate
\be{sharp-schur}
\| \sum_{\alpha \in I} m_\alpha \|_{[k;Z]} \sim 
\sup_{\alpha \in I} \| m_\alpha \|_{[k;Z]}.
\end{equation}
\end{lemma}

\begin{proof}
By adding dummy elements to $J_1$, $J_2$ if necessary we may assume that $J_1 \cup J_2 = \{1,\ldots, k\}$.

Consider the quantity
$$ |\int_{\Gamma_k(Z)} \sum_{\alpha \in I} m_\alpha(\xi) \prod_{j=1}^k f_j(\xi_j)|.$$
By \eqref{czk-def} applied to each summand, this is less than or equal to
$$ \sum_{\alpha \in I} \| m_\alpha \|_{[k;Z]}
\prod_{s=1}^2 \| f_{J_s} \|_{L^2(\supp_{J_s}(m_\alpha))}$$
where $f_{J_s}$ is the tensor product of all the $f_j$ for $j \in J_s$.  By \eqref{czk-def} it thus suffices to show that
$$
\sum_{\alpha \in I} 
\prod_{s=1}^2 \| f_{J_s} \|_{L^2(\supp_{J_s}(m_\alpha))}
\leq (A_1 A_2)^{1/2} \| f_{J_1}\|_{L^2(Z^{J_1})} \|f_{J_2}\|_{L^2(Z^{J_2})}.
$$
From the overlap of the $\supp_{J_1}(m_\alpha)$ we have
$$
\sum_{\alpha \in I} 
\| f_{J_1} \|_{L^1(\supp_{J_1}(m_\alpha))}
\| f_{J_2} \|_{L^\infty(\supp_{J_1}(m_\alpha))}
\leq A_1 \| f_{J_1}\|_{L^1(Z^{J_1})} \|f_{J_2}\|_{L^\infty(Z^{J_2})}.$$
Similarly with the roles of 1 and 2 reversed.  The former claim then follows by interpolation.  The latter claim then follows from the former claim and the Comparison principle.
\end{proof}

Usually we shall apply this lemma when $J_1 = \{j_1\}$, $J_2 = \{j_2\}$ are singleton sets.

Informally, \eqref{sharp-schur} states that if $m$ is a multiplier and the $J_1$ and $J_2$ frequency spaces can be divided into regions such that each region in $J_1$ interacts (via $m$) with only finitely many regions in $J_2$ and conversely, then the norm of $m$ is essentially equivalent to the norm of $m$ restricted to a pair of interacting regions.  Generally, the regions of frequency space to use will be dictated\footnote{Informally, one can locate the regions to use by choosing $1 \leq i,j \leq k$ and setting up an equivalence relation $\xi_i \sim \xi_j$ whenever $m(\xi)$ is non-zero. Any pair of $\xi_i$ which are connected by a chain of these relations, whose length is even and $O(1)$, should then belong to the same region, and similarly for the $\xi_j$.}
 by the geometry of the support of $m$.  For instance, if $k=3$ and $m(\xi_1,\xi_2,\xi_3)$ is supported on the region $\xi_1 + \xi_3 = O(R)$, then one should decompose the variables $\xi_1$ and $\xi_3$ into balls of radius $R$.  If instead the multiplier is supported on the region where $\angle(\xi_1,\xi_2) = O(\theta)$, where $\angle(\xi_1, \xi_2)$ is the angle between $\xi_1$ and $\xi_2$, then one should decompose the variables $\xi_1$ and $\xi_2$ into sectors of angular width $\theta$.  

In typical applications of Lemma \ref{schur}, the regions of frequency space will be ``boxes'', which we now pause to define.

\begin{definition}\label{box-def}
A \emph{box covering} of $Z$ is a partitioning of $Z$ into disjoint sets 
$(R + \eta)_{\eta \in \Sigma}$, where the \emph{fundamental domain} $R$ is a subset of $Z$ with non-zero (possibly infinite) measure which is symmetric around (and contains) the origin, and the
\emph{tiling lattice} $\Sigma$ is a discrete subgroup of $\Z$ such that the set $R+R$ can be covered by $O(1)$ boxes in the box covering.  We refer to the sets $R+\eta$ in the box covering as \emph{boxes}.
\end{definition}

A typical example of a box covering is when $Z := \R^n$, $R$ is the unit cube $R := [-1/2,1/2]^n$, and $\Sigma$ is just the integer lattice $\Sigma := \Z^n$.  Another example is when $Z := \Z^2$, $R := \Z \times \{0\}$, and $\Sigma := \{0\} \times \Z$.

Note that an induction argument shows that the sets $(R+\ldots+R+\eta)_{\eta \in \Sigma}$ have an overlap of $O(C^k)$ and have volume $O(C^k |R|)$, where there are $k$ copies of $R$ in the summation.

\begin{corollary}[Box localization]\label{boxes}
Suppose $(R + \eta)_{\eta \in \Sigma}$ is a box covering of $\Z$, and $m$ is a $[k;Z]$ multiplier such that each $\supp_j(m)$ is contained in a box in this covering for all $1 \leq j \leq k-2$.  Then 
\be{box-est}
\| m \|_{[k;Z]} \sim \sup_{\eta_{k-1}, \eta_k \in \Sigma}
\| m(\xi) \chi_{R+\eta_{k-1}}(\xi_{k-1}) \chi_{R + \eta_k}(\xi_k)
\|_{[k;Z]}.
\end{equation}
The implicit constants may depend on $k$.  Similar statements hold if we permute the indices $1, \ldots, k$.
\end{corollary}

In other words, if all but two of the $j$-supports of a multiplier are restricted to a box, then we can also restrict the other two $j$-supports to a similar box.

\begin{proof}
The lower bound for \eqref{box-est} follows from Lemma \ref{comparison}, so it suffices to show the upper bound.

Write $m = \sum_{\eta_{k-1}, \eta_k \in \Sigma} m_{\eta_{k-1}, \eta_k}$, where
$$ m_{\eta_{k-1},\eta_k} = m(\xi) \chi_{R+\eta_{k-1}}(\xi_{k-1}) \chi_{R + \eta_k}(\xi_k).$$
From the support properties of $m$ we see that for fixed $\eta_{k-1}$ there are at most $O(C^k)$ values of $\eta_k$ for which $m_{\eta_{k-1},\eta_k}$ does not vanish, and similarly with the roles of $\eta_{k-1}$ and $\eta_k$ reversed.  Restricting $(\eta_{k-1},\eta_k)$ to those pairs for which $m_{\eta_{k-1},\eta_k}$ does not vanish and applying Lemma \ref{schur}, the claim follows.
\end{proof}

As a consequence of the above theory we can obtain sharp bounds on $\| m \|_{[3;Z]}$ if $m$ has a sufficiently simple tensor product structure.  More precisely, we have

\begin{lemma}[Tensored box lemma]\label{box}
Suppose $(R + \eta)_{\eta \in \Sigma}$ is a box covering of $Z$, and $m(\xi)$ is a function from $Z$ to $\R$.  Then for any $\eta \in \Sigma$, we have
\be{box-est-2}
\| m(\xi_1) \chi_{R + \eta}(\xi_2) \|_{[3;Z]} \sim 
\sup_{\eta' \in \Sigma} \| m\|_{L^2(R + \eta')}.
\end{equation}
Also, we have
\be{box-max}
\| m(\xi_1)\|_{[3;Z]} = \|m\|_2.
\end{equation}
\end{lemma}

\begin{proof}  By a limiting argument we may assume that $|R|$ is finite.
From \eqref{box-est} and Lemma \ref{comparison} we have
$$ \| m(\xi_1) \chi_{R + \eta}(\xi_2) \|_{[3;Z]} \lesssim
\sup_{\eta' \in \Sigma} \| m(\xi_1) \chi_{R + \eta'}(\xi_1) \|_{[3;Z]}.$$
Applying Lemma \ref{cauchy-schwarz} we obtain the $\lesssim$ side of \eqref{box-est-2}.  To obtain the $\gtrsim$ side, apply \eqref{czk-def} with $f_1(\xi) = \overline{m(\xi)} \chi_{R + \eta'}(\xi)$, $f_2 = \chi_{R + \eta}$, and $f_3 = \chi_{R + R - \eta - \eta'}$.  Finally, \eqref{box-max} comes from Lemma \ref{cauchy-schwarz} and testing \eqref{czk-def} with $f_1 = \overline m$ and $f_2 = f_3$ being large characteristic functions.
\end{proof}

We shall develop some more specialized tools of the above type in later sections.  For now, we give a simple (and well-known) application of the above theory which already illustrates many of the techniques we shall use to tackle expressions such as \eqref{mess}:

\begin{proposition}\label{sobolev}  Let $d \geq 1$, and let $s_1, s_2, s_3$ be such that
\be{sum}
\min(s_1+s_2, s_2+s_3, s_3+s_1) \geq 0, \quad s_1 + s_2 + s_3 \geq \frac{d}{2},
\end{equation}
with at least one of the above two inequalities being strict.  Then
\be{inhomog}
\| \frac{1}{\langle \xi_1 \rangle^{s_1} \langle \xi_2 \rangle^{s_2} \langle \xi_3 \rangle^{s_3}} \|_{[3, \R^d]} \leq C_{d,s_1,s_2,s_3}.
\end{equation}
If we specialize to the case
$$ \min(s_1+s_2, s_2+s_3, s_3+s_1) > 0, \quad s_1 + s_2 + s_3 = \frac{d}{2}$$
then we have the homogeneous version
$$ \| \frac{1}{|\xi_1|^{s_1} |\xi_2|^{s_2} |\xi_3|^{s_3}} \|_{[3, \R^d]} \leq C_{d,s_1,s_2,s_3}.$$
\end{proposition}

\begin{proof}
We shall just prove the inhomogeneous version \eqref{inhomog}; the homogeneous version easily follows from \eqref{inhomog} and a scaling and limiting argument based on Lemma \ref{scaling}.

By Lemma \ref{comparison} and symmetry we may restrict to the region $|\xi_1| \leq |\xi_2| \leq |\xi_3|$.  Since $\xi_3 = -\xi_1 - \xi_2$, we thus have $|\xi_3| \sim |\xi_2|$.  We now dyadically decompose the left-hand side of \eqref{inhomog} as
$$ \| \sum_{N_1,N_3: 0 < N_1 \lesssim N_3} \frac{ \chi_{|\xi_1| \sim N_1} \chi_{|\xi_2| \sim N_3} \chi_{|\xi_3| \sim N_3} } {\langle N_1 \rangle^{s_1} \langle N_3 \rangle^{s_2 + s_3} } \|_{[3,\R^d]}.$$
The multiplier inside the summation has essentially disjoint $\xi_2$ and $\xi_3$ supports as $N_3$ varies dyadically.  By Lemma \ref{schur}, we can thus estimate the above by
$$\sup_{N_3} \| \sum_{N_1: 0 < N_1 \lesssim N_3} \frac{ \chi_{|\xi_1| \sim N_1} \chi_{|\xi_2| \sim N_3} \chi_{|\xi_3| \sim N_3} } {\langle N_1 \rangle^{s_1} \langle N_3 \rangle^{s_2 + s_3} } \|_{[3,\R^d]}.$$
By the triangle inequality and Comparison principle we may estimate this by
$$\sup_{N_3} \sum_{N_1: 0 < N_1 \lesssim N_3} \| \chi_{|\xi_1| \sim N_1} \|_{[3,\R^d]}
\langle N_1 \rangle^{-s_1} \langle N_3 \rangle^{-s_2 + s_3} .$$
By the Tensored Box lemma we have
$$ \| \chi_{|\xi_1| \sim N_1} \|_{[3,\R^d]} \sim N_1^{d/2}.$$
The claim then follows from the assumptions \eqref{sum} and a simple computation.
\end{proof}

By duality we thus have the well-known

\begin{corollary}[Sobolev multiplication law]\label{sobolev-mult}  Let $d \geq 1$, and let $s_1, s_2, s$ be such that
$$ s_1 + s_2 \geq 0, \quad s \leq s_1,s_2, \quad s < s_1 + s_2 - \frac{d}{2}$$
or
$$ s_1 + s_2 > 0, \quad s < s_1,s_2, \quad s \leq s_1 + s_2 - \frac{d}{2}.$$
Then
$$ \| \phi\psi \|_{H^s(\R^d)} \lesssim \| \phi \|_{H^{s_1}(\R^n)} \| \psi \|_{H^{s_2}(\R^d)}.$$
If we specialize to the case
$$ s_1 + s_2 > 0, \quad s < s_1,s_2, \quad s = s_1 + s_2 - \frac{d}{2}$$
then we have the homogeneous version
$$ \| \phi\psi \|_{\dot H^s(\R^d)} \lesssim \| \phi \|_{\dot H^{s_1}(\R^n)} \| \psi \|_{\dot H^{s_2}(\R^d)}.$$
The implicit constants depend on $d$, $s_1$, $s_2$, $s_3$.
\end{corollary}

We remark that the condition $s_1 + s_2 \geq 0$ is necessary in order for $\phi \psi$ to make sense even as a distribution.  The condition $s \leq s_1,s_2$ reflects the fact that $\phi \psi$ cannot possibly be any smoother than $\phi$ or $\psi$ individually, while the condition $s \leq s_1 + s_2 - \frac{n}{2}$ arises from scaling considerations.

\section{$X^{s,b}$ estimates}\label{simple-sec}

Let $Z$ be either $\R^d$ or $\Z^d$ for some $d \geq 1$, and let $h_1$, $h_2$, $h_3$ be three functions from $Z$ to $\R$.  Let $m$ be a $[3;Z]$ multiplier; usually $m$ will be a symbol in the variables $\xi_1$, $\xi_2$, $\xi_3$.  Finally, let $b_1, b_2, b_3$ be real numbers.

We parameterize $Z \times \R$ by $(\xi,\tau)$.
In this section we study the problem of controlling the expression
\be{m-bound}
\| \frac{m(\xi_1,\xi_2,\xi_3)}{\prod_{j=1}^3 \langle \lambda_j \rangle^{b_j}} \|_{[3; Z \times \R]}
\end{equation}
where $\lambda_j$ is as in \eqref{lj-def}.

The discussion in this section will be fairly general, but in later sections we shall specialize to the dispersion relations $h_j$ which arise in KdV, Schr\"odinger, and wave equations.  These techniques could surely be adapted to hybrid systems (such as Zakharov, or gauge field equations such as Yang-Mills in the temporal gauge), or for other dispersive equations such as the KP-I and KP-II equations, but we will not pursue these matters here.

The function $h: \Gamma_3(Z) \to \R$ is defined by
$$ h(\xi_1,\xi_2,\xi_3) := h_1(\xi_1) + h_2(\xi_2) + h_3(\xi_3) = -\lambda_1 - \lambda_2 - \lambda_3$$
plays a fundamental role; it measures to what extent the spatial frequencies $\xi_1$, $\xi_2$, $\xi_3$ can resonate with each other.  Because of this, we shall refer to $h$ as the \emph{resonance function}.

Heuristically, we expect two types of frequency interactions $(\xi_1,\xi_2,\xi_3)$ to give a significant contribution to \eqref{m-bound}.  The first major contribution comes from \emph{resonant interactions} 
when the resonance function $h$ is zero, or close to zer. 
From a PDE viewpoint, a resonance describes two plane wave solutions to the linear problem which combine to form a third plane wave solution. Roughly speaking, when the zero set of $h$ is sufficiently simple, one can obtain good bounds on \eqref{m-bound} simply by applying the above tools, and performing dyadic decompositions away from the zero set.  When the zero set is more complicated, there is no single prescription for obtaining efficient bounds, but one must adapt the techniques to the geometry.  

A second major contribution comes from \emph{coherent interactions}, when one has $\nabla h_i(\xi_i) = \nabla h_j(\xi_j)$ for some $1 \leq i < j \leq 3$.  Equivalently, a coherent interaction occurs when at least two of the surfaces $\tau_j = h_j(\xi_j)$ fail to be transverse. From a PDE viewpoint, a coherence occurs when two parallel travelling wave solutions interact, or what is essentially equivalent, when many pairs of plane waves interact to create essentially the same plane wave output.  In dispersive situations, coherent interactions are rare (generally one only encounters this problem when $\xi_i = \pm \xi_j$, or at worst if $\xi_i$, $\xi_j$ are linearly dependent), but can still dominate \eqref{m-bound}, especially in low dimensions.

Generally speaking, $X^{s,b}$ norms are better at controlling the effects of resonance, whereas physical space norms (such as mixed Lebesgue norms) are better at controlling the effect of coherence.  In some situations (e.g. gauge field theories close to the critical regularity) it has been necessary to use a combination of both types of norm; see e.g. \cite{kl-mac:null3}, \cite{kl-tar:yang-mills}.  It is not clear at present what the best way to combine these two types of norms is, or whether completely new norms are needed.

In the KdV situation
\be{kdv}
h_1(\xi), h_2(\xi), h_3(\xi) := \xi^3, Z := \R \hbox{ or } \Z
\end{equation}
we have $h(\xi) = 3\xi_1 \xi_2 \xi_3$ and $\nabla h(\xi) = \xi^2$.  Thus one only has resonance when one of the $\xi_j$ vanishes, and one only has coherence when $\xi_i = \pm \xi_j$ for some $i,j$.  These are fairly simple criteria, and one should not need to apply any sophisticated techniques beyond a dyadic decomposition of the $\xi_1$, $\xi_2$, and $\xi_3$ variables.

In the non-periodic wave situation\footnote{Periodic wave equations have essentially the same behaviour as non-periodic wave equations when time is localized, thanks to finite speed of propagation.  The Klein-Gordon relation $h(\xi) = \pm (m^2 + |\xi|^2)^{1/2}$ also behaves similarly to \eqref{wave} when $m$ is bounded, time is localized, and frequency is large, but exhibits behaviour more reminiscent of the Schr\"odinger relation \eqref{schrodinger} in other situations.  We will not discuss multilinear estimates for the Klein-Gordon relation in this paper, but refer the reader to \cite{delort:fang} for further discussion.}
\be{wave}
h_1(\xi), h_2(\xi), h_3(\xi) := \pm |\xi|, Z := \R^d.
\end{equation}
we have coherence if and only if $\xi_1, \xi_2, \xi_3$ are constant multiples of each other. When the signs $\pm$ are not all the same, one also has resonance in this case.  When all the signs are the same, there are no resonant interactions except at the origin.  This suggests that one needs to perform a dyadic decomposition based upon the angular separation of $\xi_1$ and $\xi_2$ in addition to the more usual dyadic decomposition of $\xi_1$, $\xi_2$, and $\xi_3$ separately.

In the non-periodic Schr\"odinger situation
\be{schrodinger}
h_1(\xi), h_2(\xi), h_3(\xi) := \pm |\xi|^2, Z := \R^d .
\end{equation}
we have coherence only when $\pm \xi_i = \pm \xi_j$.  When all the signs agree there are no resonant interactions except at the origin.  However, if (for instance) $h_1$ and $h_2$ are positive and $h_3$ is negative, then one has resonance when $\xi_1$ and $\xi_2$ are perpendicular.  When $d=2$ the effect of resonance can be satisfactorily controlled by an angular dyadic decomposition of $\xi_1$ and $\xi_2$.  In higher dimensions the region of resonance is more complicated, and requires an induction on scales argument that we present in Section \ref{scale-sec}.  We will not discuss the periodic case in this paper as some substantial number theoretic issues arise in this case.

For more complicated equations (Zakharov, KP-I, KP-II, etc.) one can have a far more complicated zero set; cf. the discussion in \cite{cb:zak}.  To obtain sharp results in these equations one probably needs techniques such as those in Section \ref{scale-sec}.  There has been much recent progress on the well-posedness of these equations, see e.g. \cite{cks}, \cite{saut-tz:kp}, \cite{gtv}, \cite{tak-tz:kpii}.

We now make some general remarks concerning estimates of the form \eqref{m-bound}.  We first observe that we may restrict the multiplier to the region
\be{tj-large}
| \lambda_j | \gtrsim 1
\end{equation}
since the general case then follows by an averaging over unit time scales (\eqref{convolve} and the Comparison principle).  If $m$ is a symbol, then one can often assume
\be{ximax-large}
\max(|\xi_1|, |\xi_2|, |\xi_3|) \gtrsim 1
\end{equation}
for similar reasons, because the $|\xi_j| \lesssim 1$ behaviour of $m$ is usually identical to its $|\xi_j| \sim 1$ behaviour.  (If $m$ has mild singularities for $|\xi_j| \ll 1$, then one can usually still reduce to \eqref{ximax-large}; see Corollary \ref{smooth}).

By dyadic decomposition of the variables $\xi_j$, $\lambda_j$, as well as the function $h(\xi)$, we have
\be{m-2}
\eqref{m-bound} \lesssim 
\| \sum_{N_{max} \gtrsim 1} \sum_H \sum_{L_1,L_2,L_3 \gtrsim 1}
\frac{\tilde m(N_1,N_2,N_3)}{L_1^{b_1} L_2^{b_2} L_3^{b_3}} X_{N_1,N_2,N_3;H;L_1,L_2,L_3} \|_{[3, Z \times \R]}
\end{equation}
where $X_{N_1,N_2,N_3;H;L_1,L_2,L_3}$ is the multiplier
\be{x-def}
X_{N_1,N_2,N_3;H;L_1,L_2,L_3}(\xi,\tau) :=
\chi_{|h(\xi)| \sim H} \prod_{j=1}^3 \chi_{|\xi_j| \sim N_j} \chi_{|\lambda_j| \sim L_j}
\end{equation}
and
$$ \tilde m(N_1,N_2,N_3) := \sup_{|\xi_j| \sim N_j \forall j=1,2,3} |m(\xi_1,\xi_2,\xi_3)|.$$

The quantities $N_j$ and $L_j$ thus measure the spatial frequency of the $j^{th}$ wave and how closely it resembles a free solution respectively, while the quantity $H$ measures the amount of resonance.  

From the identities
$$ \xi_1 + \xi_2 + \xi_3 = 0$$
and
$$ \lambda_1 + \lambda_2 + \lambda_3 + h(\xi) = 0$$
on the support of the multiplier, we see that $X_{N_1,N_2,N_3;H;L_1,L_2,L_3}$ vanishes unless
\be{n-comp}
N_{max} \sim N_{med}
\end{equation}
and
\be{t-comp}
L_{max} \sim \max(H, L_{med}).
\end{equation}
Thus we may implicitly assume \eqref{n-comp}, \eqref{t-comp} in the summations. These reductions have the effect of simplifying \eqref{x-def} to roughly a tensor product of two functions (as opposed to a tensor product of three functions), so that the estimates of the previous section apply. 

Suppose for the moment that $N_1 \geq N_2 \geq N_3$, so by \eqref{n-comp} we have $N_1 \sim N_2 \gtrsim 1$.  Then as $N_1$ ranges over the dyadic numbers, the symbol
in the summation in \eqref{m-2} draws upon essentially disjoint regions of frequency space in both the 1 and 2 variables, and so Schur's test \eqref{sharp-schur} applies.  Similarly for permutations of $\{1,2,3\}$.  Applying this, we end up with
\bas
\eqref{m-bound} \lesssim 
\sup_{N \gtrsim 1} 
\| \sum_{ N_{max} \sim N_{med} \sim N} &\sum_H \sum_{L_{max} \sim \max(H,L_{med})}\\
&\frac{\tilde m(N_1,N_2,N_3)}{L_1^{b_1} L_2^{b_2} L_3^{b_3} } X_{N_1,N_2,N_3;H;L_1,L_2,L_3} \|_{[3, Z \times \R]}.
\end{align*}

In light of \eqref{t-comp} and the triangle inequality\footnote{In some endpoint cases one would use Lemma \ref{schur} instead of the triangle inequality at this juncture to replace one or more of the summations with a supremum; see e.g. the discussion of \eqref{endpoint}.}, we thus see that at least one of the inequalities
\be{m-bound-1-tri}
\begin{split}
\eqref{m-bound} \lesssim 
\sum_{ N_{max} \sim N_{med} \sim N} &
\sum_{L_1,L_2,L_3 \gtrsim 1}\\
&\frac{\tilde m(N_1,N_2,N_3)}{L_1^{b_1} L_2^{b_2} L_3^{b_3}} \| X_{N_1,N_2,N_3;L_{max};L_1,L_2,L_3} \|_{[3, Z \times \R]}
\end{split}
\end{equation}
or
\be{m-bound-2-tri}
\begin{split}
\eqref{m-bound} \lesssim 
\sum_{ N_{max} \sim N_{med} \sim N} &
\sum_{ L_{max} \sim L_{med}} \sum_{H \ll L_{max}}\\
&\frac{\tilde m(N_1,N_2,N_3)}{L_1^{b_1} L_2^{b_2} L_3^{b_3}} \| X_{N_1,N_2,N_3;H;L_1,L_2,L_3} \|_{[3, Z \times \R]}.
\end{split}
\end{equation}
hold for some $N \gtrsim 1$.  Of the two right-hand sides, \eqref{m-bound-2-tri} is generally easier to estimate, as the modulation $L_j$ are so large that one does not need to use much geometrical information about the surfaces $\tau_j = h_j(\xi_j)$, and also one usually has some decay\footnote{Indeed, one must have $b_i + b_j \geq 0$ and $b_1 + b_2 + b_3 \geq 1/2$, otherwise \eqref{m-bound} is automatically infinite for the same reasons that Proposition \ref{sobolev} is sharp.} arising from the denominator $L_1^{b_1} L_2^{b_2} L_3^{b_3}$.

\begin{figure}[htbp] \centering
\ \psfig{figure=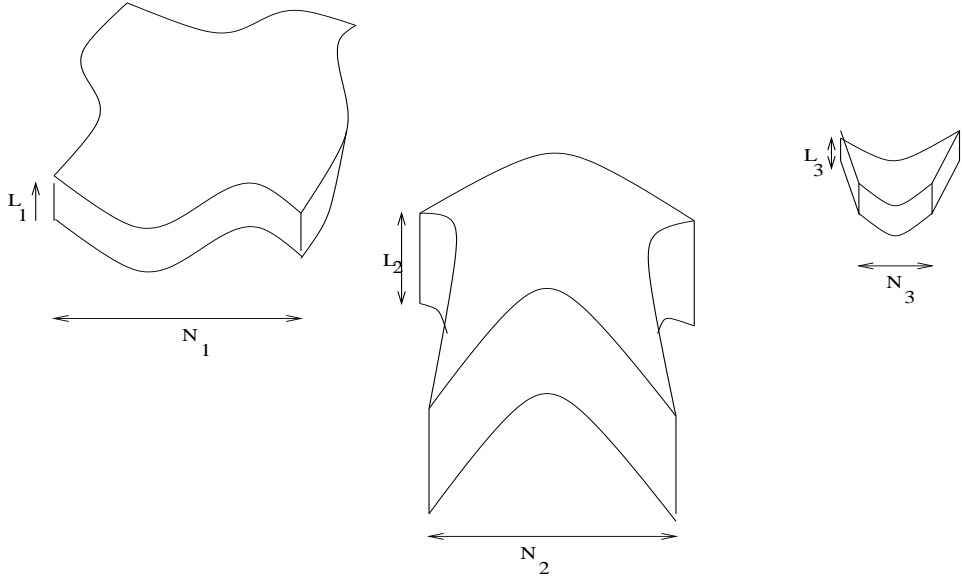,width=5in,height=3in}
\caption{A depiction of the multiplier in \eqref{h-gen}.  The three regions displayed are the sets on which the frequencies $(\xi_j, \tau_j)$ are constrained for $j=1,2,3$, although for sake of exposition we have drawn the annuli $|\xi_j| \sim N_j$ somewhat inaccurately as squares, and similarly for the constraint $|\lambda_j| \sim L_j$.  Because of the relations $\xi_1 + \xi_2 + \xi_3 = 
\tau_1 + \tau_2 + \tau_3 = 0$, it is often the case that some of these constraints are essentially redundant.  Of course, one expects the geometry of the dispersion relations $\tau_j = h_j(\xi_j)$ to play a major role in computing the multiplier norm of this object. }
\end{figure}

One is thus led to consider the expression
\be{h-gen}
\| X_{N_1,N_2,N_3;H;L_1,L_2,L_3} \|_{[3, Z \times \R]}
\end{equation}
in the low modulation case
\be{h-dom}
H \sim L_{max} 
\end{equation}
and the high modulation case
\be{t-dom}
L_{max} \sim L_{med} \gg H.
\end{equation}
Once one has good bounds on \eqref{h-gen}, the estimation of \eqref{m-bound} would then follow from \eqref{m-bound-1-tri} or \eqref{m-bound-2-tri} and some tedious computation of dyadic sums. 

The high modulation case \eqref{t-dom} is easier to handle.  For this discussion let us suppose that $L_1 \geq L_2 \geq L_3$.  In this case, the constraints $|\lambda_j| \sim L_j$ for $j=1,2$ are so weak that they have essentially no effect on the $[3;Z \times \R]$ norm.  With this philosophy in mind, we use the Comparison principle (Lemma \ref{comparison}) to estimate \eqref{h-gen} by
$$
\eqref{h-gen} \leq \| \chi_{|\lambda_3| \sim L_3} \chi_{|h(\xi)| \sim H} \prod_{j=1}^3 \chi_{|\xi_j| \sim N_j} \|_{[3,Z \times \R]}.$$
For fixed $\xi$, we have the one-dimensional estimate
$$ \| \chi_{|\lambda_3| \sim L_3} \|_{[3;\R]} \sim L_3^{1/2} \sim L_{min}^{1/2}$$
by the Tensored Box lemma.  By Lemma \ref{tensor}, we thus have
\be{t-dom-1}
\eqref{h-gen} \lesssim L_{min}^{1/2} 
\| \chi_{|h(\xi)| \sim H} \prod_{j=1}^3 \chi_{|\xi_j| \sim N_j} \|_{[3,Z]}.
\end{equation}
Although we derived \eqref{t-dom-1}  assuming $L_1 \geq L_2 \geq L_3$, it is clear from symmetry that \eqref{t-dom-1} in fact holds whenever \eqref{t-dom} does.

If we crudely estimate the multiplier in \eqref{t-dom-1} by $\chi_{|\xi_j| \sim N_{min}}$, where $N_j = N_{min}$, and use the Tensored Box lemma, we may estimate the above by
\be{t-dom-2}
\eqref{h-gen} \lesssim L_{min}^{1/2} 
| \{ \xi \in Z: |\xi| \sim N_{min} \} |^{1/2}.
\end{equation}
This rather crude estimate works surprisingly well in many cases (especially when the restriction $|h(\xi)| \sim H$ is redundant or very weak).

The case \eqref{h-dom} is more interesting, as it requires some geometric information about the surfaces $\tau = h_j(\xi)$.  As such, we only give a very general discussion here.

Suppose for the moment that $N_1 \geq N_2 \geq N_3$.  The $\xi_3$ variable in \eqref{h-gen} is currently localized to the annulus $\{ |\xi_3| \sim N_{min}\}$.  By a finite partition of unity we can restrict it further to a ball $\{ |\xi_3 -\xi_3^0| \ll N_{min} \}$ for some $|\xi_3^0| \sim N_{min}$.  But then by Box Localization (Lemma \ref{boxes}) we may localize $\xi_1, \xi_2$ similarly to regions
$$ \{ \xi_1: |\xi_1 - \xi_1^0| \ll N_{min} \}; \quad \{ \xi_2: |\xi_2 - \xi_2^0| \ll N_{min} \}$$
where $\xi_j^0 \sim N_j$.  We may assume that $|\xi_1^0 + \xi_2^0 + \xi_3^0| \ll N_{min}$ since the symbol vanishes otherwise.  We may summarize this symmetrically as
\be{wander}
\eqref{h-gen} \lesssim
\|
\chi_{|h(\xi)| \sim H} \prod_{j=1}^3 \chi_{|\xi_j - \xi_j^0| \ll N_{min}} \chi_{|\lambda_j| \sim L_j}
\|_{[3;Z \times \R]}
\end{equation}
for some $\xi_1^0, \xi_2^0, \xi_3^0$ satisfying 
\be{xio-cond}
|\xi_j^0| \sim N_j \hbox{ for } j=1,2,3; \quad |\xi_1^0 + \xi_2^0 + \xi_3^0| \ll N_{min}.
\end{equation}
Although we derived \eqref{wander} assuming that $N_1 \geq N_2 \geq N_3$, it is clear that one in fact has \eqref{wander} for all choices of $N_1$, $N_2$, $N_3$.

\begin{figure}[htbp] \centering
\ \psfig{figure=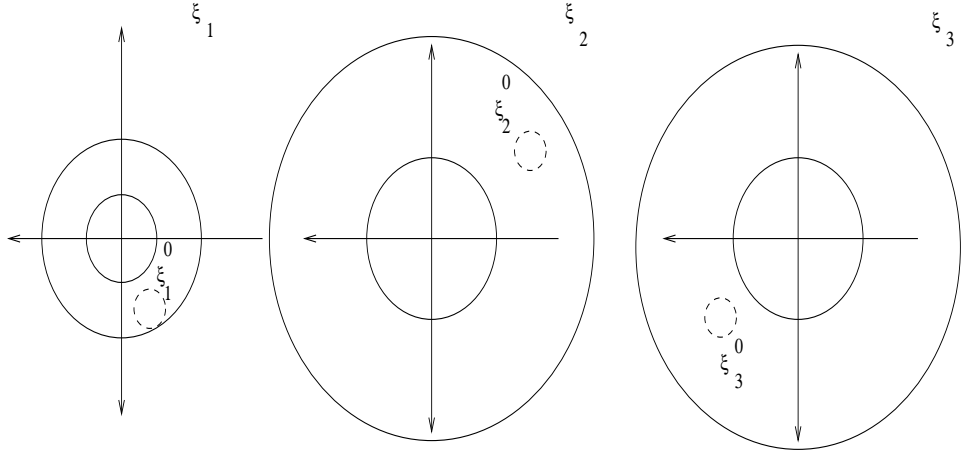,width=5in,height=2.4in}
\caption{The effect of \eqref{wander} is to localize the $\xi_1$, $\xi_2$, $\xi_3$ variables to balls of radius $\ll N_{min}$.  We can assume that the centers of these balls (essentially) add up to zero. }
\end{figure}

Suppose now that $L_1 \geq L_2 \geq L_3$.  The condition $|\lambda_1| \sim L_1 \sim H$ is so weak as to almost be redundant, so we shall use the Comparison principle to remove it.  The constraint $|\xi_1 - \xi_1^0| \ll N_{min}$ is also almost redundant given the similar constraints on $\xi_2$, $\xi_3$.  Also, there is little point now in distinguishing between $|\lambda_j| \sim L_j$ and $|\lambda_j| \lesssim L_j$, so we have
\be{wander-2}
\eqref{h-gen} \lesssim
\|
\chi_{|h(\xi)| \sim H} \prod_{j=2}^3 \chi_{|\xi_j - \xi_j^0| \ll N_{min}} \chi_{|\lambda_j| \lesssim L_j}
\|_{[3;Z \times \R]}
\end{equation}
The right-hand side of \eqref{wander-2} is almost a tensor product of two characteristic functions, which suggests discarding the $|h| \sim H$ constraint and using the Characteristic function estimate \eqref{char-zar}.  In simple situations such as the KdV case \eqref{kdv} or the $(+++)$ Schr\"odinger case \eqref{schrodinger} this procedure will give sharp results, because the constraint $|h(\xi)| \sim H$ is redundant in those cases.  However, in the other situations listed above, the constraint $|h(\xi)| \sim H$ has a non-trivial effect, and further treatment (e.g. angular frequency decompositions, or Lemma \ref{cauchy-schwarz}) is needed.

We observe the estimate

\begin{lemma}\label{char-1}
Let $A, B$ be subsets of $\Z$, and $L_1, L_2 > 0$.  Then
\bas
\| \chi_A(\xi_1) \chi_B(\xi_2)&
 \chi_{|\lambda_1| \lesssim L_1; |\lambda_2| \lesssim L_2} 
\|_{[3; Z \times \R]}
\lesssim
\min(L_1,L_2)^{1/2}\\
& | \{  \xi_1 \in A:  \xi - \xi_1 \in B; 
h_1(\xi_1) + h_2(\xi - \xi_1) = \tau + O(\max(L_1,L_2)) \}|^{1/2}
\end{align*}
for some $\xi \in Z$, $\tau \in \R$.
\end{lemma}

The right-hand side measures the size of the intersection of the hypersurface $\tau_1 = h_1(\xi_1)$ with some inverted translate of the surface  $\tau_2 = h_2(\xi_2)$.

\begin{proof}
From \eqref{char-zar} we can estimate the left-hand side by
$$
| \{ (\xi_1,\tau_1) \in Z: \xi_1 \in A; \xi - \xi_1 \in B; 
\tau_1 = h_1(\xi_1) + O(L_1);  \tau - \tau_1 = h_2(\xi - \xi_1) + O(L_2)) \}|^{1/2}.$$
For fixed $\xi_1$, the set of possible $\tau_1$ ranges in an interval of length $O(\min(L_1,L_2))$, and vanishes unless $h_1(\xi_1) + h_2(\xi - \xi_1) = \tau + O(\max(L_1,L_2))$.  The claim follows.
\end{proof}

If we can afford to ignore the constraint $|h(\xi)| \sim H$ in \eqref{wander-2}, we thus have

\begin{corollary}\label{char-2}
Let $N_1, N_2, N_3 > 0$, $L_1 \geq L_2 \geq L_3$, and $\xi^1_0, \xi^2_0, \xi^3_0$ satisfy \eqref{xio-cond}.  Then we have
$$
\eqref{h-gen}
\lesssim
L_3^{1/2} | \{ \xi_2 \in Z: |\xi_2 - \xi_2^0| \ll N_{min}; h_2(\xi_2) + h_3(\xi - \xi_2) = \tau + O(L_2) \}|^{1/2}
$$
for some $\tau \in \R$ and $\xi \in Z$ with $|\xi + \xi_1^0| \ll N_{min}$.  Similar statements hold with the roles of the indices 1,2,3 permuted.
\end{corollary}

This estimate is already enough to accurately estimate \eqref{h-gen} in many cases, such as the KdV case and the $(+++)$ Schr\"odinger case (see below).

\section{An averaging argument}

In the previous section we reduced the problem of computing norms such as \eqref{m-bound} to that of computing dyadic summations such as \eqref{m-bound-1-tri}, \eqref{m-bound-2-tri}.  Although these dyadic summations are always computable (once the quantities \eqref{h-gen} have been evaluated), the sheer number of possible cases depending on the relative sizes of $N_1, N_2, N_3$ and of $L_1, L_2, L_3$ makes the evaluation of these summations somewhat tedious.  Fortunately, in many cases one can eliminate many of these cases if the exponents $b_j$ are sufficiently large.

In particular, if $b_j > 1/2$ for some $j$, then one expects to only need to consider the case $L_j \sim 1$, if one adopts the heuristic that $X^{s,1/2+}$ functions behave for short time like free solutions.  This heuristic is implicitly used throughout the literature, and appears to have been first noted by Bourgain \cite{borg:hsd}.  A rigorous version of this heuristic is:

\begin{proposition}\label{average}
Let $Z$ be an abelian group, and let $Z \times \R$ be parameterized by $(\xi,\tau)$ for $\xi \in Z$, $\tau \in \R$.  Let $h_1, h_2: Z \to \R$ be functions, and let $m: (Z \times \R)^k \to \R$ be a non-negative function which is constant in the $\tau_1$ and $\tau_2$ variables.  Let $b_1, b_2 \in \R$ be such that
\be{av-hyp}
b_1 > 1/2, 1/2 + b_2; \quad b_1 \geq -b_2.
\end{equation}
Then
\be{av-eq}
\| \frac{m(\xi,\tau)}{\langle \lambda_1 \rangle^{b_1}
\langle \lambda_2 \rangle^{b_2}} \|_{[k;Z \times \R]}
\sim
\| \frac{m(\xi,\tau) \chi_{|\lambda_1| \sim 1}}{
\langle \lambda_2 \rangle^{b_2}} \|_{[k;Z \times \R]}
\end{equation}
The implicit constants depend on $b_1$, $b_2$.
\end{proposition}

The $b_1 \geq -b_2$ condition is necessary as the left-hand side of \eqref{av-eq} is automatically infinite otherwise.  The condition $b_1 > 1/2 + b_2$ is somewhat unsatisfactory, as in applications one often has $b_1, b_2$ close to $1/2$.  For instance, this condition is responsible for Proposition \ref{qij} being 1/4 of a derivative away from the critical regularity.  However, there are still several situations in which this Proposition can reduce the number of cases substantially, as it essentially restricts one or more of the $L_j$ variables to equal 1.

\begin{proof} For brevity we shall denote the right-hand side of \eqref{av-eq} as $X$.  The lower bound for $X$ follows immediately from the Comparison principle, so it suffices to show the upper bound.  

The idea shall be to decompose $\lambda_1$, $\lambda_2$ into dyadic shells, and then use Lemma \ref{translation} to move $\lambda_1$ to be $\sim 1$.  This may move $\lambda_2$ into another dyadic shell, depending on the relative sizes of $\lambda_1$ and $\lambda_2$.

By \eqref{convolve} we may assume that $|\lambda_1|, |\lambda_2| \gtrsim 1$.  We then split \eqref{av-eq} into three pieces determined by the regions
$$ 1 \lesssim |\lambda_1| \ll |\lambda_2|$$
$$ 1 \lesssim |\lambda_2| \ll |\lambda_1|$$
$$ 1 \lesssim |\lambda_1| \sim |\lambda_2|.$$

In the first case we dyadically decompose $|\lambda_1|$ and use the triangle inequality to estimate the left-hand side of \eqref{av-eq} by
$$ \sum_{L \gtrsim 1} L^{-b_1}
\| \frac{m(\xi,\tau) \chi_{|\lambda_1| \sim L} \chi_{|\lambda_2| \gg L}
}{\langle \lambda_2 \rangle^{b_2}} \|_{[k;Z \times \R]}.$$
We subdivide $\chi_{|\lambda_1| \sim L}$ into $O(L)$ regions of the form $\chi_{|\lambda_1- \lambda| \sim 1}$ for integers $|\lambda| \sim L$.  By Lemma \ref{schur} we can thus estimate the previous by
$$ \sum_{L \gtrsim 1} L^{-b_1} L^{1/2} \sup_{|\lambda| \sim L}
\| \frac{m(\xi,\tau) \chi_{|\lambda_1 - \lambda| \sim 1} \chi_{|\lambda_2| \gg L}
}{\langle \lambda_2 \rangle^{b_2}} \|_{[k;Z \times \R]}.$$
By Lemma \ref{translation} and the Comparison principle the expression inside the sup is $O(X)$.  Since $b_1 > 1/2$, the claim follows for this case.

In the second case we repeat the above arguments, eventually estimating this contribution by
$$ \sum_{L \gtrsim 1} L^{-b_1} L^{1/2} \sup_{|\lambda| \sim L}
\| \frac{m(\xi,\tau) \chi_{|\lambda_1 - \lambda| \sim 1} \chi_{|\lambda_2| \ll L}
}{\langle \lambda_2 \rangle^{b_2}} \|_{[k;Z \times \R]}.$$
If we use Lemma \ref{translation} to shift $\tau_1$ down by $\lambda$ and $\tau_2$ up by $\lambda$, and use the crude estimate
$$
\langle \lambda_2 \rangle^{b_2} \gtrsim \min(1, L^{b_2}) \sim \min(L^{-b_2},1)
\langle \lambda_2 + \lambda\rangle^{b_2},$$
we can estimate the expression inside the sup by $O(\max(L^{b_2},1) X)$.  Since $b_1 > 1/2, 1/2 + b_2$, the claim then follows for this case.

It remains to consider the third case.  We can dyadically decompose into pieces $|\lambda_1| \sim |\lambda_2| \sim L$ for $L \gtrsim 1$.  By Schur's test \eqref{sharp-schur} it suffices to control the contribution of a single $L$, which we now fix.

The quantity $\langle h_1(\xi_1) + h_2(\xi_2) \rangle$ fluctuates between $1$ and $L$.  We then dyadically decompose in this quantity and use the triangle inequality to estimate the contribution of this case by
$$
\sum_{1 \lesssim M \lesssim L} L^{-b_1-b_2}
\| m(\xi,\tau) \chi_{|\lambda_1| \sim L} \chi_{|\lambda_2| \sim L} \chi_{\langle h_1(\xi_1) + h_2(\xi_2) \rangle \sim M}
\|_{[k;Z \times \R]}.$$
We now decompose the condition $|\lambda_1| \sim L$ into $|\lambda_1 - \lambda| \sim 1$ for $|\lambda| \sim L$ as before.  Using Schur's test we can estimate the previous by
$$
\sum_{1 \lesssim M \lesssim L} L^{-b_1-b_2} M^{1/2}
\sup_{|\lambda| \sim M}
\| m(\xi,\tau) \chi_{|\lambda_1 - \lambda| \sim 1} \chi_{|\lambda_2| \sim L} \chi_{\langle h_1(\xi_1) + h_2(\xi_2) \rangle \sim M}
\|_{[k;Z \times \R]}.$$
By Lemma \ref{translation} the expression inside the supremum is $O(M^{b_2}X)$.  Since $b_1 > 1/2$ and $b_1 \geq -b_2$, the claim then follows.
\end{proof}

If only some of the hypotheses in \eqref{av-hyp} hold then we can achieve some partial reductions.  For instance, if $b_1 > 1/2$ and $b_1 \geq -b_2$, then we can reduce to one of the two cases $L_1 \sim 1$ or $L_1 \gg L_2$.  Conversely, if one only assumes $b_1 > 1/2 + b_2$, then one can eliminate the case $L_1 \gg L_2$ but must still deal with the cases when $L_1 \lesssim L_2$.  These reductions have a slight simplifying effect on many of the dyadic summations under consideration, but we shall not exploit these in the sequel. 

A variant on the above theme states that if $b_i + b_j > 0$ for all $i \neq j$, then one can eliminate the case \eqref{t-dom} and reduce to \eqref{h-dom}, except when $H \ll 1$ of course.  Although this heuristic can be made rigorous by a variant of the above arguments, we will not do so here, especially since the case \eqref{t-dom} can usually be dealt with quite easily, and in any event one still has to deal with \eqref{t-dom} in the $H \ll 1$ case. 

\section{Estimates related to the KdV equation}\label{kdv-sec}

In this section we specialize \eqref{m-bound} to the KdV dispersion relationship \eqref{kdv} when $Z = \R$ or $Z = \Z$.  This is the easiest of the three cases to study as space is one-dimensional, and so angular issues do not arise.  Also, since the cubic $\tau = \xi^3$ is odd, one does not need to distinguish between $(+++)$ and $(++-)$ cases.

Following the general philosophy of Section \ref{simple-sec}, we begin with the study of \eqref{h-gen}.  From the resonance identity
\be{kdv-fund}
h(\xi) = \xi_1^3 + \xi_2^3 + \xi_3^3 = 3 \xi_1 \xi_2 \xi_3
\end{equation}
we see that we may assume that 
\be{h-comp-kdv}
H \sim N_1 N_2 N_3
\end{equation}
since the multiplier in \eqref{h-gen} vanishes otherwise.  The constraint $\chi_{|h(\xi)| \sim H}$ is now redundant and will be discarded.  

We can now compute \eqref{h-gen} easily from the discussion of the previous section.  The coherent cases $\xi_i = \pm \xi_j$ are exceptional, and the estimates are rather unfavorable, but all the cases are relatively easy to compute.  To unify the periodic and non-periodic cases, we adopt the notation that $\langle x \rangle_Z$ is $|x|$ if $Z = \R$ and $1 + |x|$ if $Z = \Z$.  Note that the measure of an interval $\{ x \in Z: a \leq x \leq b\}$ is $O(\langle b-a \rangle_Z)$ for both choices of $Z$.

\begin{figure}[htbp] \centering
\ \psfig{figure=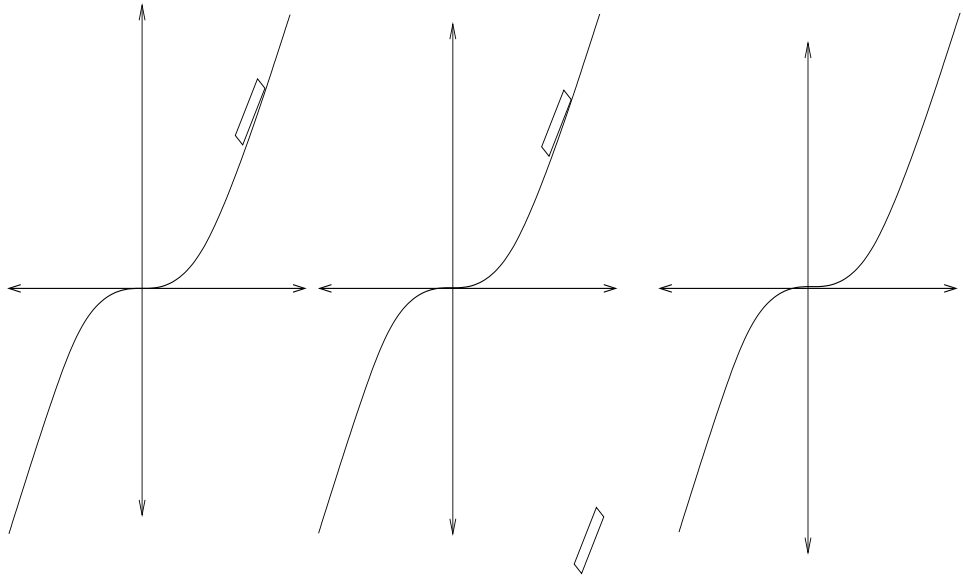,width=5in,height=3in}
\caption{The Knapp example which shows that \eqref{kdv-excep} is sharp in the (++)-coherent case.}
\end{figure}

\begin{proposition}\label{kdv-core}
Let $H, N_1, N_2, N_3, L_1, L_2, L_3 > 0$ obey \eqref{n-comp}, \eqref{t-comp}, \eqref{h-comp-kdv}. 
\begin{itemize}
\item ((++) Coherence) If $N_{max} \sim N_{min}$ and $L_{max} \sim H$, then we have
\be{kdv-excep}
\eqref{h-gen} \lesssim 
L_{min}^{1/2}
\langle N_{max}^{-1/4} L_{med}^{1/4} \rangle_Z.
\end{equation}
\item ((+-) Coherence) If $N_2 \sim N_3 \gg N_1$ and $H \sim L_1 \gtrsim L_2, L_3$, then
\be{kdv-weird}
\eqref{h-gen} \lesssim L_{min}^{1/2}
\langle 
N_{max}^{-1} \min(H, \frac{N_{max}}{N_{min}} L_{med})^{1/2}
\rangle_Z.
\end{equation}
Similarly for permutations.
\item In all other cases, we have
\be{kdv-standard}
\eqref{h-gen} \lesssim L_{min}^{1/2}
\langle 
N_{max}^{-1} \min(H,L_{med})^{1/2}
\rangle_Z.
\end{equation}
\end{itemize}
\end{proposition}

\begin{figure}[htbp] \centering
\ \psfig{figure=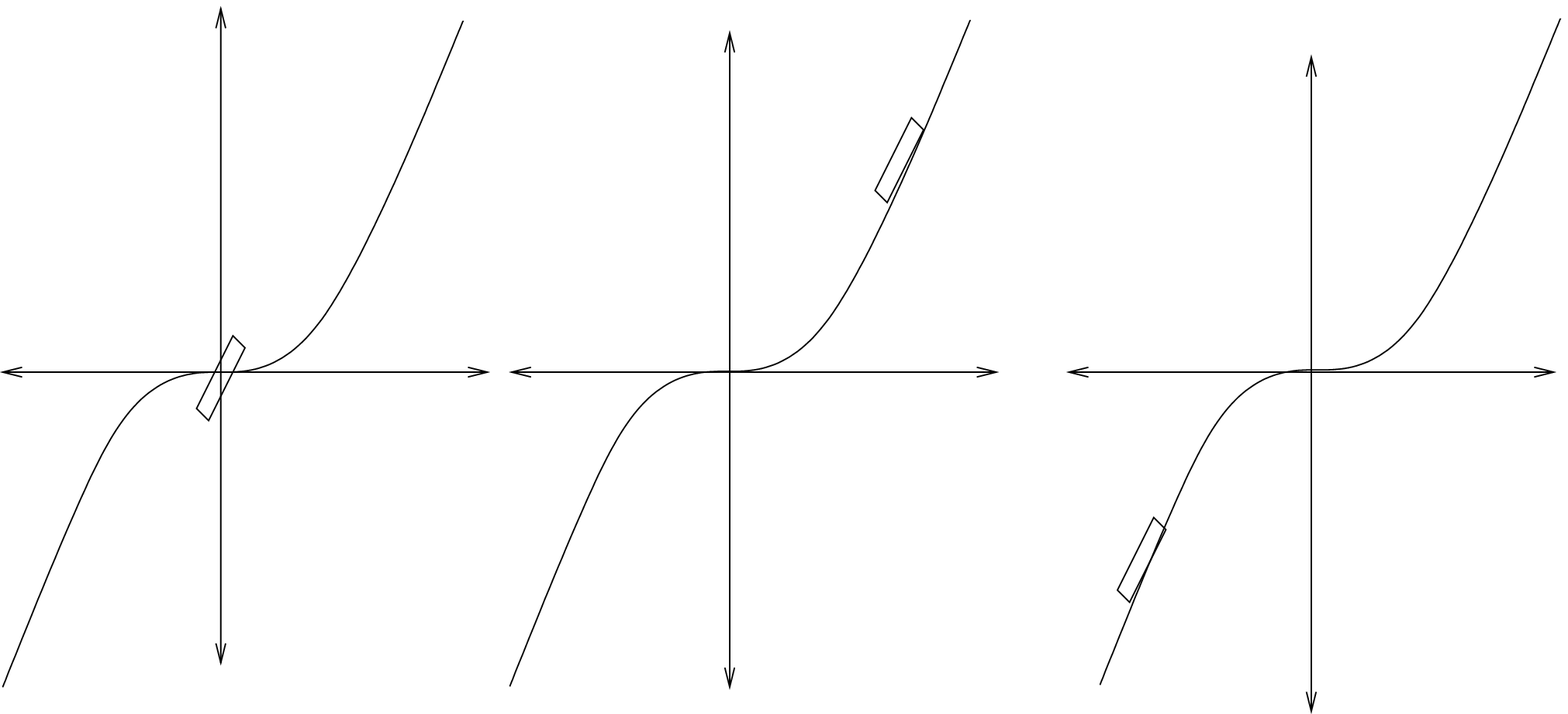,width=5in,height=3in}
\caption{The example which shows that \eqref{kdv-weird} is sharp in the (+-)-coherent case.}
\end{figure}

\begin{proof}  
In the case \eqref{t-dom} we have
$$ \eqref{h-gen} \lesssim L_{min}^{1/2} \langle N_{min}^{1/2}\rangle_Z$$
by \eqref{t-dom-2}, from which the claimed bounds follow.  It thus remains to consider the case \eqref{h-dom}.
By symmetry we may assume that $L_1 \geq L_2 \geq L_3$.

By Corollary \ref{char-2}, we have
$$  \eqref{h-gen} \lesssim L_3^{1/2} | 
\{ \xi_2: |\xi_2 - \xi_2^0| \ll N_{min}; |\xi - \xi_2 - \xi_3^0| \ll N_{min};
(\xi_2)^3 + (\xi - \xi_2)^3 = \tau_2 + O(L_2)
\}|^{1/2}$$
for some $\xi$, $\xi_1^0$, $\xi_2^0$, $\xi_3^0$ satisfying \eqref{xio-cond} and $|\xi + \xi_1^0| \ll N_{min}$.

To compute the right-hand side of this expression we shall use the identity
\be{trio}
(\xi_2)^3 + (\xi - \xi_2)^3 = 3\xi (\xi_2 - \frac{\xi}{2})^2 + \frac{\xi^3}{4}.
\end{equation}

We need only consider three cases: $N_1 \sim N_2 \sim N_3$, $N_1 \sim N_2 \gg N_3$, and $N_2 \sim N_3 \gg N_1$.  (The case $N_1 \sim N_3 \gg N_2$ then follows by symmetry).

If $N_1 \sim N_2 \sim N_3$, we see from \eqref{trio} that $\xi_2$ variable is contained in the union of two intervals of length $O(N_1^{-1/2} L_2^{1/2})$ at worst, and \eqref{kdv-excep} follows.  

If $N_1 \sim N_2 \gg N_3$, we must have $|\xi_2 - \frac{\xi}{2}| \sim N_1$, so \eqref{trio} shows that $\xi$ is contained in the union of two intervals of length $O(N_1^{-2} L_2)$, and \eqref{kdv-standard} follows.

If $N_2 \sim N_3 \gg N_1$, then we must have $|\xi_2 - \frac{\xi}{2}| \sim N_2$, so \eqref{trio} shows that $\xi$ is contained in the union of two intervals of length $O(N_1^{-1} N_2^{-1} L_2)$.  But $\xi_2$ is also contained in an interval of length $\ll N_1$.  The claim \eqref{kdv-weird} follows.
\end{proof}

\begin{figure}[htbp] \centering
\ \psfig{figure=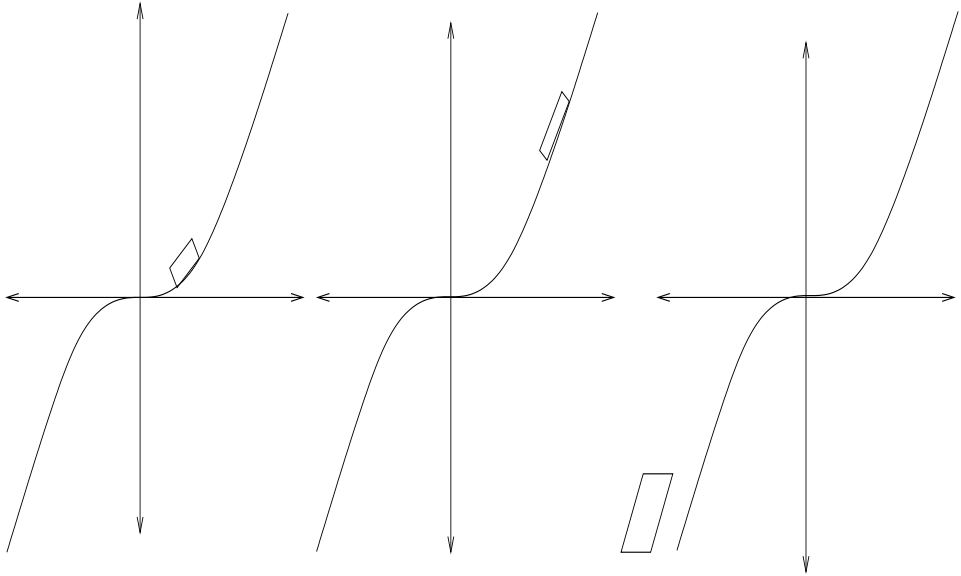,width=5in,height=3in}
\caption{The example which shows that \eqref{kdv-standard} is sharp in the non-coherent case.}
\end{figure}

The bounds listed above are sharp.  For \eqref{kdv-standard}, this is obtained by testing \eqref{czk-def} with 
$$ f_j(\xi,\tau) := \chi_{|\xi| \sim N_j; |\lambda| \sim L_j}$$
for $j=1,2,3$, with $\lambda := \tau - \xi^3$ as usual.  For \eqref{kdv-excep}, we use the ``Knapp example''
\bas
f_1(\xi,\tau) &:= \chi_{|\xi - N_1| \lesssim N_1^{-1/2} L_2^{1/2}; |\lambda| \lesssim L_1}\\
f_2(\xi,\tau) &:= \chi_{|\xi - N_1| \lesssim N_1^{-1/2} L_2^{1/2}; |\lambda| \lesssim L_2}\\
f_3(\xi,\tau) &:= \chi_{|\xi - 2N_1| \lesssim N_1^{-1/2} L_2^{1/2}; |\tau - \xi^3/4| \lesssim L_2}
\end{align*}
in the case $L_1 \leq L_2 \leq L_3$, and similarly for the permutations of this case.  Finally, for \eqref{kdv-weird}, we use 
\bas
f_1(\xi,\tau) &:= \chi_{|\xi| \sim N_1; |\tau - 3N_2^2 \xi| \lesssim N_1^2 N_2}
\\
f_2(\xi,\tau) &:= \chi_{|\xi - N_2| \lesssim N_1; |\lambda| \lesssim L_2}\\
f_3(\xi,\tau) := \chi_{|\xi + N_2| \lesssim N_1; |\lambda| \lesssim L_3}
\end{align*}
in the case $N_2 \sim N_3 \gtrsim N_1$ and $H \sim L_1 \gtrsim L_2 \gtrsim L_3$, and similarly for permutations of this case.  We omit the details.

With Proposition \ref{kdv-core} one can now prove sharp bilinear estimates in both the periodic and non-periodic setting.  We illustrate this with an asymmetric bilinear estimate on the real line, which can be viewed as a bilinear improvement to the Strichartz embedding $X^{-1/8,1/2+\eps}_{\tau = \xi^3}(\R \times \R) \subset L^4(\R \times \R)$.  We will then use this bilinear estimate to derive a trilinear estimate.

\begin{proposition}\label{kpv-2}
For all $u,v$ on $\R \times \R$ and $0 < \eps \ll 1$, we have
$$ \| uv\|_{L^2(\R \times \R)} \lesssim
\| u\|_{X^{-1/2,1/2-\eps}_{\tau = \xi^3}(\R \times \R)}
\| v\|_{X^{1/4,1/2+\eps}_{\tau = \xi^3}(\R \times \R)}$$
\end{proposition}

\begin{proof}
By Plancherel it suffices to show that
\be{kpv-2-est}
\|
\frac{\langle \xi_2 \rangle^{1/2} }
{
\langle \xi_1 \rangle^{1/4}
\langle \tau_2 - \xi_2^3 \rangle^{1/2 - \eps}
\langle \tau_1 - \xi_1^3 \rangle^{1/2+\eps}
}
\|_{[3;\R \times \R]} \lesssim 1.
\end{equation}
From \eqref{m-bound-1-tri} or \eqref{m-bound-2-tri} it suffices to show that
\be{kdv-3-tri}
\begin{split}
\sum_{ N_{max}\sim N_{med} \sim N} 
&\sum_{L_1,L_2,L_3 \gtrsim 1}
\frac{ \langle N_2 \rangle^{1/2} }{ \langle N_1\rangle^{1/4}
L_1^{1/2+\eps} L_2^{1/2-\eps} }\\
&\| X_{N_1,N_2,N_3;L_{max};L_1,L_2,L_3} \|_{[3, Z \times \R]}
\lesssim 1
\end{split}
\end{equation}
and
\be{kdv-4-tri}
\begin{split}
\sum_{ N_{max} \sim N_{med} \sim N} 
&\sum_{ L_{max} \sim L_{med}} \sum_{H \ll L_{max}}
\frac{ \langle N_2 \rangle^{1/2} }{ \langle N_1\rangle^{1/4}
L_1^{1/2+\eps} L_2^{1/2-\eps} }\\
&\| X_{N_1,N_2,N_3;H;L_1,L_2,L_3} \|_{[3, Z \times \R]}
\lesssim 1
\end{split}
\end{equation}
for all $N \gtrsim 1$.  This will be accomplished by Lemma \ref{kdv-core} and some tedious summation.

Fix $N$.  We first prove \eqref{kdv-4-tri}.  We may assume \eqref{h-comp-kdv}.  By \eqref{kdv-standard} we reduce to
\bas
\sum_{ N_{max} \sim N_{med} \sim N} 
&\sum_{ L_{max} \sim L_{med} \gtrsim N_1 N_2 N_3}\\
&\frac{ \langle N_2 \rangle^{1/2} }{ \langle N_1\rangle^{1/4}
L_1^{1/2+\eps} L_2^{1/2-\eps} }
L_{min}^{1/2} N^{-1} (N_1 N_2 N_3)^{1/2}
\lesssim 1.
\end{align*}
Estimating
$$ \frac{\langle N_2 \rangle^{1/2} } {\langle N_1 \rangle^{1/4} }
\lesssim \frac{N^{1/2}}{\langle N_{min} \rangle^{1/4}}; \quad
L_1^{1/2+\eps} L_2^{1/2-\eps} \gtrsim 
L_{min}^{1/2+\eps} L_{med}^{1/2-\eps} 
\sim L_{min}^{1/2+\eps} (N_1 N_2 N_3)^{1/2-\eps}$$
and then performing the $L$ summations, we reduce to
$$\sum_{ N_{max} \sim N_{med} \sim N} 
\frac{ N^{1/2} }{ \langle N_{min}\rangle^{1/4} N^{-1} (N_1 N_2 N_3)^\eps }
\lesssim 1.$$
which is certainly true (with about $N^{-1/2}$ to spare).

Now we show \eqref{kdv-3-tri}.  We may assume $L_{max} \sim N_1 N_2 N_3$.

We first deal with the contribution where \eqref{kdv-excep} holds.  In this case we have $N_1, N_2, N_3 \sim N \gtrsim 1$, so we reduce to
$$
\sum_{ L_{max} \sim N^3}
\frac{N^{1/2} }{ N^{1/4}
L_{min}^{1/2+\eps} L_{med}^{1/2-\eps} }
L_{min}^{1/2} N^{-1/4} L_{med}^{1/4} \lesssim 1.
$$
But this is easily verified.

Now we deal with the cases where \eqref{kdv-weird} applies. We do not have perfect symmetry and must consider the cases
\bas
 N \sim N_1 \sim N_2 \gg N_3&; H \sim L_3 \gtrsim L_1, L_2\\
 N \sim N_2 \sim N_3 \gg N_1&; H \sim L_1 \gtrsim L_2, L_3\\
 N \sim N_1 \sim N_3 \gg N_2&; H \sim L_2 \gtrsim L_1, L_3
\end{align*}
separately.

In the first case we reduce by \eqref{kdv-weird} to
$$
\sum_{N_3 \ll N } 
\sum_{1 \lesssim L_1,L_2 \lesssim N^2 N_3} 
\frac{ N^{1/2} }{ N^{1/4}
L_1^{1/2+\eps} L_2^{1/2-\eps} }
L_{min}^{1/2} N^{-1}
\min(N^2 N_3, \frac{N}{N_3} L_{med})^{1/2}
\lesssim 1.$$
Performing the $N_3$ summation we reduce to
$$
\sum_{1 \lesssim L_1,L_2 \lesssim N^3} 
\frac{ N^{1/2} }{ N^{1/4}
L_1^{1/2+\eps} L_2^{1/2-\eps} }
L_{min}^{1/2} N^{-1}
N^{3/4} L_{med}^{1/4}
\lesssim 1$$
which is easily verified.

To unify the second and third cases we replace $L_1^{1/2+\eps}$ by $L_1^{1/2-\eps}$.  By asymmetry it suffices now to show the second case. We simplify using the first half of \eqref{kdv-weird} to
$$
\sum_{N_1 \ll N } 
\sum_{1 \lesssim L_2,L_3 \ll N^2 N_1}
\frac{N^{1/2}}{ \langle N_1\rangle^{1/4}
(N^2 N_1)^{1/2-\eps} L_2^{1/2-\eps} }
L_{min}^{1/2}
N_1^{1/2}
\lesssim 1.
$$
We may assume $N_1 \gtrsim N^{-2}$ since the inner sum vanishes otherwise.  Performing the $L$ summation we reduce to
$$
\sum_{N^{-2} \lesssim N_1 \ll N } 
\frac{N^{1/2}}{ \langle N_1\rangle^{1/4}
(N^2 N_1)^{1/2-2\eps} }
N_1^{1/2}
\lesssim 1
$$
which is easily verified (with about $N^{-1/2}$ to spare).

To finish the proof of \eqref{kdv-3-tri} it remains to deal with the cases where \eqref{kdv-standard} holds.  This reduces to
$$
\sum_{ N_{max} \sim N_{med} \sim N} 
\sum_{ L_{max} \sim N_1 N_2 N_3}
\frac{ \langle N_2 \rangle^{1/2} }{ \langle N_1\rangle^{1/4}
L_1^{1/2+\eps} L_2^{1/2-\eps} }
L_{min}^{1/2} N^{-1} L_{med}^{1/2} \lesssim 1.
$$
Performing the $L$ summations, we reduce to
$$
\sum_{ N_{max} \sim N_{med} \sim N} 
\frac{ \langle N_2 \rangle^{1/2} (N_1 N_2 N_3)^{\eps}}{ \langle N_1\rangle^{1/4} } N^{-1} \lesssim 1$$
which is easily verified (with about $N^{-1/2}$ to spare).
\end{proof}

One can of course prove many other bilinear $X^{s,b}$ estimates of KdV type from Lemma \ref{kdv-core}; for instance, the bilinear estimates in \cite{kpv:kdv}, \cite{cst}, \cite{ckstt} can also be deduced by the above techniques.  We will not attempt to give an exhaustive characterization here of all such bilinear estimates due to the prohibitive number of cases (especially at the endpoints).

From the above bilinear estimate one can deduce the following trilinear estimate:

\begin{corollary}\label{kpv-3}
For all $u_1,u_2,u_3$ on $\R \times \R$ and $0 < \eps \ll 1$, we have
$$
\| (u_1 u_2 u_3)_x \|_{X^{1/4,-1/2+\eps}_{\tau = \xi^3}(\R \times \R)} \lesssim
\prod_{j=1}^3 \| u_j \|_{X^{1/4,1/2+\eps}_{\tau = \xi^3}(\R \times \R)}
$$
with the implicit constant depending on $\eps$.
\end{corollary}

This estimate can be used to give an alternate proof of the local well-posedness of the mKdV equation in $H^s(\R)$ for $s \geq 1/4$, which was first shown in \cite{kpv:periodic} using maximal function and Kato smoothing estimates.  The $1/4$ exponent is sharp; see \cite{kpv:kdv}.  

\begin{proof}
By duality and Plancherel it suffices to show that
$$
\|
\frac{(\xi_1 + \xi_2 + \xi_3) \langle \xi_4 \rangle^{1/4}}
{
\langle \tau_4 - \xi_4^3 \rangle^{1/2 - \eps}
\prod_{j=1}^3 \langle \xi_j \rangle^{1/4}
\langle \tau_j - \xi_j^3 \rangle^{1/2+\eps}
}
\|_{[4;\R \times \R]} \lesssim 1.
$$
We estimate $|\xi_1 + \xi_2 + \xi_3|$ by $\langle \xi_4\rangle$.  We then apply the inequality\footnote{This inequality is essentially a special case of the fractional Leibnitz rule, viewed on the Fourier transform side in a dualized form.}
$$ \langle \xi_4 \rangle^{5/4}
\lesssim \langle \xi_4 \rangle^{1/2}
\sum_{j=1}^3 \langle \xi_j \rangle^{3/4}$$
and symmetry to reduce to
$$
\|
\frac{\langle \xi_4 \rangle^{1/2} \langle \xi_2 \rangle^{1/2} }
{
\langle \xi_1 \rangle^{1/4} \langle \xi_3 \rangle^{1/4}
\langle \tau_4 - \xi_4^3 \rangle^{1/2 - \eps}
\prod_{j=1}^3
\langle \tau_j - \xi_j^3 \rangle^{1/2+\eps}
}
\|_{[4;\R \times \R]} \lesssim 1.
$$
We may minorize $\langle \tau_2 - \xi_2^3 \rangle^{1/2+\eps}$ by
$\langle \tau_2 - \xi_2^3 \rangle^{1/2-\eps}$.  But then the estimate follows from \eqref{kpv-2-est} and the $TT^*$ identity (Lemma \ref{compose}).
\end{proof}

We now give some examples in the periodic setting.  We begin with a proof of Bourgain's $L^4$ Strichartz estimate for the periodic KdV equation.

\begin{proposition}\label{borg-est}\cite{borg:hsd}  For any function $u$ in $\T \times \R$, we have
\be{borg-strich}
\| u \|_{L^4(\T \times \R)} \lesssim \| u \|_{X^{0,1/3}_{\tau = \xi^3}(\T \times \R)}.
\end{equation}
\end{proposition}

\begin{proof}
As this estimate is linear rather than multilinear we shall be able to apply some additional techniques such as Littlewood-Paley theory (and such mundane tools as the triangle inequality) in order to simplify the argument substantially (basically by preventing cross-terms when we finally pass to the bilinear setting).  

We first observe that if $u$ is constant in $x$, then this estimate follows immediately from the one-dimensional Sobolev embedding $H^{1/3}(\R) \subset L^4(\R)$.  Thus we may subtract off the mean and assume that $\hat u(0,\tau)$ is identically zero.  By dividing the spatial frequency into regions $\xi < 0$ and $\xi > 0$ and using symmetry of the cubic $\tau = \xi^3$ we may assume further that $\hat u$ is supported on the half-plane $\xi > 0$.  

For each dyadic $N \gtrsim 1$, let $u_N$ be a smooth localization to the frequency range $\xi \sim N$.  Since we have
$$ \| u  \|_{X^{0,1/3}_{\tau = \xi^3}(\T \times \R)} \sim (\sum_N \| u_N \|_{X^{0,1/3}_{\tau = \xi^3}(\T \times \R)}^2)^{1/2}$$
and the Littlewood-Paley estimate
$$ \| u \|_{L^4} \sim \| (\sum_N |u_N|^2)^{1/2} \|_4 \leq (\sum_N \|u_N\|_4^2)^{1/2}$$
it suffices to show that
$$ \| u_N \|_{L^4(\T \times \R)} \lesssim \| u_N \|_{X^{0,1/3}_{\tau = \xi^3}(\T \times \R)}.$$
uniformly in $N$.  (This trick works in general for all Strichartz estimates; see e.g. \cite{sogge:wave}).

Fix $N \gtrsim 1$.  Squaring the above estimate and using duality, we reduce to showing that
$$ |\int\int u_N^2 w\ dx dt| \lesssim \| u_N \|_{X^{0,1/3}_{\tau = \xi^3}(\T \times \R)}^2
\| w \|_{L^2(\T \times \R)}.$$
By Plancherel, this will follow if we can show that
$$ \| \frac{ \chi_{\xi_1, \xi_2 \sim N} }{
\langle \tau_1 - \xi_1^3 \rangle^{1/3}
\langle \tau_2 - \xi_2^3 \rangle^{1/3} } \|_{[3; \Z \times \R]} \lesssim 1.$$
Since $\xi_1, \xi_2$ are positive and comparable to $N$, we see that $\xi_3$ is negative and comparable to $N$.

We may assume \eqref{tj-large}. By a dyadic decomposition it thus suffices to show that
$$ \sum_{L_1,L_2,L_3 \gtrsim 1} \| \frac{ X_{N,N,N;N^3;L_1,L_2,L_3} }{
L_1^{1/3} L_2^{1/3} } \|_{[3; \Z \times \R]} \lesssim 1$$
We may assume $L_1 \leq L_2 \leq L_3$ as the other cases are similar or better.

First consider the contribution when $L_3 \sim N^3$.  In this case we use \eqref{kdv-excep} to estimate the above by
$$
\sum_{N^3 \gtrsim L_2 \gtrsim L_1} \frac{ L_1^{1/2} \langle N_1^{-1/4} L_2^{1/4} \rangle }{
L_1^{1/3} L_2^{1/3} },$$
and this easily sums to $O(1)$ as desired.

By \eqref{t-comp} it remains only to consider the case when $L_3 \sim L_2 \gg N^3$.  We then apply \eqref{kdv-standard} to reduce to
$$ \sum_{L_3 \gtrsim N^3} \sum_{L_1} \frac{ L_1^{1/2} \langle N^{-1} (N^3)^{1/2} \rangle }{
L_1^{1/3} L_3^{1/3} } \lesssim 1.$$
But this is easily verified also.
\end{proof}

From Proposition \ref{borg-est} we have

\begin{corollary}\cite{kpv:kdv} We have
\be{kpv-t-est}
\| (uv)_x \|_{X^{-1/2,-1/2}_{\tau = \xi^3}(\T \times \R)} 
\lesssim
\| u \|_{X^{-1/2,1/2}_{\tau = \xi^3}(\T \times \R)} 
\| v \|_{X^{-1/2,1/2}_{\tau = \xi^3}(\T \times \R)}
\end{equation}
whenever $u$ and $v$ satisfy the mean zero condition
$$ \int u(x,t)\ dx = \int v(x,t)\ dx = 0$$
for all $t$. 
\end{corollary}

This estimate is a key ingredient in the local well-posedness theory of the KdV equation in $H^s(\T)$ for $s \geq -1/2$, see \cite{kpv:kdv}. 
One can also obtain \eqref{kpv-t-est} directly from Proposition \ref{kdv-core} and summing, but in order to avoid logarithmic divergence problems one must either use some very delicate orthogonality arguments (including some non-trivial uses of Lemma \ref{schur}), or some inspired applications of Corollary \ref{conjugate}.  (In \cite{kpv:kdv}, one has to show that various logarithmic integrals are in fact convergent).  Here, we give a more direct proof based on Proposition \ref{borg-est}.

\begin{proof}
Applying duality and Plancherel as before, the estimate \eqref{kpv-t-est} is equivalent to
$$
\| \frac{(\xi_1 + \xi_2) \chi_{\xi_1 \xi_2 \neq 0} \langle \xi_1 \rangle^{1/2} \langle\xi_2\rangle^{1/2} \langle \xi_3 \rangle^{-1/2}}{
\langle \tau_1 - \xi_1^3 \rangle^{1/2}
\langle \tau_2 - \xi_2^3 \rangle^{1/2}
\langle \tau_3 - \xi_3^3 \rangle^{1/2} }
\|_{[3; \Z \times \R]} \lesssim 1.$$
We may re-arrange the numerator, and write this more symmetrically as
$$
\| \frac{\chi_{\xi_1 \xi_2 \xi_3 \neq 0} |\xi_1|^{1/2} |\xi_2|^{1/2} |\xi_3|^{1/2}}{
\langle \tau_1 - \xi_1^3 \rangle^{1/2}
\langle \tau_2 - \xi_2^3 \rangle^{1/2}
\langle \tau_3 - \xi_3^3 \rangle^{1/2} }
\|_{[3; \Z \times \R]} \lesssim 1.$$
From the resonance identity \eqref{kdv-fund} we have
$$ \sum_{j=1}^3 \tau_j - \xi_j^3 = - 3 \xi_1 \xi_2 \xi_3$$
and hence
$$ 1 \lesssim \sum_{j=1}^3 \frac{|\tau_j - \xi_j^3|^{1/2}}{|\xi_1|^{1/2} |\xi_2|^{1/2} |\xi_3|^{1/2}}.$$
Inserting this estimate into the above and using symmetry, we reduce to showing that
$$
\| \frac{1}{
\langle \tau_2 - \xi_2^3 \rangle^{1/2}
\langle \tau_3 - \xi_3^3 \rangle^{1/2} }
\|_{[3; \Z \times \R]} \lesssim 1.$$
We may of course replace $\langle \tau_j - \xi_j^3 \rangle^{1/2}$ by
$\langle \tau_j - \xi_j^3 \rangle^{1/3}$.  The above estimate is then equivalent to the trilinear estimate
$$ |\int\int uvw\ dx dt| \lesssim \| u\|_{L^2_{x,t}} \| v \|_{X^{0,1/3}_{\tau = \xi^3}} \| w \|_{X^{0,1/3}_{\tau = \xi^3}}.$$
But this is an immediate consequence of Proposition \ref{borg-est} and H\"older.
\end{proof}

To obtain sharp higher order multilinear estimates in the periodic setting, one would inevitably be led (e.g. by Lemma \ref{cauchy-schwarz}) to the number-theoretic problem of counting integer solutions to Diophantine equations such as $n = \xi_1^3 + \ldots + \xi_k^3$, $\xi_1 + \ldots + \xi_k = 0$; see e.g. the $L^6$ theory in \cite{borg:hsd}.  As we have seen, though, the $k=3$ theory is simpler (mainly thanks to the resonance identity \eqref{kdv-fund}), and does not require much number theory.

\section{Transverse intersections}\label{transverse-sec}

In the next two sections we develop some additional tools which will help us in the wave and Schr\"odinger computations.

Corollary \ref{char-2} allows us to estimate \eqref{h-gen} in terms of measures of certain sets.  If we applied these techniques to the wave equation relations $h(\xi) = \pm |\xi|$ then one would eventually be forced to compute the measures of neighbourhoods of ellipsoids and hyperboloids.  While this can be done (see e.g. \cite{damiano:null}), we prefer to give a slightly different treatment which relies only on the transversality of the surfaces $\tau_j = h_j(\xi_j)$ and not on specific geometric facts about ellipsoids and hyperboloids.  More precisely, we will use

\begin{figure}[htbp] \centering
\ \psfig{figure=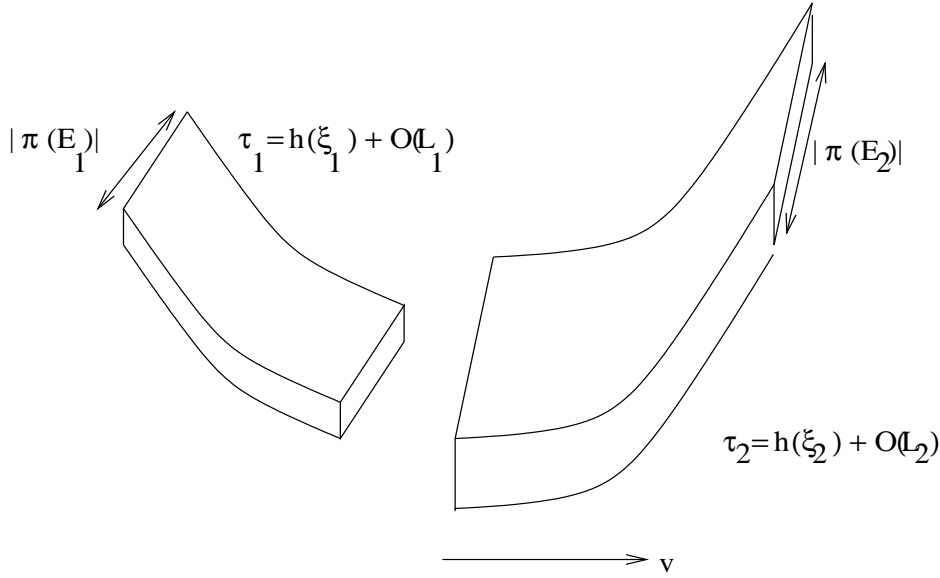,width=5in,height=3in}
\caption{The situation in Lemma \ref{transverse}.  Because the slopes of $h_1$ and $h_2$ in the $v$ direction differ by at least $\theta$, the generic intersection between the two regions will be something like a box with height $\min(L_1,L_2)$, width $\max(L_1,L_2)/\theta$, and with the remaining dimensions having measure $\min(|\pi_v(E_1)|, |\pi_v(E_2)|)$.}
\end{figure}

\begin{lemma}\label{transverse}
Let $E_1, E_2$ be open subsets of $\R^d$, and let $v$ be a unit vector in $\R^d$.  Let $\theta > 0$, and let $h_1: E_1 \to \R$, $h_2: E_2 \to \R$ be smooth functions which satisfy the transversality condition
\be{trans-cond}
|D_v h_1(\xi_1) - D_v h_2(\xi_2)| \gtrsim \theta
\end{equation}
for all $\xi_1 \in E_1, \xi_2 \in E_2$, where $D_v$ is the directional derivative in the direction $v$.  Then for any $L_1, L_2 > 0$ we have
\be{trans-main}
\| \prod_{j=1}^2 \chi_{E_j}(\xi_j) \chi_{|\lambda_j| \lesssim L_j}
\|_{[3,\R^d \times \R]}
\lesssim L_1^{1/2} L_2^{1/2} \theta^{-1/2} \min(|\pi_v(E_1)|, |\pi_v(E_2)|)^{1/2},
\end{equation}
where $|\pi_v(E)|$ is the $d-1$-dimensional measure of the projection of $E$ onto the orthogonal complement of $v$.  Also, we have the crude bound
\be{trans-triv}
\| \prod_{j=1}^2 \chi_{E_j}(\xi_j) \chi_{|\lambda_j| \lesssim L_j}
\|_{[3,\R^d \times \R]}
\lesssim \min(L_1,L_2)^{1/2} \min(|E_1|,|E_2|)^{1/2}.
\end{equation}
\end{lemma}

\begin{proof}
By Lemma \ref{char-1} we may estimate the left-hand side of \eqref{trans-main}, \eqref{trans-triv} by
$$ \min(L_1, L_2)^{1/2}
|\{ \xi_1 \in E_1: \xi - \xi_1 \in E_2; h_1(\xi_1) + h_2(\xi - \xi_1) = \tau + O(\max(L_1,L_2)) \}|^{1/2}$$
for some $\xi \in \R^d$, $\tau \in \R$.  The claim \eqref{trans-triv} is now clear.  To obtain \eqref{trans-main}, we observe from \eqref{trans-cond} that
$h_1(\xi_1) + h_2(\xi - \xi_1)$ has a derivative of $\gtrsim \theta$ in the $v$ direction of $\xi_1$.  In particular, $h_1(\xi_1) + h_2(\xi - \xi_1)$ is monotone in the $v$ direction, and for fixed values of $\pi_v(\xi_1)$ the $v$ co-ordinate of $\xi_1$ is constrained inside an interval of length $O(\theta^{-1} \max(L_1,L_2))$.  The claim \eqref{trans-main} (for $|\pi_v(E_1)|$) then follows by Fubini's theorem; the claim for $|\pi_v(E_2)|$ follows by symmetry.
\end{proof}

This lemma can also be used to give alternate proofs of much of Propositions \ref{kdv-core}, \ref{sppp-core}, and \ref{sppm-core}, although it has some difficulty dealing with the coherent interactions (because the transversality parameter $\theta$ then degenerates to zero).  Note that if one attempted to extend this lemma to the periodic setting in higher dimensions one would immediately encounter the difficult number-theoretic problem of accurately estimating the number of lattice points which lie near a prescribed hypersurface.

One could certainly consider other transversality conditions than \eqref{trans-cond}, such as a control on a mixed partial derivative (cf. \cite{ccw}).  We will not consider these matters here.  (In any event, the explicit nature of the functions $h_j(\xi)$ ensures that these quantities can always be computed accurately in the non-periodic setting).

\section{Separating the coarse and fine scales}\label{split-sec}

Another tool to analyze the norm $\|m\|_{[k;Z]}$ involves separating the fine scales from the coarse scales, in the spirit of \eqref{semi-direct}.  To make this precise we use the box covering notation from before, with $R$ containing the fine scales and $\Sigma$ describing the coarse scales\footnote{Here ``fine'' and ``coarse'' are in the context of frequency space.  In physical space of course the two notions are reversed.}.

\begin{lemma}\label{split}
Let $(R + \eta)_{\eta \in \Sigma}$ be a box covering of $Z$ with $|R| < \infty$, and let $m$ be a $[k;Z]$ multiplier.  Define the function $M: \Sigma^k \to \R^+$ by
\be{M-def}
M(\eta_1, \ldots, \eta_k) := \| |m(\xi)| \prod_{j=1}^k \chi_{R + \eta_j}(\xi_j) \|_{[k;Z]}.
\end{equation}
Then
\be{M-eq}
\| m \|_{[k;Z]} \lesssim |R|^{1 - \frac{k}{2}} \| \sum_{\eta \in \Sigma^k} M(\eta) 
\prod_{j=1}^k \chi_{kR + \eta_j}(\xi_j) \|_{[k;Z]}
\end{equation}
with the implicit constant dependent on $k$, and $kR$ being the sum of $k$ copies of $R$.
\end{lemma}

The $|R|^{1-\frac{k}{2}}$ is a natural scaling factor, which also appears in Lemma \ref{scaling}.  Heuristically, the right-hand side of \eqref{M-eq} is like $\| M \|_{[k;\Sigma]}$, and so \eqref{M-eq} can be viewed as a variant of \eqref{semi-direct} if one accepts that $Z$ is approximately isomorphic to
$\Sigma \times R$.

\begin{proof}
To prove the claim it suffices by \eqref{czk-def} to show that
\be{split-init}
|\int_{\Lambda_k(Z)} m(\xi) \prod_{j=1}^k f_j(\xi_j)|
\lesssim |R|^{1 - \frac{k}{2}} \| \tilde M \|_{[k;Z]}
\end{equation}
for all $L^2$-normalized $f_1, \ldots, f_k$, where $\tilde M(\xi)$ is the coarse scale multiplier
$$ \tilde M(\xi) := \sum_{\eta \in \Sigma^k} M(\eta) 
\prod_{j=1}^k \chi_{kR + \eta_j}(\xi_j).$$

Fix the $f_j$.  Decompose $f_j = \sum_{\eta_j \in \Sigma} f_{j,\eta_j}$, where
$f_{j,\eta_j} := f_j \chi_{R + \eta_j}$ is the restriction of $f_j$ to the box $R + \eta_j$.  By \eqref{czk-def} and \eqref{M-def}, we may estimate the left-hand side of \eqref{split-init} by
$$ \sum_{\eta \in \Sigma^k}
M(\eta) \prod_{j=1}^k \| f_{j,\eta_j} \|_2.
$$
We may estimate this by
$$ |R|^{1-k} \int_{\Lambda_k(Z)} \tilde M(\xi)
\prod_{j=1}^k F_j(\xi_j)
$$
where the functions $F_1, \ldots, F_k$ on $Z$ are defined by
$$ F_j(\xi_j) := \sum_{\eta_j \in \Sigma} \|f_{j,\eta_j}\|_2 \chi_{kR + \eta_j}(\xi_j).$$
The claim then follows since
$ \| F_j\|_2 \lesssim |R|^{1/2} \|f_j\|_2.$
\end{proof}

This lemma allows one to smooth out mild singularities in the multiplier.  A typical application is

\begin{corollary}\label{smooth}  Let $\R^d \times \R$ be parameterized by $(\xi,\tau)$ for $\xi \in \R^d$, $\tau \in \R$.
Let $m$ be a $[3;\R^d \times \R]$, and let $s_1, \ldots, s_3 \geq 0$ and $\alpha_1, \ldots, \alpha_3 \geq 0$ be such that $s_1 + \ldots + s_3 < d/2$ and $\alpha_1 + \ldots + \alpha_3 < 1/2$.  Then
$$ \| \prod_{j=1}^3 
(\frac{\langle \xi_j \rangle}{|\xi_j|})^{s_j} 
(\frac{\langle \tau_j \rangle}{|\tau_j|})^{\alpha_j} m \|_{[3,\R^d \times \R]}
\lesssim \| \tilde m \|_{[3, \R^d \times \R]}$$
where the implicit constant depends on $k$, the $s_j$, and the $\alpha_j$, and
$$ \tilde m(\xi,\tau) := \sup_{\xi' = \xi + O(1); \tau' = \tau + O(1)}
|m(\xi',\tau')|.$$
\end{corollary}

\begin{proof}
We apply the previous Lemma with $R$ being the unit cube in $\R^d \times \R$, and $\Sigma$ being the lattice $\Z^d \times \Z$.  We observe that for any $\eta_1, \eta_2, \eta_3 \in \Sigma$ we have
$$
\| \prod_{j=1}^3 
(\frac{\langle \xi_j \rangle}{|\xi_j|})^{s_j} 
(\frac{\langle \tau_j \rangle}{|\tau_j|})^{\alpha_j} \chi_{R + \eta_j}(\xi_j) \|_{[3,\R^d \times \R]} \lesssim 1$$
thanks to the Comparison principle, Lemma \ref{tensor}, and Corollary \ref{sobolev-mult}.  The claim then follows from Lemma \ref{split}.
\end{proof}

One can also use this lemma to perform induction on scale arguments, and thus obtain some sharp multiplier estimates which do not seem to be accessible by more elementary techniques; see Section \ref{scale-sec}.  

\section{Estimates related to the wave equation}\label{wave-sec}

We now consider estimates of the form \eqref{mess} for the wave equation in $\R^d$ for\footnote{The one-dimensional case essentially reduces to product Sobolev theory in null co-ordinates, see e.g. \cite{keel:wavemap}, and the estimates can be deduced from Lemma \ref{tensor} and Proposition \ref{sobolev}.  We omit the details.} $d \geq 2$.  The situation for the wave is slightly different than the first order equations discussed earlier in that there are two dispersion relations, $\tau = |\xi|$ and $\tau = -|\xi|$.  In practice, this means that one must modify the definition of $X^{s,b}$ norm slightly\footnote{In some of the wave equation literature one uses the slightly different weight $\langle |\xi| + |\tau| \rangle$ instead of $\langle |\xi| \rangle$.  Although that weight is more natural from Lorentz invariance considerations, it is technically more complicated, and has some difficulty dealing with Coulomb gauge conditions.  In the main, though, the norms are equally capable of proving well-posedness results.} to
$$ \| u \|_{X^{s,b}_{\tau = \pm |\xi|}} :=
\| \langle \xi \rangle^s \langle |\tau| - |\xi| \rangle^b \hat u \|_{L^2_{\xi,\tau}}.$$
The problem of obtaining good bilinear estimates in $X^{s,b}$ then reduces to that of controlling expressions of the form\footnote{From the homogeneous bilinear estimates in \cite{damiano:null} and such arguments as Proposition \ref{average}, Corollary \ref{smooth}, and Corollary \ref{convex} one can already obtain a large class of estimates of the form \eqref{m-bound-wave} which seem to be mostly adequate for applications.  See \cite{kl-sel:survey}.}
\be{m-bound-wave}
\| \frac{m(\xi_1,\xi_2,\xi_3)}{\prod_{j=1}^3 \langle |\tau_j| - |\xi_j| \rangle^{b_j}} \|_{[3; \R^d \times \R]}.
\end{equation}

The dyadic decomposition of this expression into building blocks like \eqref{h-gen} proceeds slightly differently from before.  For any $N_1,N_2,N_3,L_1,L_2,L_3,H > 0$, $\epsilon_1,\epsilon_2,\epsilon_3 = \pm 1$, we consider the quantities
\be{w-gen}
\| X_{N_1,N_2,N_3;H;L_1,L_2,L_3;\epsilon_1,\epsilon_2,\epsilon_3}
\|_{[3; \R^d \times \R]}
\end{equation}
where 
$$ X_{N_1,N_2,N_3;H;L_1,L_2,L_3;\epsilon_1,\epsilon_2,\epsilon_3}
:= 
\chi_{|h(\xi)| \sim H}
\prod_{j=1}^3 \chi_{|\xi_j| \sim N_j; \epsilon_j \tau_j \geq 0; |\lambda_j| \sim L_j }.$$
and $\lambda_j$, $h(\xi)$ are defined as before with $h_j(\xi) := \epsilon_j |\xi|$.

If all the signs $\eps_j$ agree, then \eqref{w-gen} vanishes since $\tau_1 + \tau_2 + \tau_3 = 0$.  Thus it suffices by symmetry and time reversal (using e.g. Lemma \ref{scaling}) to consider the $(++-)$ case
\be{eps-set}
\epsilon_1 = \epsilon_2 = +1; \quad \epsilon_3 = -1,
\end{equation}
and we shall assume this for the remainder of this discussion.

The resonance function $h(\xi): \Gamma_3(\R^d) \to \R$ is now given by
$$ h(\xi_1,\xi_2,\xi_3) := |\xi_1| + |\xi_2| - |\xi_3|.$$
From the algebraic identity
$$ h(\xi_1,\xi_2,\xi_3) = 2 \frac{|\xi_1| |\xi_2| - \xi_1 \cdot \xi_2}{|\xi_1| + |\xi_2| + |\xi_3|}$$
we obtain the resonance ``identity''
\be{h-angle}
|h| \sim \frac{N_1 N_2}{N_{max}} \angle(\xi_1, \xi_2)^2 \sim \min(N_1,N_2) \angle(\xi_1, \xi_2)^2
\end{equation}
when $|\xi_j| \sim N_j$, and 
\be{cosine}
\angle(\xi_1, \xi_2) := \cos^{-1} \frac{ \xi_1 \cdot \xi_2 }{|\xi_1| |\xi_2|}
\end{equation}
is the angle between $\xi_1$ and $\xi_2$.  Thus one has resonance when $\xi_1$ and $\xi_2$ are close to parallel.  Because of this, our multilinear expressions shall be analyzed by angular decompositions of the $\xi_1$ and $\xi_2$ variables.

As before, the expression \eqref{w-gen} vanishes unless \eqref{n-comp} and \eqref{t-comp} both hold.  Also, from \eqref{h-angle} we must have
\be{h-w}
H \lesssim \min(N_1,N_2).
\end{equation}

\begin{figure}[htbp] \centering
\ \psfig{figure=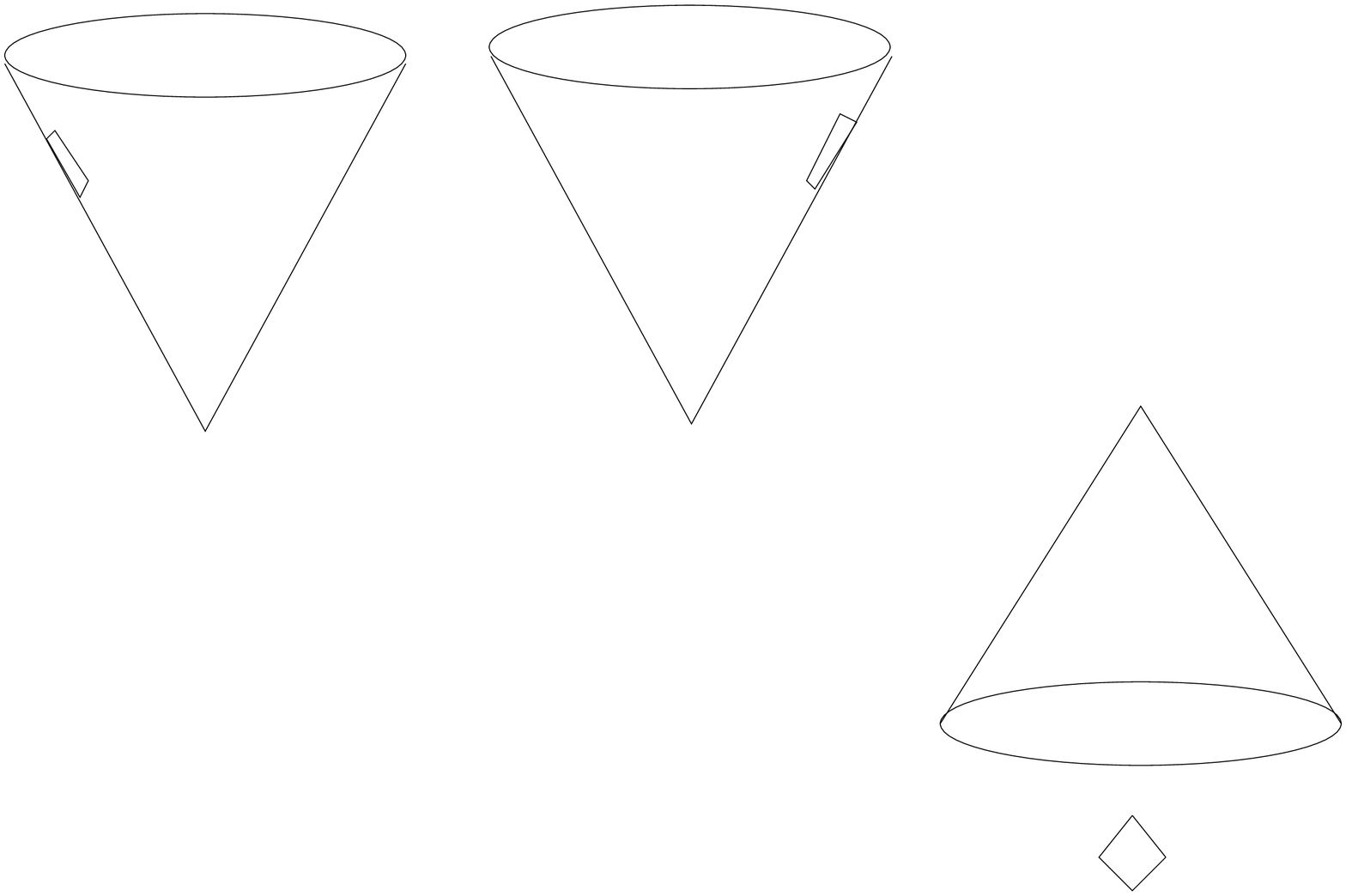,width=5in,height=3in}
\caption{The example which shows that \eqref{wpp} is sharp in the (++) high-high interaction case.}
\end{figure}

\begin{proposition}\label{wave-core}
Let $N_1, N_2, N_3, L_1, L_2, L_3, H > 0$ obey \eqref{n-comp}, \eqref{t-comp}, \eqref{h-w}.  
\begin{itemize}
\item ($(++)$ high-high interactions) If $N_1 \sim N_2 \gg N_3$, then
\be{wpp}
\eqref{w-gen} \lesssim L_{min}^{1/2} \min(L_{med},N_3)^{1/2} N_3^{(d-1)/2}.
\end{equation}
Also, \eqref{w-gen} vanishes unless $H \sim N_1$, and if $L_3 \ll N_1$, then \eqref{w-gen} vanishes unless $L_1, L_2 \sim N_1$. 
\item  (High-low interactions) If $N_1 \sim N_3 \gtrsim N_2$ and $L_2 \ll L_{max}$, then
\be{w-standard}
\eqref{w-gen} \lesssim (\frac{H}{N_2})^{(d-3)/4} L_{min}^{1/2} \min(H,L_{med})^{1/2} N_2^{(d-1)/2}.
\end{equation}
\item ($(+-)$ high-high interactions) If $N_1 \sim N_3 \gtrsim N_2$ and $L_2 \sim L_{max}$, then
\be{w-pm}
\eqref{w-gen} \lesssim (\frac{H}{N_2})^{(d-3)/4} L_{min}^{1/2} \min(H, \frac{N_1}{N_2} L_{med})^{1/2} N_2^{(d-1)/2}.
\end{equation}
\end{itemize}
Similar statements hold with the role of $1$ and $2$ reversed.
\end{proposition}

\begin{figure}[htbp] \centering
\ \psfig{figure=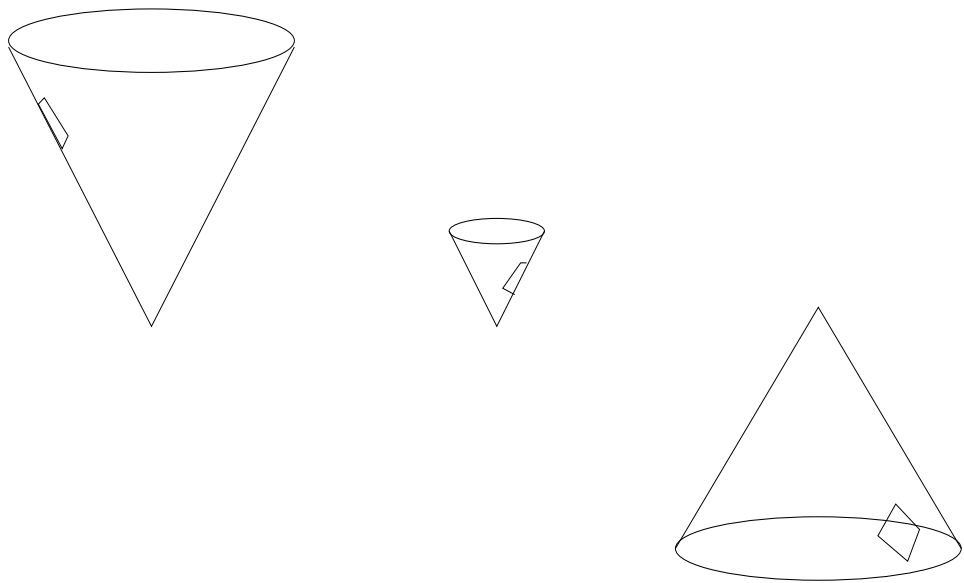,width=5in,height=3in}
\caption{The example which shows that \eqref{w-standard} is sharp in the high-low interaction case (with angular separation $\sim 1$, so that $H \sim N_2$).  One can use Lorentz transforms to essentially generate the full family of counterexamples to cover the small angle interaction case, when $H \ll N_2$.}
\end{figure}

The exponent $(d-3)/4$ is a familiar aspect of wave equation estimates, and can be related to Lorentz invariance considerations.  Variants of these estimates appear elsewhere in the literature, for instance estimates for \eqref{w-gen} for the Klein-Gordon equation with large mass appear in \cite{delort:fang}, and estimates for homogeneous solutions to the wave equation (which essentially corresponds to the case when $L_1, L_2 \ll N_{min}, L_3$) appear in \cite{damiano:null}.

\begin{figure}[htbp] \centering
\ \psfig{figure=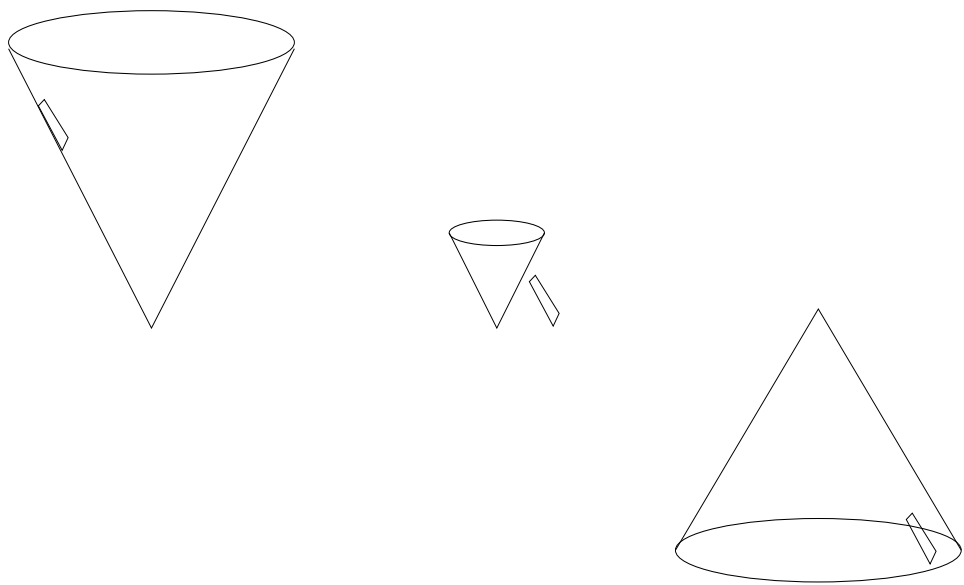,width=5in,height=3in}
\caption{The example which shows that \eqref{w-pm} is sharp in the (+-) high-high interaction case.}
\end{figure}

\begin{proof}
We treat the three cases of the Proposition separately.

{\bf Case 1: ($(++)$ high-high interactions) $N_1 \sim N_2 \gg N_3$.}

The constraint $H \sim N_1$ is clear, since $|\xi_1| + |\xi_2| - |\xi_3| \sim N_1 + O(N_3) \sim N_1$.  To verify the second constraint, suppose $L_3 \ll N_1$.  Then $|\tau_3| \ll N_1$.  Since $\tau_1 + \tau_2 = -\tau_3$ and $\tau_1, \tau_2$ are positive, we must have $|\tau_1|, |\tau_2| \ll N_1$, hence $L_1, L_2 \sim N_1$.  This shows the second constraint.

Now we prove \eqref{wpp}.  If \eqref{t-dom} holds then the claim follows from \eqref{t-dom-2}.  In light of this, \eqref{t-comp}, and the preceding discussion, the only remaining possibility is that $L_3 \sim H \sim N_1 \gtrsim L_1, L_2$.  By symmetry we may assume that $L_2 \geq L_1$.

In this case $(\xi_1, \tau_1)$ and $(\xi_2,\tau_2)$ are near opposite sides of the upper light cone, while $(\xi_3, \tau_3)$ is near the time axis.  The intersection of one light cone with a translated inverse of the other will basically be a graph over an ellipsoid of bounded eccentricity.

The variable $\xi_3$ can be localized to a ball of radius $\ll N_3$, so one can localize $\xi_1$ and $\xi_2$ to similar balls by Lemma \ref{boxes}.  We thus have
$$
\eqref{w-gen} \lesssim
\| \prod_{j=1}^3 \chi_{|\xi_j - \xi_j^0| \ll N_3} \chi_{|\lambda_j| \lesssim L_j} \|_{[3; \R^d \times \R]},
$$
where $\xi_1^0, \xi_2^0, \xi_3^0$ are such that $|\xi_j^0| \sim N_j$.  We may of course assume that $|\xi_1^0 + \xi_2^0 + \xi_3^0| \ll N_3$.

By Lemma \ref{char-1} we obtain
$$
\eqref{w-gen} \lesssim L_1^{1/2}
| \{ \xi_1: |\xi_1 - \xi_1^0| \ll N_3; |\xi_1| + |\xi - \xi_1| = \tau + O(L_2) \}|^{1/2}$$
Since $|\xi_1| \sim N_1$ and $|\xi| \sim N_3 \ll N_1$, we that the above set is within $O(L_2/N_1)$ of an ellipsoid of bounded eccentricity and principal radii $\sim N_1$.  Since the set is also contained in a ball of radius $\ll N_3$, it must have measure $O(\min(N_3, L_2/N_1) N_3^{d-1})$, and the claim follows.

{\bf Case 2: (High-low interactions) $N_1 \sim N_3 \gtrsim N_2$; $L_2 \ll L_{max}$.}

We either have $L_3 \gtrsim L_2, L_1$ or $L_1 \gtrsim L_2, L_3$.

{\bf Case 2(a): ($(++)$ case) $L_3 \gtrsim L_2, L_1$.}

By repeating the Case 1 argument we have
$$
\eqref{w-gen} \lesssim
\| \chi_{|h(\xi)| \sim H} \prod_{j=1}^2 \chi_{|\xi_j - \xi_j^0| \ll N_2} \chi_{|\lambda_j| \lesssim L_j} \|_{[3; \R^d \times \R]}$$
where $\xi_1^0, \xi_2^0, \xi_3^0$ are such that $|\xi_j^0| \sim N_j$.  We may of course assume that $|\xi_1^0 + \xi_2^0 + \xi_3^0| \ll N_2$.

From \eqref{h-angle} we have $\angle(\xi_1, \xi_2) \sim (H/N_2)^{1/2}$. This relation motivates the following partition.  Let $C$ be a large number, let $\Omega$ be a maximal $C^{-1} (H/N_2)^{1/2}$-separated subset of the sphere $S^{d-1}$, and for each $\omega \in \Omega$ let $\Gamma(\omega)$ denote the cone
\be{gam-def}
\Gamma(\omega) := \{ \xi \in \R^d \backslash \{0\}: \angle(\xi, \omega) \leq C^{-1/2} (H/N_2)^{1/2} \}.
\end{equation}
We then partition the $\xi_1$ and $\xi_2$ variable into these cones:
$$ 1 \lesssim \sum_{\omega_1} \sum_{\omega_2} \chi_{\Gamma(\omega_1)}(\xi_1) \chi_{\Gamma(\omega_2)}(\xi_2).$$
The contribution of a single term vanishes unless $\angle(\omega_1, \omega_2) \sim (H/N_2)^{1/2}$.  Thus each cone $\Gamma(\omega_1)$ interacts with at most $O(1)$ cones $\Gamma(\omega_2)$, and conversely.  By \eqref{sharp-schur} we thus have
$$
\eqref{w-gen} \lesssim
\| \prod_{j=1}^2 \chi_{|\xi_j - \xi_j^0| \ll N_2} \chi_{\Gamma(\omega_j)}(\xi_j) \chi_{|\lambda_j| \lesssim L_j} \|_{[3; \R^d \times \R]}$$
for one such pair $\Gamma(\omega_1), \Gamma(\omega_2)$ of cones.

Fix $\omega_1$, $\omega_2$.  Let $v$ be a unit vector which is coplanar with $\omega_1$, $\omega_2$, is perpendicular with $\omega_1$, and obeys $\angle(v,\omega_2) \sim 1$. (This latter condition is redundant unless $H \sim N_2$).  Observe that the conditions of Lemma \ref{transverse} apply with 
$$ E_j := \{ \xi: |\xi - \xi_j^0| \ll N_2; \xi \in \Gamma(\omega_j) \},$$
$h_1(\xi) = h_2(\xi) = |\xi|$, and $\theta = (H/N_2)^{1/2}$.  Since $|\pi_v(E_2)| \lesssim N_2 (HN_2)^{(d-2)/2}$, we thus have from \eqref{trans-main} that
$$
\eqref{w-gen} \lesssim L_1^{1/2} L_2^{1/2} N_2^{(d-1)/2}
(H/N_2)^{(d-3)/4}.$$
Also, from \eqref{trans-triv} and the fact that $|E_2| \lesssim N_2 (HN_2)^{(d-1)/2}$ we have
$$
\eqref{w-gen} \lesssim L_{min}^{1/2} N_2^{d/2} (H/N_2)^{(d-1)/4}.$$
Combining the two estimates we obtain the claim \eqref{w-standard}.

{\bf Case 2(b): ($(+-)$ interactions) $L_3 \gtrsim L_2, L_1$.}

We modify the Case 2(a) argument as follows.  We begin with
$$
\eqref{w-gen} \lesssim
\| \chi_{|h(\xi)| \sim H} \prod_{j=2}^3 \chi_{|\xi_j - \xi_j^0| \ll N_2} \chi_{|\lambda_j| \lesssim L_j} \|_{[3; \R^d \times \R]}$$
where the $\xi_j^0$ are as before.  

\begin{figure}[htbp]\label{fig:sine} \centering
\ \psfig{figure=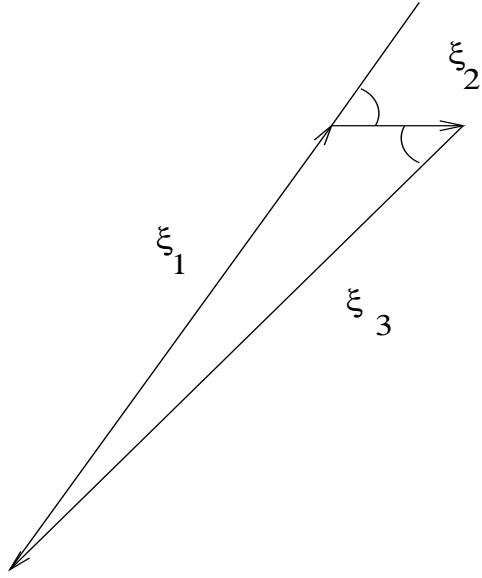,width=2.5in,height=3in}
\caption{If $|\xi_1| \sim |\xi_3|$, then the two displayed angles are comparable in size.  The angle subtended by $\xi_1$ and $\xi_3$ is about $|\xi_1|/|\xi_3|$ the size of the other two displayed angles.}
\end{figure}

Since $\angle(\xi_1, \xi_2) \sim (H/N_2)^{1/2}$ and $|\xi_1| \sim |\xi_3|$, the sine rule (see Figure \ref{fig:sine}) shows that $\angle(-\xi_3, \xi_2) \sim (H/N_2)^{1/2}$.  We then introduce the cones $\Gamma(\omega)$ as before and partition
$$ 1 \lesssim \sum_{\omega_2} \sum_{\omega_3} 
\chi_{\Gamma(\omega_2)}(\xi_2) \chi_{\Gamma(\omega_3)}(\xi_3).$$
The contribution of a single term vanishes unless $\angle(\omega_2, -\omega_3) \sim (H/N_2)^{1/2}$.  By Schur's test \eqref{sharp-schur} we reduce to
$$
\eqref{w-gen} \lesssim
\| \prod_{j=2}^3 \chi_{|\xi_j - \xi_j^0| \ll N_2} \chi_{\Gamma(\omega_j)}(\xi_j) \chi_{|\lambda_j| \lesssim L_j} \|_{[3; \R^d \times \R]}.$$
One then applies Lemma \ref{transverse} as before; the function $h_3(\xi)$ is $-|\xi|$ instead of $+|\xi|$, and we have $\angle(\omega_2, -\omega_3) \sim (H/N_2)^{1/2}$ rather than $\angle(\omega_1, \omega_2) \sim (H/N_2)^{1/2}$, but the two changes essentially cancel each other out, and we obtain the same bounds. We omit the details. 

{\bf Case 3: ($(+-)$ high-high interactions) $N_1 \sim N_3 \gtrsim N_2$; $L_2 \sim L_{max}$.}

We begin with
$$
\eqref{w-gen} \lesssim
\| \chi_{|h(\xi)| \sim H} \prod_{j=1,3} \chi_{|\xi_j - \xi_j^0| \ll N_2} \chi_{|\lambda_j| \lesssim L_j} \|_{[3; \R^d \times \R]}$$
where the $\xi_j^0$ are as before.  

Since $\angle(\xi_1,\xi_2) \sim (H/N_2)^{1/2}$ and $|\xi_2| \sim \frac{N_2}{N_1} |\xi_3|$, the sine rule shows that $\angle(-\xi_1, \xi_3) \sim \frac{N_2}{N_1} (H/N_2)^{1/2}$.    We can then introduce cones as before, except with angular width $C^{-1} \frac{N_2}{N_1} (H/N_2)^{1/2}$ rather than $C^{-1} (H/N_2)^{1/2}$.  We then apply Lemma \eqref{transverse} with $h_1(\xi) = h_3(\xi) = |\xi|$ and $\theta = (N_2/N_1) (H/N_2)^{1/2}$.  The claim \eqref{w-pm} follows.
\end{proof}

The bounds listed above are sharp.  For \eqref{wpp}, this is obtained by testing \eqref{czk-def} with 
$$ f_1(\xi,\tau) := \chi_{|\xi - N_1 e_1| \lesssim N_3} \chi_{|\lambda| \lesssim L_1}$$
$$ f_2(\xi,\tau) := \chi_{|\xi + N_1 e_1| \lesssim N_3} \chi_{|\lambda| \lesssim L_2}$$
$$ f_3(\xi,\tau) := \chi_{|\xi| \lesssim N_3} \chi_{|\lambda| \lesssim L_3} \chi_{|\tau + 2N_1| \lesssim L_2 + L_1 + N_3}.$$
For \eqref{w-standard}, the example is given by
$$ f_j(\xi,\tau) := \chi_{|\xi| \sim N_j; \angle(\xi,\eps_j e_1) \lesssim (H/N_2)^{1/2}} \chi_{|\lambda| \lesssim L_j}.$$
This also shows that \eqref{w-pm} is sharp when $L_{med} \gtrsim \frac{N_2}{N_1} H$.  When $L_{med} \lesssim \frac{N_2}{N_1}$, one uses
$$ f_j(\xi,\tau) := \chi_{|\xi| \sim N_j; \angle(\xi,\eps_j e_1) \lesssim \frac{N_2}{N_1} (H/N_2)^{1/2}} \chi_{|\lambda| \lesssim L_j}$$
for $j=1,3$, and
$$ f_2(\xi,\tau) := \chi_{|\xi| \sim N_2; \angle(\xi,e_1) \lesssim (H/N_2)^{1/2}} \chi_{\tau = \xi \cdot e_1 + O(\frac{N_2}{N_1}H)}$$
instead.  We omit the details.

One can use Proposition \ref{wave-core} to prove various null form estimates.  We will not give a comprehensive list of null form estimates (cf. the program initiated in \cite{damiano:null}), but content ourselves with some examples.  We first consider the null form $Q_{ij}$ defined for $1 \leq i < j \leq d$ by
$$ Q_{ij}(\phi, \psi) = \partial_i \phi \partial_j \psi - \partial_j \phi \partial_i \psi.$$

\begin{proposition}\label{qij}
Let $d \geq 3$, and let $D^{-1}$ be a Fourier multiplier on $\R^d$ which is a homogeneous symbol of order $-1$.  Then for all $0 < \eps \ll 1$ and $s := \frac{d}{2}-\frac{3}{4}+C\eps$ we have
$$ \| D^{-1} Q_{ij}(\phi,\psi) \|_{X^{s-1,-1/4+\eps}_{\tau = \pm |\xi|}(\R^d \times \R)} 
\lesssim
\| \phi \|_{X^{s,3/4+\eps}_{\tau = \pm |\xi|}(\R^d \times \R)} 
\| \psi \|_{X^{s,3/4+\eps}_{\tau = \pm |\xi|}(\R^d \times \R)}$$
and
$$ \| Q_{ij}(D^{-1} \phi,\psi) \|_{X^{s-1,-1/4+\eps}_{\tau = \pm |\xi|}(\R^d \times \R)} 
\lesssim
\| \phi \|_{X^{s,3/4+\eps}_{\tau = \pm |\xi|}(\R^d \times \R)} 
\| \psi \|_{X^{s,3/4+\eps}_{\tau = \pm |\xi|}(\R^d \times \R)}$$
\end{proposition}

Morally speaking, this Proposition suggests\footnote{The above estimates are not quite sufficient by themselves, because one also needs some ``elliptic'' estimates to deal with the gauge condition; see e.g. \cite{kl-mac:yang}.  Also, these well-posedness results are not optimal.  In $d=3,4$ one can improve down to the scaling regularity $\frac{d}{2}-1+\eps$ in $d=3,4$ \cite{kl-mac:null2}, \cite{kl-mac:null3}, \cite{kl-tar:yang-mills}, and presumably one can do so for $d \geq 5$ by the ideas in \cite{tat:5+1}.  However, the $s > \frac{d}{2} - \frac{3}{4}$ result appears to be the best one can achieve from $X^{s,b}$ norms alone.} that the Maxwell-Klein-Gordon and Yang-Mills equations in (say) the Coulomb gauge are locally well-posed in $H^s$ for $s > \frac{d}{2} - \frac{3}{4}$.  In $d=3$ this computation was carried through for the full Maxwell-Klein-Gordon equations in Coulomb gauge in \cite{cuccagna}; see also \cite{keel:mkg}, \cite{tao:yang-mills}.

\begin{proof}
The symbols for $D^{-1} Q_{ij}(\phi,\psi)$ and $Q_{ij}(D^{-1} \phi,\psi)$ can both be majorized by
$$ \frac{ |\xi_1 \wedge \xi_2| }{\min(|\xi_1|, |\xi_2|, |\xi_3|)}.$$
Thus it suffices to show that
$$
\| \frac{ |\xi_1 \wedge \xi_2| \langle \xi_3 \rangle^{s-1} }
{ \min(|\xi_1|, |\xi_2|, |\xi_3|)
\langle \xi_1 \rangle^s \langle \xi_2 \rangle^s 
\langle |\tau_1| - |\xi_1| \rangle^{3/4+\eps}
\langle |\tau_2| - |\xi_2| \rangle^{3/4+\eps}
\langle |\tau_3| - |\xi_3| \rangle^{1/4-\eps} } \|_{[3;\R^d \times \R]} \lesssim 1.
$$

By Corollary \ref{smooth} (splitting into three cases depending on which of the $|\xi_j|$ is smallest) and the hypothesis $d \geq 3$, we may replace the
term $\min(|\xi_1|, |\xi_2|, |\xi_3|)$ by $\min(\langle \xi_1 \rangle, \langle \xi_2\rangle, \langle \xi_3\rangle)$.  By two applications of Proposition \ref{average} (and splitting into regions $\tau_j > 0$ and $\tau_j < 0$)
one can replace $\frac{1}{\langle \tau_j - \xi_j \rangle^{3/4+\eps}}$ by $\chi_{||\tau_j| - |\xi_j|| \sim 1}$ for $j=1,2$).  We have thus reduced to
$$
\| \frac{ |\xi_1 \wedge \xi_2| \langle \xi_3 \rangle^{s-1} 
\prod_{j=1}^2 \chi_{||\tau_j| - |\xi_j|| \sim 1} }
{ \min(\langle \xi_1 \rangle, \langle \xi_2\rangle, \langle\xi_3 \rangle)
\langle \xi_1 \rangle^s \langle \xi_2 \rangle^s 
\langle |\tau_3| - |\xi_3| \rangle^{1/4-\eps} } \|_{[3;\R^d \times \R]} \lesssim 1.
$$
We may assume \eqref{ximax-large} and \eqref{tj-large} as usual.  We can break this up dyadically as
$$
\| \sum_{ N_{max} \gtrsim 1}
\sum_{L_3 \gtrsim 1} \sum_H
\frac{ |\xi_1 \wedge \xi_2| \langle N_3 \rangle^{s-1} }
{ \langle N_{min} \rangle
\langle N_1 \rangle^s \langle N_2 \rangle^s 
L_3^{1/4-\eps} } 
X_{N_1,N_2,N_3;H;1,1,L_3;\epsilon_1,\epsilon_2,\epsilon_3}
 \|_{[3; \R^d \times \R]} \lesssim 1.$$
We may assume that two of the $\epsilon_j$ is positive and the third is negative.  We may assume \eqref{n-comp} and \eqref{t-comp}, so in particular $L_3 \sim \max(1,H)$.  

Suppose for the moment that we are in the situation \eqref{eps-set}.  Then from the estimate
$$ |\xi_1 \wedge \xi_2| \lesssim N_1 N_2 \angle(\xi_1,\xi_2)$$
and \eqref{h-angle} we have
$$ |\xi_1 \wedge \xi_2| \lesssim N_1 N_2 \frac{H^{1/2}}{\min(N_1,N_2)^{1/2}}
\lesssim N_{max} N_{min}^{1/2} H^{1/2}.$$
Also, from the permutation invariance of $|\xi_1 \wedge \xi_2|$ we have
$$ |\xi_1 \wedge \xi_2| \lesssim N_{max} N_{min}.$$
By symmetry we thus have
$$|\xi_1 \wedge \xi_2| \lesssim N_{max} N_{min}^{1/2} \min(H,N_{min})^{1/2}$$ for all choices of $\epsilon_j$.  We thus reduce to
\be{wave-schur}
\| \sum_{ N_{max} \sim N_{med} \gtrsim 1}
\sum_H
\frac{ N_{max} \min(H,N_{min})^{1/2} \langle N_3 \rangle^{s-1} }
{ \langle N_{min} \rangle^{1/2}
\langle N_1 \rangle^s \langle N_2 \rangle^s 
\langle H \rangle^{1/4-\eps} } 
X_{N_1,N_2,N_3;H;1,1,\max(1,H);\epsilon_1,\epsilon_2,\epsilon_3}
 \|_{[3; \R^d \times \R]} \lesssim 1.
\end{equation}
By Schur's test we may assume that $N_{max} \sim N_{med} \sim N$ for some fixed $N \gtrsim 1$.

Since
$$ \min(H, N_{min})^{1/2} \leq H^{1/4-2\eps} N_{min}^{1/4+2\eps}
\leq \langle H \rangle^{1/4 - \eps} \langle N_{min} \rangle^{1/4+2\eps} \min(H^\eps, H^{-\eps})$$ 
and
$$
\frac{ \langle N_3 \rangle^{s-1} }
{ \langle N_1 \rangle^s \langle N_2 \rangle^s 
} \lesssim N^{-1} \langle N_{min} \rangle^{-s}
$$
we have
$$
\frac{ N_{max} \min(H,N_{min})^{1/2} \langle N_3 \rangle^{s-1} }
{ \langle N_{min} \rangle^{1/2}
\langle N_1 \rangle^s \langle N_2 \rangle^s 
\langle H \rangle^{1/4-\eps} } 
\lesssim \langle N_{min} \rangle^{-s-\frac{1}{4}+2\eps}
\min(H^\eps, H^{-\eps}).$$
Also, from the various cases of Proposition \ref{wave-core} and the hypothesis $d \geq 3$ we have
$$
\| X_{N_1,N_2,N_3;H;1,1,\max(1,H);\epsilon_1,\epsilon_2,\epsilon_3} \|_{[3; \R^d \times \R]} \lesssim N_{min}^{(d-1)/2}.$$
From these estimates and the hypothesis $s > \frac{d}{2} - \frac{3}{4} + C\eps$ we easily see that \eqref{wave-schur} holds uniformly in $N$.
\end{proof}

One could also treat the null form
$$ Q_0(\phi,\psi) = \phi_t \psi_t - \nabla \phi \cdot \nabla \psi$$
by these techniques.  The symbol for this null form is essentially given by
$$ \tau_1 \tau_2 - \xi_1 \cdot \xi_2$$
and this can be estimated by, e.g.
$$ |\tau_1 \tau_2 - \xi_1 \cdot \xi_2| \lesssim L_{max} (N_{max} + L_{max}).$$
In principle, one can use this to recover estimates such as those in \cite{kman.selberg}, but we shall not do so here.

One can also prove product estimates of the form
$$ \| \phi\psi \|_{X^{s,b}_{\tau = \pm |\xi|}(\R^d \times \R)} 
\lesssim
\| \phi \|_{X^{s_1,b_1}_{\tau = \pm |\xi|}(\R^d \times \R)} 
\| \psi \|_{X^{s_2,b_2}_{\tau = \pm |\xi|}(\R^d \times \R)}$$
with $b_1, b_2$ large (say $b_1, b_2 > 1/2$) by the above techniques.
Because of the $(\frac{H}{N_2})^{(d-3)/4}$ factor in Proposition \ref{wave-core}, one can only get efficient estimates when $b \geq -\frac{d-3}{4}$.  When $b > -\frac{d-3}{4}$ there is no difficulty in deriving these types of estimates from Proposition \ref{wave-core}, but one can run into difficulties involving logarithmic divergences at the endpoint $b = -\frac{d-3}{4}$.  In some ``double endpoint'' cases one cannot remove the logarithm entirely (see \cite{damiano:null}, or \cite{ntt:kdv-counter} for some analogues for KdV and Schr\"odinger), but in less extreme cases one can use Schur's test \eqref{sharp-schur} to eliminate the divergence.  We illustrate this\footnote{This estimate is proven in several places, e.g. \cite{damiano:null}, \cite{tataru:wave1}.} with a simplified model estimate, namely
\be{endpoint}
\| 
\chi_{|\tau_1 - |\xi_1||, |\tau_2 - |\xi_2|| \sim 1}
\chi_{|\xi_1|, |\xi_2|, |\xi_3| \sim N} \langle \tau_3 + |\xi_3| \rangle^{\frac{d-3}{4}}
\|_{[3;\R^d \times \R]} \lesssim 
N^{(3d-5)/4}
\end{equation}
for fixed $N \gg 1$.  A routine dyadic decomposition using \eqref{w-standard} would yield $N^{(3d-5)/4} \log(N)$.  To eliminate this logarithm, we dyadically decompose $|\tau_3 - |\xi_3||$ (ignoring the easily dealt with region $|\tau_3 - |\xi_3|| \lesssim 1$) as
$$
\| 
\sum_{1 \ll L_3 \lesssim N}
X_{N,N,N;L_3;1,1,L_3;+1,+1,-1}
L_3^{\frac{d-3}{4}}
\|_{[3;\R^d \times \R]} \lesssim 
N^{(3d-5)/4}.
$$
Now for each $L_3$, we let $\Omega(L_3)$ be a maximal $C^{-1} (L_3/N)^{1/2}$-separated subset of the sphere $S^{d-1}$, and for each $\omega \in \Omega(L_3)$ let $\Gamma_{L_3}(\omega)$ be the cones with direction $\omega$ and angular width $C^{-1/2} (L_3/N)^{1/2}$ as in \eqref{gam-def}.  We can thus split the above as
$$
\| 
\sum_{1 \ll L_3 \lesssim N}
\sum_{\omega, \omega' \in \Omega(L_3)} m_{L_3,\omega,\omega'}
\|_{[3;\R^d \times \R]} \lesssim 
N^{(3d-5)/4}
$$
where
$$
m_{L_3,\omega,\omega'} :=
X_{N,N,N;L_3;1,1,L_3;+1,+1,-1}
\chi_{\Gamma_{L_3}(\omega)}(\xi_1) \chi_{\Gamma_{L_3}(\omega')}(\xi_2)
L_3^{\frac{d-3}{4}}.
$$
From \eqref{t-comp} and \eqref{h-angle} we see that $m_{L_3,\omega,\omega'}$ vanishes unless
$N^{-1/2} + \angle(\omega,\omega') \sim (L_3/N)^{1/2}$, and we may implicitly assume this in the summation.  

Using the notation of Lemma \ref{schur}, we have
$$ \supp_{1,2}(m_{L_3,\omega,\omega'}) \subset (\Gamma_{L_3}(\omega) \times \R) \times (\Gamma_{L_3}(\omega') \times \R).$$
Also, we have
$$ \supp_3(m_{L_3,\omega,\omega'}) \subset 
(\tilde \Gamma_{L_3}(\omega) \times \R) \cap \{ (\xi,\tau): ||\tau| - |\xi|| \sim L_3 \}$$
where $\tilde \Gamma_{L_3}(\omega)$ is some enlargement of $\Gamma_{L_3}(\omega)$.  Applying Lemma \ref{schur} with $J_1 = \{1,2\}$ and $J_3 = \{3\}$, it thus suffices to show that
$$ \| m_{L_3, \omega, \omega'} \|_{[3;\R^d \times \R]} \lesssim N^{(3d-5)/4}$$
for each $L_3$, $\omega$, $\omega'$.  But this easily follows by discarding the constraints $\xi_1 \in \Gamma_{L_3}(\omega)$, $\xi_2 \in \Gamma_{L_3}(\omega')$ and using \eqref{w-standard}.

\section{Orthogonal interactions}\label{scale-sec}

We present the following estimate controlling the interaction of nearly orthogonal frequencies, which will be of importance in the study of higher-dimensional Schr\"odinger estimates.  The technique of induction of scales which is used in the proof may also be of application to other expressions $\|m(\xi)\|_{[k;\R^d]}$ in which $m$ concentrates near a smooth hypersurface in $\Gamma_k(\R^d)$.

\begin{proposition}\label{rotate-prop}
Let $d \geq 2$, $R \geq r > 0$, and $0 < \theta \lesssim 1$.  Then
\be{rotate-eq}
\| \chi_{|\xi_1| \sim R} \chi_{|\xi_2| \sim r} \chi_{\angle(\xi_1,\xi_2) = \frac{\pi}{2} + O(\theta)} 
\|_{[3;\R^d]} \lesssim 
r^{d/2} \theta^{1/2} \min(1, \frac{R\theta}{r})^{\frac{1}{2}-\eps}
\end{equation}
for any $\eps > 0$, with the implicit constant depending on $\eps$.  When $d=2$ the $\eps$ can be removed.
\end{proposition}

Apart from the $\eps$, this estimate is sharp, as one can see by testing \eqref{czk-def} with 
$$ f_1(\xi) = \chi_{|\xi - Re_1| \lesssim r}; \quad f_2(\xi) = \chi_{|\xi| \lesssim r; |\xi \cdot e_1| \lesssim \theta + \frac{r}{R}}; f_3(\xi) = \chi_{|\xi + Re_1| \lesssim r}.$$

\begin{proof}  We recommend reading this proof initially assuming the simplifying assumption $r=R$, as this case already contains the main idea of the argument.  The reader may also wish to rescale $r=R=1$.

First suppose that $\theta \gtrsim r/R$.  Then we estimate the left-hand side of \eqref{rotate-eq} by
$$
\| \chi_{|\xi_2| \sim r} \chi_{\angle(\xi_1, \xi_2) = \frac{\pi}{2} + O(\theta)} \|_{[3;Z]}.$$
The claim then follows by Lemma \ref{cauchy-schwarz} with $j=1$.  Note that this argument also gives the $\theta \ll r/R$ case with $\eps = 1/2$.

Henceforth we will assume $\theta \ll r/R$.

We now give the more elementary $d=2$ argument, in which $\eps = 0$.  We remark that this argument is essentially in \cite{cdks}.

\begin{figure}[htbp] \centering
\ \psfig{figure=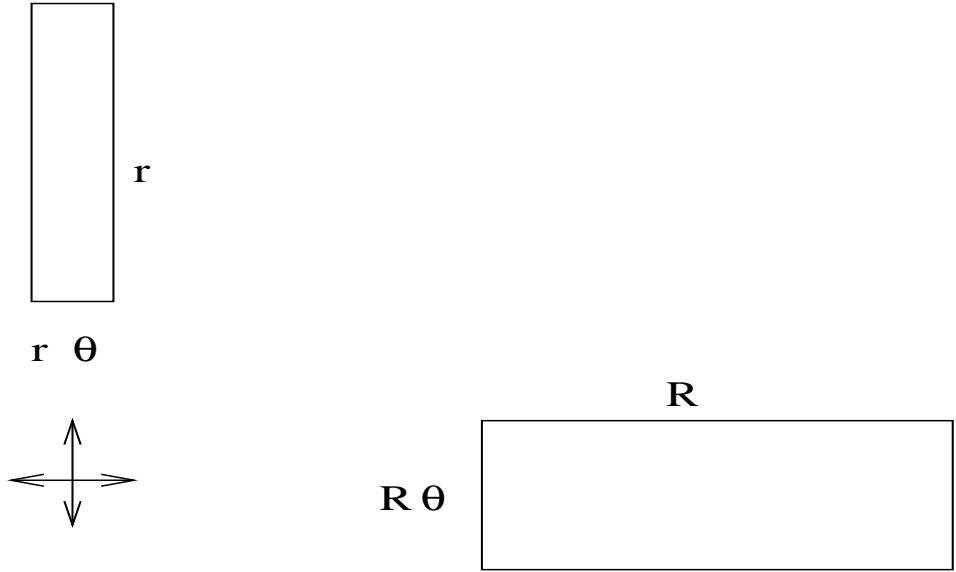,width=5in,height=3in}
\caption{The $d=2$ situation.  One uses Lemma \ref{schur} to localize $\xi_1$ and $\xi_2$ to the sectors displayed.  The generic intersection between these two regions is a $r\theta \times R\theta$ rectangle.}
\end{figure}

Divide $R^2$ into sectors $\Gamma$ of angular width $\sim \theta$ for some large constant $C$, with vertex at the origin.  We can estimate \eqref{rotate-eq} as
$$ \| \sum_{\Gamma, \Gamma'}
\chi_{|\xi_1| \sim R} \chi_{|\xi_2| \sim r} \chi_{\Gamma}(\xi_1)
\chi_{\Gamma'}(\xi_2) \chi_{\angle(\xi_1,\xi_2) = \frac{\pi}{2} + O(\theta)} \|_{[3;Z]}.$$
The summands vanish unless $\angle(\Gamma, \Gamma') \geq \frac{\pi}{2} - O(\theta)$.  Thus each $\Gamma$ interacts with at most $O(1)$ sectors $\Gamma'$, and conversely.  By Schur's test (Lemma \ref{schur}) it thus suffices to show that
$$
\| \chi_{|\xi_1| \sim R} \chi_{|\xi_2| \sim r} \chi_{\Gamma}(\xi_1)
\chi_{\Gamma'}(\xi_2) \|_{[3;Z]}
\lesssim
r^{1/2} R^{1/2} \theta.$$
But this follows from \eqref{char-zar} and elementary geometry.

We now give the general $d \geq 2$ argument.  We shall assume that \eqref{rotate-eq} has already been proven for some $\eps$, and show that this implies \eqref{rotate-eq} with $\eps$ replaced by $\eps/2$.  Since we have already proven \eqref{rotate-eq} for $\eps = 1/2$, the claim then follows by iteration.

Roughly speaking, the point is as follows.  The condition $\angle(\xi_1,\xi_2) \geq \frac{\pi}{2} - \theta$ has too much curvature in the $d>2$ case to be usefully decomposed into boxes or sectors.  However, at scales $\sqrt{rR\theta}$ and below, the condition no longer depends on $d-2$ of the dimensions, and we can use Lemma \ref{tensor} to reduce to the $d=2$ case. One then uses Lemma \ref{split} and the iteration hypothesis \eqref{rotate-eq} to handle the coarse scales, thus introducing the $\eps/2$ loss.

We turn to the details.  Fix $R$, $r$, $\theta$, $\eps$.  We may localize $\xi_2$ to a ball of the form $|\xi_2 - \xi_2^0| \ll r$ for some $|\xi_2^0| \sim r$.  By Lemma \ref{boxes} we can thus localize $\xi_1$ to a similar ball
$|\xi_1 - \xi_1^0| \ll r$.  By a mild rotation and scaling we may assume that $\xi_1^0 = Re_1$. We thus have to show that
\be{rotate-targ}
\| m(\xi_1,\xi_2) \|_{[3;\R^d]} \lesssim 
r^{d/2} \theta^{1/2} (\frac{R\theta}{r})^{\frac{1-\eps}{2}}
\end{equation}
where
$$ m(\xi_1,\xi_2) := \chi_{|\xi_1 - \xi_1^0| \ll r} \chi_{|\xi_2 - \xi_2^0| \ll r} \chi_{\angle(\xi_1, \xi_2) = \pi/2 + O(\theta) }.$$
Let $Q$ denote the box\footnote{The reason for this choice of $Q$ is that the angular condition $\angle(\xi_1, \xi_2) = \pi/2 + O(\theta)$ becomes two-dimensional on translates of $Q$, and that $Q$ is essentially maximal with respect to this property.}  centered at the origin with sides parallel to the axes, and all side-lengths equal to $\sqrt{rR\theta}$ except for the $e_1$ side-length, which is $\frac{r}{R}\sqrt{rR\theta}$.  Let $\Sigma$ be the canonical tiling lattice of $Q$, so that $(Q + \eta)_{\eta \in \Sigma}$ is thus a box covering.

\begin{figure}[htbp] \centering
\ \psfig{figure=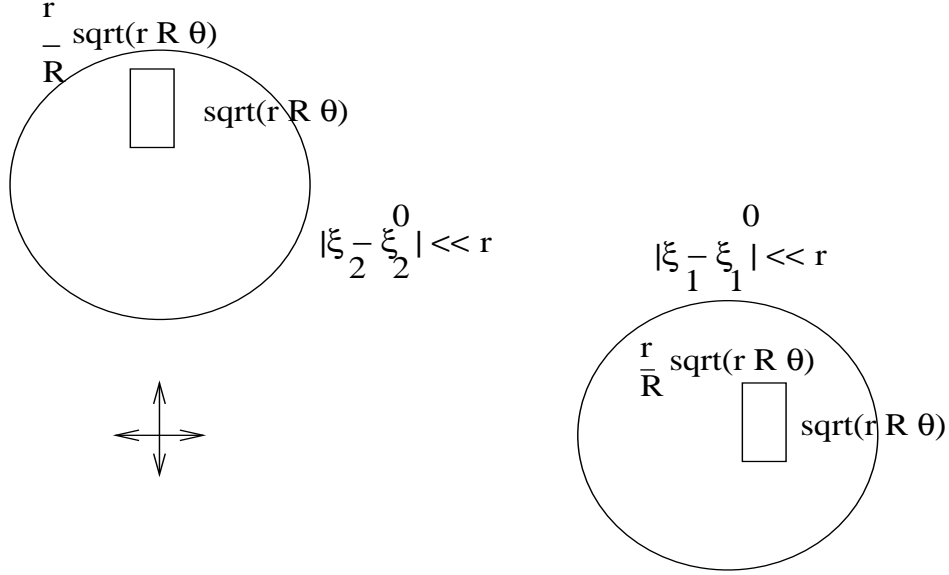,width=5in,height=3in}
\caption{The $d>2$ situation.  One uses box localization to localize $\xi_1$, $\xi_2$ to the balls displayed.  One then subdivides space into boxes of the displayed dimension.  Note that two boxes will only interact if they subtend an angle of $\pi/2 + O(\sqrt{r\theta/R})$ at the origin.  Inside each box the angular condition $\angle(\xi_1,\xi_2) = \pi/2 + O(\theta)$ becomes essentially two-dimensional.}
\end{figure}

We wish to apply Lemma \ref{split}.  We begin by proving the fine-scale estimate
\be{m-local}
\| m(\xi_1,\xi_2) \prod_{j=1}^3 \chi_{Q + \eta_j}(\xi_j) \|_{[3;\R^d]}
\lesssim r^{\frac{d}{4}} R^{\frac{d}{4}} \theta^{\frac{d+2}{4}}
\end{equation}
for all $\eta_1,\eta_2,\eta_3 \in \Sigma$.  We shall achieve this by reducing to the $d=2$ case already proven; alternatively one can modify the $d=2$ argument to prove this estimate directly.

Fix $\eta_1$, $\eta_2$, $\eta_3$.  In order for the left-hand side of \eqref{m-local} to be non-zero, there must exist $\xi_1 \in Q + \eta_1$ and $\xi_2 \in Q + \eta_2$ such that $m(\xi_1,\xi_2) \neq 0$.  From this and elementary geometry we see that
\be{eta-cond}
|\eta_1| \sim R; \quad |\eta_2| \sim r; \quad \angle(\eta_1, \eta_2) = \frac{\pi}{2}
+ O(\sqrt{\frac{r\theta}{R}}).
\end{equation}

For any $\xi \in \R^d$, write $\xi = (\xi',\xi'')$, where $\xi'$ is the orthogonal projection to $span(\eta_1,\eta_2)$ and $\xi''$ is the complement of this projection.  We observe that when $\xi_j \in Q + \eta_j$ for $j=1,2$, we have the estimates
$$|\xi'_1| \sim R; \quad |\xi'_2| \sim r; \quad |\xi''_1| \lesssim \sqrt{Rr\theta}; |\xi''_2| \lesssim \sqrt{Rr\theta}$$
which implies that
$$ \xi_1 \cdot \xi_2 = \xi'_1 \cdot \xi'_2 + \xi''_1 \cdot \xi''_2
= \xi'_1 \cdot \xi'_2 + O(Rr\theta).$$
From this and \eqref{cosine} we see that the condition $\angle(\xi_1,\xi_2) = \pi/2 + O(\theta)$ implies that $\angle(\xi'_1,\xi'_2) = \pi/2 + O(\theta)$.  In other words, the components $\xi''_1$ and $\xi''_2$ do not significantly affect the angle.  We can thus estimate the left-hand side of \eqref{m-local} by
$$
\|  \chi_{\angle(\xi'_1, \xi'_2) = \pi/2 + O(\theta) } \chi_{|\xi'_1| \sim R} \chi_{|\xi'_1| \sim r} \prod_{j=1}^3 \chi_{Q + \eta_j}(\xi_j) \|_{[3;\R^d]}.
$$
Recall that the variables $\xi''_1$, $\xi''_2$ are constrained to a ball of radius $\sqrt{rR\theta}$.  Applying \eqref{direct} and \eqref{box-est-2} we can therefore estimate the previous by
$$
(\sqrt{rR\theta})^{(d-2)/2}
\|  \chi_{\angle(\xi'_1, \xi'_2) = \pi/2 + O(\theta) } \chi_{|\xi'_1| \sim R} \chi_{|\xi'_2| \sim r} \|_{[3;\R^2]},
$$
and \eqref{m-local} then follows from the $d=2$ case of \eqref{rotate-eq} with $\eps = 0$, which we have proven previously.

From \eqref{m-local} and Lemma \ref{split} we have
$$
\| m(\xi_1,\xi_2) \|_{[3;\R^d]} \lesssim |Q|^{-\frac{1}{2}} 
\| \sum_{\eta \in \Sigma^k} M(\eta) \prod_{j=1}^k \chi_{kQ + \eta_j}(\xi_j) \|_{[k;Z]}$$
where 
$$M(\eta) := r^{\frac{d}{4}} R^{\frac{d}{4}} \theta^{\frac{d+2}{4}}$$
whenever \eqref{eta-cond} holds, and $M(\eta) = 0$ otherwise.  Applying the hypothesis \eqref{rotate-prop} with $\theta$ replaced by $\sqrt{r\theta/R}$ (the geometric mean of $\theta$ and $r/R$) we thus obtain
$$
\| m(\xi_1,\xi_2) \|_{[3;\R^d]} \lesssim |Q|^{-\frac{1}{2}} 
r^{\frac{d}{4}} R^{\frac{d}{4}} \theta^{\frac{d+2}{4}}
r^{d/2} (\sqrt{r\theta/R})^{1/2} (\frac{R\sqrt{r\theta/R}}{r})^{\frac{1}{2}-\eps}$$
which simplifies to \eqref{rotate-targ} with $\eps$ replaced by $\eps/2$, as desired. 
\end{proof}

The author conjectures that the $\eps$ can be removed in all dimensions, possibly by exploiting the $L^p$ theory of Radon transforms and related objects (see e.g. \cite{stein:radon}).  This may be related to the arguments in \cite{tvv:bilinear}, where the $L^p$ theory of circular averages was used to produce bilinear restriction estimates for the paraboloid.

\section{Estimates related to the Schr\"odinger equation}\label{schro-sec}

We now study \eqref{m-bound} assuming the Schr\"odinger dispersion relation $h_j(\xi_j) = \pm |\xi_j|^2$.  We will restrict our attention to the non-periodic case, as a sharp treatment of the periodic case\footnote{The one-dimensional periodic case, however, is similar to the periodic KdV situation, and can be dealt with by the techniques in Section \ref{kdv-sec}.  The two-dimensional semi-periodic case is also tractable, see \cite{tak-tz:rxt}.} seems to require some non-trivial number-theoretic information regarding the representations of a number as sums of squares. 

The intersections of one paraboloid with another are always graphs over spheres or hyperplanes.  As such, the geometry is simple enough that one can compute the measures of these intersections easily without recourse to Lemma \ref{transverse}, though one could of course use that Lemma to obtain the same estimates.  

Up to symmetry, there are only two possibilities for the $h_j$: the $(+++)$ case
\be{sppp}
h_1(\xi) = h_2(\xi) = h_3(\xi) = |\xi|^2
\end{equation}
and the $(++-)$ case
\be{sppm}
h_1(\xi) = h_2(\xi) = |\xi|^2; \quad h_3(\xi) = -|\xi|^2.
\end{equation}
The $(+++)$ case corresponds to estimates of the form
\be{sppp-est}
\| \overline{\phi} \overline{\psi} \|_{X^{s,b}_{\tau = |\xi|^2}}
\lesssim \| \phi \|_{X^{s_1,b_1}_{\tau = |\xi|^2}} \| \psi \|_{X^{s_2,b_2}_{\tau = |\xi|^2}}
\end{equation}
while the $(++-)$ case and its permutations are similar but treat $\phi \psi$, $\overline{\phi} \psi$, or $\phi \overline{\psi}$ instead of $\overline{\phi} \overline{\psi}$.

Of the two cases, the $(+++)$ case is substantially easier, because the resonance function 
$$h(\xi) := |\xi_1|^2 + |\xi_2|^2 + |\xi_3|^2$$ 
does not vanish except at the origin.  We will be able to obtain sharp estimates in this case in all dimensions.

The $(++-)$ case is more delicate, because the resonance function 
$$h(\xi) := |\xi_1|^2 + |\xi_2|^2 - |\xi_3|^2$$ 
can vanish when $\xi_1$ and $\xi_2$ are orthogonal, and we will need to call upon Proposition \ref{rotate-prop} to obtain sharp results (although we lose an epsilon for $d>2$).

We now consider the $(+++)$ case \eqref{sppp}.  As before, we begin with the treatment of \eqref{h-gen}.  Clearly we may assume that
\be{h-comp-sppp}
H \sim N_{max}^2
\end{equation}
since the symbol vanishes otherwise.  

\begin{proposition}\label{sppp-core}
Let $H, N_1, N_2, N_3, L_1, L_2, L_3 > 0$ obey \eqref{n-comp}, \eqref{t-comp}, \eqref{h-comp-sppp}.  Let $Z = \R^d$ for any $d \geq 1$, and let the dispersion relations be given by \eqref{sppm}.
\begin{itemize}
\item If $d \geq 2$, then we have
\be{sppp-standard}
\eqref{h-gen} \lesssim 
L_{min}^{1/2}
N_{max}^{-1/2}
N_{min}^{(d-1)/2}
\min(N_{max} N_{min}, L_{med})^{1/2}.
\end{equation}
\item If $d = 1$, then we also have \eqref{sppp-standard} except in the exceptional case $N_{max} \sim N_{min}$ and $L_{max} \sim H$, in which case
\be{sppp-excep}
\eqref{h-gen} \lesssim 
L_{min}^{1/2} L_{med}^{1/4}.
\end{equation}
\end{itemize}
\end{proposition}

The case \eqref{sppp-excep} corresponds to the coherent case $\xi_i = \xi_j$.  For $d \geq 2$ this case no longer dominates; heuristically, this states that two Schr\"odinger waves will generically not be coherent when the dimension is large.  A similar phenomenon holds for the wave equation, although the effect is shifted upwards by one dimension due to the fact that solutions to the wave equation must propagate along null rays.

\begin{proof}
This will be a reprise of the proof of Proposition \ref{kdv-core}.  (Alternatively, one could proceed using Lemma \ref{transverse}).

The case \eqref{t-dom} follows from \eqref{t-dom-2}, so we may assume that \eqref{h-dom} holds.  We may also assume that $L_1 \geq L_2 \geq L_3$, hence $L_1 \sim N_{max}^2$.  By Corollary \ref{char-2} we have
$$
\eqref{h-gen}
\lesssim
L_3^{1/2} | \{ \xi_2 \in Z: |\xi_2 - \xi_2^0| \ll N_{min}; |\xi_2|^2 + |\xi - \xi_2|^2 = \tau + O(L_2) \}|^{1/2}$$
for some $\xi$, $\tau$, $\xi_1^0$, $\xi_2^0$, $\xi_3^0$ satisfying \eqref{xio-cond} and $|\xi + \xi_1^0| \ll N_{min}$.

\begin{figure}[htbp] \centering
\ \psfig{figure=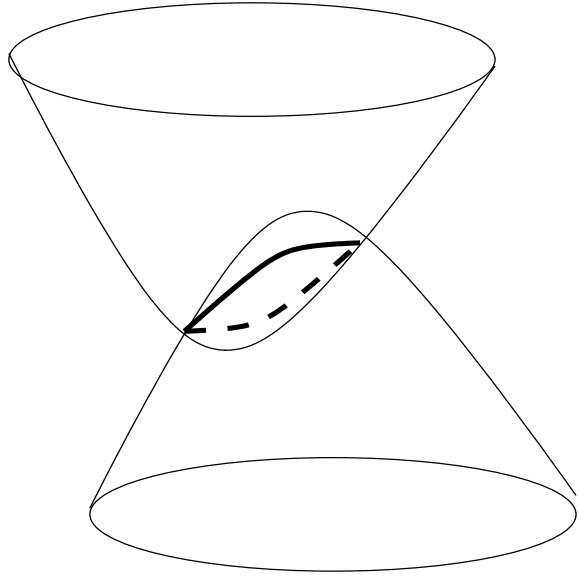,width=3in,height=3in}
\caption{The geometric interpretation of \eqref{ppp-ident}.  If one inverts an upward paraboloid and then translates it by $(\xi,\tau)$, the resulting intersection with the original upward paraboloid is a graph over a sphere centered at $\xi/2$ and radius $\sqrt{\tau/2 - |\xi|^2/4}$.}
\end{figure}

From Corollary \ref{char-2} and the identity
\be{ppp-ident}
|\xi_2 - \frac{\xi}{2}|^2 = \frac{|\xi_2|^2 + |\xi - \xi_2|^2}{2} - \frac{|\xi|^2}{4}
\end{equation}
it suffices to show that
\be{ppp-alt}
| \{ \xi_2: |\xi_2 - \xi_2^0| \ll N_{min}; 
|\xi_2 - \frac{\xi}{2}|^2 = \frac{\tau}{2} - \frac{|\xi|^2}{4} + O(L_2) \} |
\lesssim N_{max}^{-1}
N_{min}^{d-1} \min(L_2, N_{min} N_{max}),
\end{equation}
with the right-hand side replaced by $L_2^{1/2}$ in the exceptional case $d=1$, $N_{min} \sim N_{max}$.

As in Proposition \ref{kdv-core}, we need only consider three cases: $N_1 \sim N_2 \sim N_3$, $N_1 \sim N_2 \gg N_3$, and $N_2 \sim N_3 \gg N_1$.  

Suppose $N_1 \sim N_2 \sim N_3$.  Define the radius $R > 0$ by 
$$ R^2 := |\frac{\tau}{2} - \frac{|\xi|^2}{4}| + L_2.$$
If $R \gtrsim L_2^{1/2}$, then the set in \eqref{ppp-alt} is contained in an annulus of radius $R$ and thickness $O(L_2/R)$.  If $R \sim L_2^{1/2}$, then the above set is contained in a ball of radius $O(R)$.  In either case the claim follows (checking the $d=1$ and $d>1$ cases separately), noting that we must have $R = O(N_1)$ and $L_2 = O(N_1^2)$.

Now suppose $N_1 \sim N_2 \gg N_3$.  We can assume that $\frac{\tau}{2} - \frac{|\xi|^2}{4} \sim N_1^2$, since the set in \eqref{ppp-alt} vanishes otherwise.  But then this set is contained in an annulus of thickness $O(L_2/N_1)$ while simultaneously being contained in a ball of radius $O(N_3)$, and the claim \eqref{sppp-standard} follows.

Now suppose $N_2 \sim N_3 \gg N_1$.
We can assume that $\frac{\tau}{2} - \frac{|\xi|^2}{4} \sim N_2^2$.  But then this set is contained in an annulus of thickness $O(L_2/N_2)$ and simultaneously in a ball of radius $O(N_1)$, and the claim \eqref{sppp-standard} follows.
\end{proof}

The examples used to show that \eqref{kdv-standard} and \eqref{kdv-excep} were sharp can be easily adapted to show that the estimates in Proposition \ref{sppp-core} are similarly sharp. 

We now consider the $(++-)$ case \eqref{sppm}.  In one dimension $d=1$, Proposition \ref{sppp-core} holds unchanged in most cases except with \eqref{sppp} replaced by \eqref{sppm} and \eqref{h-comp-sppp} replaced by
$$ H \sim N_1 N_2,$$
with the additional change that in the cases $N_1 \sim N_{min}, L_1 \sim L_{max} \sim H$ or
$N_2 \sim N_{min}, L_2 \sim L_{max} \sim H$, the bound \eqref{sppp-standard} must be weakened\footnote{We thank
Sebastian Herr and Martin Hadac for pointing out this fact, which was erroneously omitted from previous versions of this paper.} to
$$ \eqref{h-gen} \lesssim L_{min}^{1/2} N_{min}^{-1/2} L_{med}^{1/2}.$$
We omit the routine modifications of the argument necessary to obtain this bound.

Now consider the $d \geq 2$ $(++-)$ case.  In contrast to the previous cases, the quantity $H$ can now take on a non-trivial range of values even when $N_1$, $N_2$, $N_3$ are fixed.  This is indicated by the resonance identity
$$
|h(\xi)| = ||\xi_1|^2 + |\xi_2|^2 - |\xi_3|^2| = 2 |\xi_1 \cdot \xi_2| \sim |\xi_1| |\xi_2| |\pi/2 - \angle(\xi_1, \xi_2)|,
$$
In particular, we may assume
\be{sppm-h}
H \lesssim N_1 N_2,
\end{equation}
and that
\be{sppm-ident}
\angle(\xi_1, \xi_2) = \frac{\pi}{2} + O(\frac{H}{N_1 N_2}).
\end{equation}
Thus the quantity $H$ imposes an orthogonality constraint on $\xi_1$ and $\xi_2$ similar to that which appears in Proposition \ref{rotate-prop}.

Another complication in the $(++-)$ case is that we do not have perfect symmetry between the three frequency variables; indeed, only $(\xi_1,\tau_1)$ and $(\xi_2,\tau_2)$ are interchangeable.  This causes an unpleasant increase in the number of cases to consider.

\begin{proposition}\label{sppm-core}
Let $N_1, N_2, N_3 > 0$, $L_1, L_2, L_3 > 0$, $H > 0$ satisfy \eqref{n-comp}, \eqref{t-comp}, and \eqref{sppm-h}.  Let $Z = \R^d$ for some $d \geq 2$, and let the dispersion relations be given by \eqref{sppm}.  Let $\eps > 0$.
\begin{itemize}
\item ((++) case)  If $N_1 \sim N_2 \gg N_3$, then \eqref{h-gen} vanishes unless $H \sim N_1^2$, in which case one has
\be{sppm-pp}
\eqref{h-gen} \lesssim L_{min}^{1/2} N_{max}^{-1/2} N_{min}^{(d-1)/2}
\min(N_{max} N_{min}, L_{med})^{1/2}.
\end{equation}
\item ((+-) Coherence) If we have
\be{sppm-excep-case}
N_1 \sim N_3 \gg N_2; H \sim L_2 \gg L_1, L_3, N_2^2 
\end{equation} 
then we have
\be{sppm-excep}
\eqref{h-gen} \lesssim L_{min}^{1/2} N_{max}^{-1/2} N_{min}^{(d-1)/2}
\min(H,\frac{H}{N_{min}^2}L_{med})^{1/2}
\end{equation}
Similarly with the roles of 1 and 2 reversed.
\item  In all other cases, we have
\be{sppm-est}
\eqref{h-gen} \lesssim L_{min}^{1/2} N_{max}^{-1/2} N_{min}^{(d-1)/2}
\min(H,L_{med})^{1/2}
\min(1, \frac{H}{N_{min}^2})^{1/2 - \eps}
\end{equation}
\end{itemize}
The implicit constants depend on $\eps$.  When $d=2$ the $\eps$ can be removed.
\end{proposition}

The estimates in this Proposition may appear overly complicated; nevertheless, they are sharp except for the $\eps$.  To give the examples we let $\{j_1, j_2, j_3\} = \{1,2,3\}$ be such that $L_{j_1} \geq L_{j_2} \geq L_{j_3}$.  When $N_1 \sim N_2 \gg N_3$ the counter-example is given by defining 
$\prod_{j=1}^3 f_j(\xi_j,\tau_j)$ to be the characteristic function on the region
$$ \xi_1 = N_{max} e_1 + O(N_{min}); \xi_2 = -N_{max} e_1 + O(N_{min}); \xi_3 = O(N_{min})$$
$$ \lambda_{j_1} = -N^2 + O(L_{med} + N_{max} N_{min}); 
\lambda_{j_2} = O(L_{med}); \lambda_{j_3} = O(L_{min}).$$
The left-hand side of \eqref{czk-def} is then $\gtrsim N_{min}^{2d} L_{min} L_{med}$, while the right-hand side is $\lesssim \eqref{h-gen} N_{min}^{3d/2} L_{med}^{1/2} L_{min}^{1/2} (L_{med}+N_{max} N_{min})^{1/2}$, which gives the example for \eqref{sppm-pp}.

When $N_1 \sim N_3 \gtrsim N_2$ the counter-example to \eqref{sppm-est} is given using the region
$$ \xi_1 = N_{max} e_1 + O(N_{min}); \xi_2 = O(N_{min}); \xi_2 \cdot e_1 = O(\frac{H + N_{min}^2}{N_{max}}); \xi_3 = -N_{max} e_1 + O(N_{min})$$
$$ \lambda_{j} = O(L_j) \hbox{ for } j=1,2,3.$$ 
The left-hand side of \eqref{czk-def} is then $\gtrsim N_{min}^{2d} \frac{H}{N_{max} N_{min}} \min(1, \frac{H}{N_{min}^2}) L_{min} L_{med}$, while the right-hand side is $\lesssim \eqref{h-gen} N_{min}^{3d/2} (\frac{H+N_{min}^2}{N_{min} N_{max}})^{1/2} L_{min}^{1/2} L_{med}^{1/2} (L_{med}+H)^{1/2}$, and the claim follows.

In the exceptional case \eqref{sppm-excep-case} one can improve this counterexample to
$$ \xi_1 = N_1 e_1 + O(N_2); \xi_2 = O(N_2); \xi_2 \cdot e_1 = O(\frac{H}{N_1}); \xi_3 = -N_1 e_1 + O(N_2)$$
$$ \lambda_1 = O(L_1); 
\lambda_2 = -2N_1 \xi_2 \cdot e_1 + O(N_2^2 + L_{2}); 
\lambda_3 = O(L_3).$$
The left-hand side of \eqref{czk-def} is then $\gtrsim N_2^{2d-1} \frac{H}{N_1} L_{min} L_{med}$, while the right-hand side is $\lesssim \eqref{h-gen} N_2^d N_2^{(d-1)/2} (\frac{H}{N_1})^{1/2} L_{min}^{1/2} L_{med}^{1/2} (N_2^2 + L_{med})^{1/2}$.  This shows that \eqref{sppm-excep} is sharp.

\begin{proof}
In addition to \eqref{ppp-ident}, we will need the algebraic identity 
\be{ppm-ident}
\xi_2 \cdot \xi = \frac{|\xi_2|^2 - |\xi - \xi_2|^2}{2} + \frac{|\xi|^2}{2}
\end{equation}
for all $\xi_2, \xi$.

\begin{figure}[htbp] \centering
\ \psfig{figure=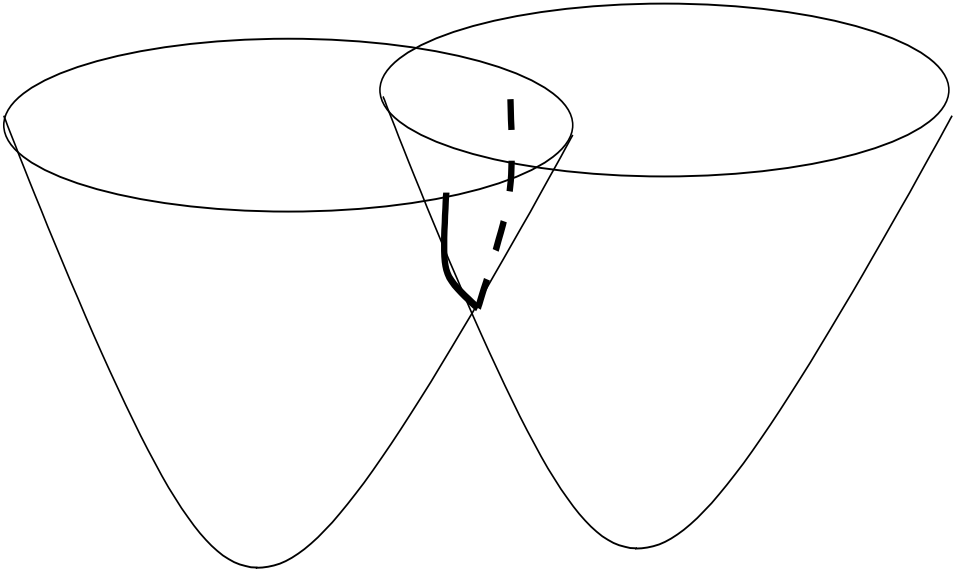,width=5in,height=3in}
\caption{The geometric interpretation of \eqref{ppm-ident}.  If one inverts an downward paraboloid and then translates it by $(\xi,\tau)$, the resulting intersection with the original upward paraboloid is a graph over a hyperplane orthogonal to $\xi_2$.}
\end{figure}

We first consider the case $N_1 \sim N_2 \gg N_3$.  In this case we have
$|\xi_3|^2 \ll |\xi_1|^2 + |\xi_2|^2$, so we may assume that $H \sim N_1^2$.
The case \eqref{t-dom} now follows from \eqref{t-dom-2}, so we may assume \eqref{h-dom}.  If $N_1^2 \sim H \sim L_3$ then we are essentially in the $(+++)$ case, since there is essentially no distinction between the constraints $|\tau - |\xi^2|| \sim N_1^2$ and $|\tau + |\xi^2|| \sim N_1^2$ when $|\xi| \sim N_3 \ll N_1$.  Thus by symmetry of the first and second variable it suffices to consider the case
$$ N_1^2 \sim H \sim L_1 \geq L_2, L_3.$$

By some permutation of Corollary \ref{char-2} and \eqref{ppm-ident} we thus have
$$ \eqref{h-gen} \lesssim L_{min}^{1/2}
| \{ \xi_2: |\xi_2 - \xi_2^0| \ll N_3;
\xi \cdot \xi_2 = \frac{\tau}{2} + \frac{|\xi|^2}{2} + O(L_{med}) \} |^{1/2}.$$
for some $\xi, \tau, \xi_1^0, \xi_2^0, \xi_3^0$ obeying \eqref{xio-cond} and $|\xi + \xi_1^0| \ll N_3$.
Since $|\xi| \sim N_1$, the measure of this set is $O(N_3^d \min(1, L_{med}/(N_1 N_3))$.  The claim \eqref{sppm-pp} follows.  (Alternatively, one can argue using Lemma \ref{transverse}).

Having disposed of the case $N_1 \sim N_2 \gg N_3$, it remains by symmetry to deal with the case
$$ N_1 \sim N_3 \gtrsim N_2.$$

First suppose that $H \sim N_1^2$, which forces $N_1 \sim N_2$.  Then we can repeat the $N_1 \sim N_2 \gg N_3$ arguments for this case.  Thus we may assume that $H \ll N_1^2$.

We first consider the case \eqref{t-dom}.  By \eqref{t-dom-1} we have
$$
\eqref{h-gen} \lesssim L_{min}^{1/2} 
\| \chi_{|h(\xi)| \lesssim H} \chi_{|\xi_1| \sim N_1} \chi_{|\xi_2| \sim N_2} \|_{[3,Z]}.
$$
The desired bound \eqref{sppm-est} then follows from \eqref{sppm-ident} and Proposition \ref{rotate-prop}.

It remains to consider the case \eqref{h-dom}.  The desired bound \eqref{sppm-est} has now simplified to
\be{sppm-simple}
\eqref{h-gen} \lesssim L_{min}^{1/2} N_2^{d/2}
(\frac{L_{med}}{N_1 N_2})^{1/2}
\min(1, \frac{H}{N_2^2})^{1/2 - \eps},
\end{equation}
though in the exceptional case \eqref{sppm-excep-case} our task is still to show \eqref{sppm-excep}.

We estimate \eqref{h-gen} by
\be{complex}
\eqref{h-gen} \lesssim \| \chi_{\angle(\xi_1,\xi_2) = \frac{\pi}{2} + O(\frac{H}{N_1 N_2})} \prod_{j=1}^3 \chi_{|\xi_j| \sim N_j} \chi_{|\lambda_j| \sim L_j} \|_{[3,\R^d \times \R]}.
\end{equation}

The main difficulty with \eqref{complex} is that there are two separate phenomena that need to be exploited to obtain a sharp bound.  The first is the transversality between the surfaces $\tau_j = h_j(\xi_j)$.  The other is the angular constraint \eqref{sppm-ident}.  To be able to exploit both simultaneously we again use Lemma \ref{split}.  The key is to cover $\R^d \times \R$ by boxes which are large enough to capture an optimum amount of transversality subject to the condition that the boxes stay small enough so as not to disturb the angular constraint.

We turn to the details.  Let $C$ be a large number ($C = 10^d$ will do) and define $A:= CN_1^2$, $B := C^{-1}\min(N_2, H/N_2)$.  Let $R$ be the box
$$ R := \{ (\xi,\tau): |\tau| \leq A; |\xi \cdot e_1| \leq H/N_1; |\xi \cdot e_j| \leq B \hbox{ for all } 2 \leq j \leq n\}.$$
Let $\Sigma$ be the canonical tiling lattice for $R$, so that $(R + (\eta,\lambda))_{(\eta,\lambda) \in \Sigma}$ is a box covering.  Let $M$ be the multiplier defined in \eqref{M-def}.  It is clear that 
$$ M((\eta_1,\lambda_1), (\eta_2, \lambda_2), (\eta_3, \lambda_3)) = 0$$
unless $|\eta_1|, |\eta_3| \sim N_1$, $|\eta_2| \sim N_2$, $|\eta_1 - N_1 e_1| \ll N_2$, $\angle(\eta_1, \eta_2) = \frac{\pi}{2} + O(\frac{H}{N_1 N_2})$, and $\lambda_1, \lambda_2, \lambda_3 = O(A)$.

\begin{figure}[htbp]\centering
\ \psfig{figure=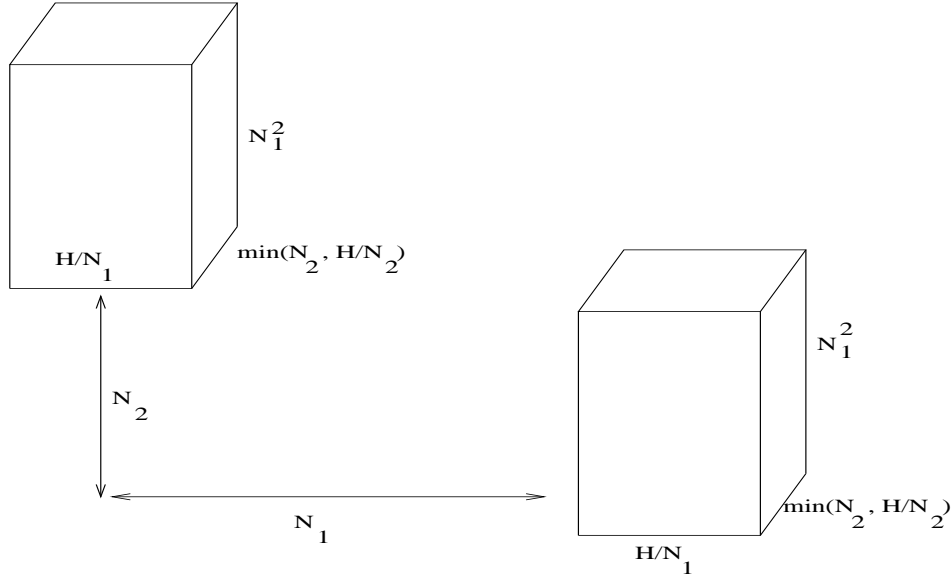,width=5in,height=3in}
\caption{The boxes used to subdivide the $(\xi_1,\tau_1)$ and $(\xi_2,\tau_2)$ variables into fine and coarse scales (ignoring factors of $C$).  Since $\tau_1, \tau_2 = O(N_1^2)$ the boxes do not extend much in the $\tau$ direction.  The boxes only interact if they subtend a spatial angle of $\pi/2 + O(H/N_1 N_2)$, which we can treat by Proposition \ref{rotate-prop}.  To treat the interaction of a single pair of boxes we can use the type of techniques used previously (basically exploiting the transversality between the two sections of the paraboloid).}
\end{figure}

We now claim that
\be{M-brane}
|M((\eta_1,\lambda_1), (\eta_2, \lambda_2), (\eta_3, \lambda_3))| \lesssim
L_{min}^{1/2} B^{(d-1)/2} (L_{med}/N_1)^{1/2}
\end{equation}
except in the exceptional case \eqref{sppm-excep-case}, in which case one must place an additional factor of $(\frac{H}{N_2^2 + L_{med}})^{1/2}$ on the right-hand side.

Assuming \eqref{M-brane} for the moment, we see from Lemma \ref{split} that
\bas
\eqref{h-gen} \lesssim &(\frac{H}{N_1} B^{d-1} A)^{-1/2}
L_{min}^{1/2} B^{(d-1)/2} (\frac{L_{med}}{N_1})^{1/2}\\
&\| \chi_{|\xi_1| \sim N_1, |\xi_2| \sim N_2} \chi_{\angle(\xi_1,\xi_2) \lesssim \frac{H}{N_1 N_2}} \chi_{|\tau_1|, |\tau_2|, |\tau_3| \lesssim A}
\|_{[3; \R^d \times \R]}.
\end{align*}
From Lemma \ref{box} (for instance) we have
$$ \|\chi_{|\tau_1|, |\tau_2|, |\tau_3| \lesssim A}
\|_{[3; \R]} \lesssim A^{1/2},$$
so by the tensor product lemma (Lemma \ref{tensor}) the above simplifies to
$$
\eqref{h-gen} \lesssim 
L_{min}^{1/2} L_{med}^{1/2} H^{-1/2}
\| \chi_{|\xi_1| \sim N_1} \chi_{|\xi_2| \sim N_2} \chi_{\angle(\xi_1,\xi_2) \lesssim H/N_1 N_2} \|_{[3; \R^d]}.$$
Applying Proposition \ref{rotate-prop} we thus obtain
$$
\eqref{h-gen} \lesssim
L_{min}^{1/2} L_{med}^{1/2} H^{-1/2}
N_2^{d/2} (\frac{H}{N_1 N_2})^{1/2}
\min(1, \frac{H}{N_2^2})^{\frac{1}{2}-\eps}$$
which is \eqref{sppm-simple} as desired.  In the exceptional case \eqref{sppm-excep-case} we repeat the above argument but acquire an extra factor of $(\frac{H}{N_2^2 + L_{med}})^{1/2}$, which ultimately yields \eqref{sppm-excep} instead of \eqref{sppm-simple}.

It remains to prove \eqref{M-brane}.  Fix $\eta_1, \eta_2, \eta_3, \sigma_1, \sigma_2, \sigma_3$ with $\eta_1 + \eta_2 + \eta_3 = 0$, $\sigma_1 + \sigma_2 + \sigma_3 = 0$, and $\angle(\eta_1, \eta_2) \lesssim H/N_1 N_2$.  It suffices to show that
$$
\|
\prod_{j=1}^3 \chi_{R + (\eta_j,\sigma_j)}(\xi_j,\tau_j) \chi_{|\xi_j| \sim N_j} \chi_{|\lambda_j| \lesssim L_j} \|_{[3,\R^d \times \R]}
\lesssim 
L_{min}^{1/2} B^{(d-1)/2} (L_{med}/N_1)^{1/2},$$
with an additional factor of $(\frac{H}{N_2^2 + L_{med}})^{1/2}$ in the right-hand side in the exceptional case \eqref{sppm-excep}.
We may assume that $|\eta_j| \sim N_j$, since the above expression vanishes otherwise.

We shall use Lemma \ref{char-1} and elementary geometry\footnote{Alternatively, one can apply Lemma \ref{transverse}.  Another method is to use Lemma \ref{tensor} to reduce to the $d=1$ case.}.  Let $\{j_1, j_2, j_3\} = \{1,2,3\}$ be such that $L_{j_1} \geq L_{j_2} \geq L_{j_3}$.  Discarding the restrictions on $\xi_{j_1}$ and $\tau_{j_1} - h_{j_1}(\xi_{j_1})$, and using Lemma \ref{char-1} we reduce to showing
\be{3-way}
\begin{split}
| \{ \xi_1: &\xi_1 = \eta_{j_2} + O(B); \quad (\xi_1 - \eta_{j_2}) \cdot e_1 = O(\frac{H}{N_1}); \\
& h_{j_2}(\xi_1) + h_{j_3}(\xi - \xi_1) = \tau + O(L_{med}) \}|
\lesssim B^{d-1} \frac{L_{med}}{N_1}
\end{split}
\end{equation}
for all $\xi, \tau$, with an additional factor of $H/(N_2^2 + L_{med})$ when \eqref{sppm-excep} holds. 

Fix $\xi,\tau$; we may assume that $(\xi,\tau) \in (-\eta_{j_1},-\sigma_{j_1}) + CR$ since the set vanishes otherwise.  In particular we have $|\xi| \sim N_{j_1}$. 

We now split into three cases depending on the value of $\{j_2,j_3\}$.  (The estimate \eqref{3-way} is symmetric with respect to interchanging $j_2$ and $j_3$).

First suppose that $\{j_2,j_3\} = \{1,2\}$.  Using \eqref{ppp-ident} we can estimate the left-hand side of \eqref{3-way} by
$$
| \{ \xi_1: \xi_1 = \eta_1 + O(B);
|\xi_1 - \frac{\xi}{2}|^2 = \frac{\tau}{2} - \frac{|\xi|^2}{4} + O(L_{med})
\}|.$$
To finish the proof of \eqref{3-way} in this case it thus suffices to show that $|\xi_1 - \frac{\xi}{2}| \sim N_1$.  Since $|\xi| \sim N_1$ and $|\xi_1| \sim N_2$, we are done unless $N_1 \sim N_2$.  But in this case we have
$$ \angle (\eta_1,\eta_2) = \frac{\pi}{2} + O(H/N_1^2) \sim 1$$
since we have assumed $H \ll N_1^2$.  Since $\xi_1$, $\xi$ are within $\ll N_2$ of $\eta_2$ and $-\eta_1-\eta_2$ respectively, the claim \eqref{3-way} follows.

Now suppose that $\{j_2,j_3\} = \{2,3\}$.  Using \eqref{ppm-ident} we can estimate the left-hand side of \eqref{3-way} by
$$
| \{ \xi_1: \xi_1 = \eta_2 + O(B);
\xi_1 \cdot \xi = \frac{\tau}{2} + \frac{|\xi|^2}{2} + O(L_{med})
\}|.$$
Since $|\xi| \sim N_{j_1} = N_1$, we thus see that this set is contained in a $O(L_{med}/N_1)$ neighbourhood of a hyperplane, and also in a ball of radius $O(B)$.  The claim \eqref{3-way} follows.

Finally, suppose that $\{j_2,j_3\} = \{1,3\}$.  Using \eqref{ppm-ident} as before, the left-hand side of \eqref{3-way} becomes
$$
| \{ \xi_1: \xi_1 = \eta_1 + O(B); (\xi_1 - \eta_{j_2}) \cdot e_1 = O(H/N_1);
\xi_1 \cdot \xi = \frac{\tau}{2} + \frac{|\xi|^2}{2} + O(L_{med})
\}|.$$
First suppose that $H \sim N_1 N_2$, so in particular $B \sim N_2$.  Since $|\xi| \sim N_{j_1} = N_2$, the above set is contained in a $O(L_{med}/N_2)$ neighbourhood of a hyperplane, and also in a ball of radius $O(B)$.  One can then bound the measure of the above set as
$$ O(\min(B^d, B^{d-1} \frac{L_{med}}{N_2})) = O( B^{d-1} \frac{L_{med}}{N_1} \frac{N_1 N_2}{N_2^2 + L_{med}})$$
which is \eqref{3-way} with the additional factor of $\frac{H}{N_2^2 + L_{med}}$.

Now suppose $H \ll N_1 N_2$, so that $\angle(\eta_1, \eta_2) \sim 1$.  From this and the estimates 
$$|\xi + \eta_1 + \eta_3| \ll N_2; |\eta_1 + \eta_2 + \eta_3| \lesssim B \ll N_2; |\eta_j| \sim N_j; |\eta_1 - N_1 e_1| \ll N_2$$
we see that $|\xi| \sim N_2$ and $\angle(\xi,e_1) \sim 1$.  From this we see that the size of the above set is at most $O(B^{d-2} H/N_1 \min(B,L_{med}/N_2))$.
When $H \lesssim N_2^2$, then $B \sim H/N_2$ and $L_{med}/N_2 \lesssim H/N_2 = B$, and the above estimate simplifies to \eqref{3-way} as desired.  When $H \gg N_2^2$ then $B \sim N_2$ and the above estimate simplifies to
$$ O(N_2^{d-2} \frac{H}{N_1} \min(N_2, \frac{L_{med}}{N_2}))
= O( N_2^{d-1} \frac{L_{med}}{N_1} \frac{H}{N_2^2 + L_{med}})$$
which is \eqref{3-way} with the additional factor of $\frac{H}{N_2^2 + L_{med}}$.
\end{proof} 

As a sample application of the above Lemma we give

\begin{proposition}\label{new-schro}
For all $u,v$ in $\R^3 \times \R$, we have
\be{ns-est}
\| u \overline{v} \|_{X^{s,-1/2 + \eps}_{\tau = |\xi|^2}(\R^3 \times \R)}
\lesssim 
\| u \|_{X^{s,1/2 - \eps}_{\tau = |\xi|^2}(\R^3 \times \R)}
\| v \|_{X^{s,1/2 - \eps}_{\tau = |\xi|^2}(\R^3 \times \R)}
\end{equation}
whenever $\eps > 0$ and $0 \geq s > -1/4 + C\eps$, with the implicit constant depending on $s$ and $\eps$.
\end{proposition}

Of course, this $d=3$ result implies the same result for $d=1,2$ by the method of descent (or the Comparison principle and Lemma \ref{tensor}).  For $d=2,3$ the exponent $-1/4$ is sharp up to epsilons; see\footnote{There is a misprint in the proof of that theorem; the quantities $|R^+|, |R^-|, |R^0|$ should be $\sim 1$ rather than $\sim N$.} \cite{staff:quadratic}, Theorem 2.2. 
From this Proposition we have local well-posedness of
$$ u_t = i \Delta u + C |u|^2$$
in $H^s(\R^3)$ for all $s > -1/4$ and any complex $C$.  We remark that the best result for this equation one can obtain via Strichartz estimates is $s \geq 0$, see \cite{cwI}.

The bilinear expressions $uv$ and $\overline{u} \overline{v}$ are better behaved than $u \overline v$; for instance, one can use the above Lemmata to derive the analogue of \eqref{ns-est} for those forms when $s > -1/2 + C\eps$.   In two dimensions the correct exponent is $s > -3/4 + C \eps$, see \cite{staff:quadratic}, \cite{cdks}.  For $d \geq 4$ one can go down to scaling $s > \frac{d-4}{2} + C\eps$, but this is inferior to Strichartz techniques, which can give well-posedness for $s \geq \frac{d-4}{2}$ for any type of quadratic nonlinearity (see \cite{cwI}).

\begin{proof}
By duality and permutation of indices the estimate is equivalent to
$$
\|
\frac{ \langle \xi_1 \rangle^{-s} \langle \xi_3 \rangle^{-s} }{\langle \xi_2 \rangle^{-s}
\langle \tau_1 - |\xi_1|^2 \rangle^{1/2-\eps}
\langle \tau_2 - |\xi_2|^2 \rangle^{1/2-\eps}
\langle \tau_3 + |\xi_3|^2 \rangle^{1/2-\eps}
} \|_{[3;\R^2 \times \R]} \lesssim 1.$$
Thus we are in the $(++-)$ case \eqref{sppm}.  By \eqref{m-bound-1-tri}, \eqref{m-bound-2-tri} it suffices to show that
\be{sppm-1-tri}
\sum_{N_{max} \sim N_{med} \sim N} 
\sum_{L_1,L_2,L_3 \gtrsim 1}
\frac{\langle N_1\rangle^{-s} \langle N_3 \rangle ^{-s} }
{ \langle N_2\rangle^{-s} (L_1 L_2 L_3)^{1/2-\eps} }
\| X_{N_1,N_2,N_3;L_{max};L_1,L_2,L_3} \|_{[3, Z \times \R]}
\lesssim 1
\end{equation}
and
\be{sppm-2-tri}
\sum_{ N_{max} \sim N_{med} \sim N} 
\sum_{ L_{max} \sim L_{med}} \sum_{H \ll L_{max}}
\frac{\langle N_1\rangle^{-s} \langle N_3 \rangle ^{-s} }
{ \langle N_2\rangle^{-s} (L_1 L_2 L_3)^{1/2-\eps} }
\| X_{N_1,N_2,N_3;H;L_1,L_2,L_3} \|_{[3, Z \times \R]}
\lesssim 1
\end{equation}
for all $N \gtrsim 1$.  

Since $s \leq 0$, we can estimate
$$ \frac{\langle N_1\rangle ^{-s} \langle N_3 \rangle ^{-s} }
{ \langle N_2\rangle^{-s} (L_1 L_2 L_3)^{1/2-\eps} } \lesssim \frac{ N^{-2s} }{\langle N_{min} \rangle^{-s} L_{min}^{1/2} L_{med}^{1/2} L_{max}^{1/2-3\eps}
}.$$
Also, from the various cases of Proposition \ref{sppm-core} we have
$$ \| X_{N_1,N_2,N_3;H;L_1,L_2,L_3} \|_{[3, Z \times \R]}
\lesssim L_{min}^{1/2} L_{med}^{1/2} N^{-1/2} N_{min}
(\frac{N_{min}^2}{H})^{\eps}
(\frac{H}{N_{min}^2})^{1/2}
$$
whenever $N_{max} \sim N$.  Combining these we obtain
$$ \frac{\langle N_1\rangle^{-s} \langle N_3 \rangle ^{-s} }
{ \langle N_2\rangle^{-s} (L_1 L_2 L_3)^{1/2-\eps} } 
\| X_{N_1,N_2,N_3;H;L_1,L_2,L_3} \|_{[3, Z \times \R]}
\lesssim \frac{ N^{-2s-1/2} N_{min}^{2\eps} H^{1/2 - \eps} }{ \langle N_{min}\rangle^{-s} L_{max}^{1/2 - 3\eps} }.$$
The claims \eqref{sppm-1-tri}, \eqref{sppm-2-tri} are now easily verified using the fact that $H \lesssim N^2$.  Note that the $L$, $H$, $N$ summations may give powers of $N^{\eps}$, but this acceptable because of $N^{-2s-1/2}$ factor and the assumption $s > -1/4 + C\eps$.
\end{proof}

\end{document}